\documentclass[EJP]{ejpecp} 

\usepackage{amsrefs}
\usepackage{enumerate}  
\usepackage{bbm}
\usepackage{mathrsfs}
\newcommand{\E}{\mathbbm{E}}

\newcommand{\N}{\mathbbm{N}}

\newcommand{\D}{\mathcal{D}}
\DeclareMathOperator{\Tr}{Tr}
\allowdisplaybreaks[1]

\SHORTTITLE{Multivariate functional Stein's method of exchangeable pairs}

\TITLE{Stein's method of exchangeable pairs in multivariate functional approximations} 

\AUTHORS{%
  Christian~D\"obler\footnote{Heinrich-Heine-Universit\"at D\"usseldorf, Germany.
    \EMAIL{christian.doebler@hhu.de}}
  \and 
  Miko{\l}aj J.~Kasprzak\footnote{University of Luxembourg, Luxembourg. \EMAIL{mikolaj.kasprzak@uni.lu}}\
  \thanks{supported by the \textbf{FNR grant FoRGES (R-AGR- 3376-10)} at Luxembourg University}}
  

\AMSSUBJ{60B10; 60F17} 
\AMSSUBJSECONDARY{0B12; 60J65; 60E05; 60E15} 

\SUBMITTED{June 1, 2020} 
\ACCEPTED{January 30, 2021} 

\VOLUME{26}
\YEAR{2021}
\PAPERNUM{28}
\DOI{10.1214/21-EJP587}

\ABSTRACT{In this paper we develop a framework for multivariate functional approximation by a suitable Gaussian process via an exchangeable pairs coupling that satisfies a suitable approximate linear regression 
property, thereby building on work by Barbour (1990)  and Kasprzak (2020). We demonstrate the applicability of our results by applying them to joint subgraph counts in an Erd\H{o}s-Renyi random graph model on the one hand and to vectors of weighted, degenerate $U$-processes on the other hand. As a concrete instance of the latter class of examples, we provide a bound for the functional approximation of a vector of 
success runs of different lengths by a suitable Gaussian process which, even in the situation of just  a single run, would be outside the scope of the existing theory.}

\begin{document}



\section{Introduction}
In his seminal paper \cite{stein}, Charles Stein introduced a method for proving normal approximations and obtained a bound on the speed of convergence to the standard normal distribution. Later, Barbour \cite{barbour} and G\"otze \cite{gotze} developed the so-called generator approach to finding Stein's equation, which made it possible to study approximations by many other probability laws. As a result, in \cite{diffusion}, the method was adapted to approximations by the (infinite-dimensional) Wiener measure.  

Moreover, the exchangeable-pair approach, first developed by Stein in his monograph \cite{stein1} in the context of univariate normal approximations, has been at the heart of many results proved using Stein's method. It was extended by \cite{RiRo97} and used in the context of non-normal approximations in \cite{chatterjee, roellin1, chatterjee1,  dobler, DP18b}. The publication of \cite{meckes, reinert_roellin, meckes09} brought a breakthrough in the understanding of the exchangeable-pair approach and made it available for applications to a wide array of multivariate normal approximation problems. The very recent paper \cite{Doe20} developed a functional analytic approach that provides a substantial extension of the method of exchangeable pairs and, in particular, makes it possible to dispense with the linear regression property in finite-dimensional settings. In \cite{kasprzak3} the method was applied to the study of functional limit results and approximations by univariate Gaussian processes, using the setup of \cite{stein1, RiRo97} and \cite{diffusion}.

In this paper we combine the functional approximation of \cite{diffusion} and the multivariate exchangeable-pair method of \cite{reinert_roellin, meckes09}. We obtain an abstract approximation theorem, which is applied in the context of weighted degenerate U-statistics, a particularly interesting example of which are homogeneous sums. The strength of the abstract approximation result is also presented in a random-graph-theoretic application.
\subsection{Motivation}
We are motivated by examples of multivariate quantities whose distance from the normal distribution can be established using Stein's method of exchangeable pairs, and whose functional equivalents have not been studied yet. Functional limit results play an important role in applied fields. Scaling limits of discrete processes can be studied using stochastic analysis and are often more robust to changes in the local details than the discrete processes themselves. That is why researchers often choose to describe discrete phenomena with continuous models. The error they make by doing this is measured by rates of convergence in functional limit results. The current paper contributes to solving the problem of bounding those rates.

The two main applications motivating the paper and considered therein are a continuous Gaussian-process approximation of a rescaled weighted U-statistic and the study of an Erd\H{o}s-Renyi random graph process.  U-statistics are central objects in the field of mathematical statistics. Due to their appealing properties, they have found numerous applications to estimation, statistical testing and other problems. They appear in decompositions of more general statistics into sums of terms of a simpler form (see, e.g. \cite[Chapter 6]{serfling} or \cite{rubin_vitale} and \cite{vitale}) and play an important role in the study of random fields (see, e.g. \cite[Chapter 4]{christofides}). Moreover, functional limit theorems for rescaled U-statistics have found applications in the field of changepoint analysis (see e.g. \cite{CH88, Ferger94, GH95,  CH_book, Ferger01, Gombay2004, HR_survey, RW19}), where it is particularly useful to know the functional limits of the related test statistics. On the other hand, the Erd\H{o}s-Renyi random graph model has found numerous applications in various fields (see \cite{models}), including epidemic modelling \cite{barbour_mollison} and modelling of evolutionary conflicts \cite{cannings}.

The first application discussed in the paper deals with the approximation of so-called weighted $U$-processes, i.e. process analogues of the class of weighted $U$-statistics. This class of processes is very wide, containing the so-called homogeneous sum processes as well as symmetric, degenerate (complete or incomplete) $U$-processes.
We derive a general result and successfully apply it to the case of homogeneous sum processes in Subsection \ref{homsums}. As a concrete example, in Subsection \ref{runs}, we provide a bound for a Gaussian approximation of a process that is defined as a vector of success runs of different lengths.
For functional limit theorems involving the class of symmetric, degenerate $U$-processes, we refer the reader to the recent paper \cite{DKP}.
Moreover, we remark that, even in the univariate case of weighted $U$-statistics, the literature about limit theory for these random quantities is quite restricted. Indeed, apart from the abundance of references on limit theorems for homogeneous sums, the majority of articles focus on the limiting behavior of so-called reduced or incomplete $U$-statistics, i.e.
weighted $U$-statistics whose weights only assume the values $0$ and $1$ (see e.g. \cite{Blom, Ja, BrKil}). Limit theorems for general weighted $U$-statistics can be found in references \cite{RiRo97, OnRe, Major}. We stress, however, that the last two references focus on non-normal limiting distributions and that, in the degenerate case, \cite{RiRo97} only considers kernels of order $2$. Moreover, the literature about functional central limit theorems (FCLTs) for weighted $U$-statistics is even scarcer. Indeed, only for homogeneous sum processes \cite{Mik, Basa} have we been able to find comparable results in the literature. We defer a discussion and comparison with our findings to Subsection \ref{homsums}.

The second example comes originally from \cite{Janson1991} and was studied using exchangeable pairs in a finite-dimensional context in \cite{reinert_roellin1}. We look at a (dynamic) Erd\H{o}s-Renyi random graph with $\lfloor nt\rfloor$ vertices, where $t$ denotes the time, and study the distance from the asymptotic distribution of the joint law of the number of edges and the number of two-stars.
 Our approach can, however, be also extended to cover the number of triangles. Those statistics are often used when approximating the clustering coefficient of a network and applied in conditionally uniform graph tests. 

\subsection{Contribution of the paper}
The main achievements of the paper are the following:
\begin{enumerate}
\item An abstract approximation theorem (Theorem \ref{theorem1}), bounding the distance between a stochastic process $\mathbf{Y}_n$ valued in $\mathbbm{R}^d$, for a fixed positive integer $d$, and a Gaussian mixture process. The estimate is derived under the assumption that that the process $\mathbf{Y}_n$ satisfies the linear regression condition 
\begin{equation}\label{cond}
Df(\mathbf{Y}_n)[\mathbf{Y}_n]=2\mathbbm{E}\left\lbrace \left.\vphantom{\sum}Df(\mathbf{Y}_n)\left[(\mathbf{Y}_n-\mathbf{Y}_n')\Lambda_n\right]\,\right|\,\mathbf{Y}_n\right\rbrace+R_f,
\end{equation}
for all $f:D\left([0,1],\mathbbm{R}^d\right)\to\mathbbm{R}$ in a certain class of test functions, a random process $\mathbf{Y}_n'$ such that $(\mathbf{Y}_n,\mathbf{Y}_n')$ is an exchangeable pair, some $\Lambda_n\in\mathbbm{R}^{d\times d}$ and some random variable $R_f=R_f(\mathbf{Y}_n)$. In \eqref{cond} (and in the entire paper) $Df$ denotes the Fr\'echet derivative of $f$. The class of test functions, with respect to which the bound in Theorem \ref{theorem1} is obtained, is so rich that the bound approaching zero fast enough implies weak convergence of the law of $\mathbf{Y}_n$ in the Skorokhod and uniform topologies on the Skorokhod space. The exact conditions under which this happens are stated in Proposition \ref{prop_m}.
\item A novel framework for continuous Gaussian process approximations of vectors of weighted, degenerate $U$-processes, presented in Section \ref{weighted}. Apart from proving a general result about those, we show how it may be applied in examples involving non-degenerate $U$-processes. In order to study such examples using our theory, one may decompose the given $U$-process into the vector of its degenerate Hoeffding components and prove a multivariate Gaussian limit theorem for this vector. Then, by applying a linear functional, one obtains a Gaussian limit for the original process. This strategy, in a quantified fashion, is exemplified by the application to the $r$-runs process, discussed in Subsection \ref{runs}. We stress that, even in the case of just one $r$-run process, the results about univariate functional approximations via exchangeable pairs from \cite{kasprzak3} would not be sufficient to obtain a Gaussian approximation. Thus, in this example, the multidimensionality of our approach proves to be absolutely vital. Moreover, both the kernels and the coefficients of the weighted $U$-processes we study in our general result may (and will in most cases) depend on the sample size $n$, hence yielding Gaussian limits even in degenerate situations.
At the same time, our methods are flexible enough in order to yield bounds for the classical results on asymptotic Gaussianity, in non-degenerate situations, when the kernels are fixed.
\item A novel quantitative functional limit theorem for the edge counts and the number of two-stars in an Erd\H{o}s-Renyi random graph $G(n,p)$ on $n$ vertices with fixed edge probability $p$. Letting $I_{i,j}$, for $i,j=1,\cdots,n$ be the indicator that edge $(i,j)$ is present in the graph, we consider the following statistics:
\[\mathbf{T}_n(t)=\frac{\lfloor nt\rfloor -2}{2n^2}\sum_{i,j=1}^{\lfloor nt\rfloor}I_{i,j},\quad\mathbf{V}_n(t)=\frac{1}{6n^2}\underset{i,j,k\text{ distinct}}{\sum_{1\leq i,j,k\leq \lfloor nt\rfloor}}I_{ij}I_{jk},\qquad t\in[0,1],\]
corresponding to the number of edges and the number of two-stars, respectively.
Theorem \ref{theorem_pre_limiting} provides a bound on the distance between the law of the process
\begin{equation}\label{graph_ex}
t\mapsto \left(\mathbf{T}_n(t)-\mathbbm{E}\mathbf{T}_n(t),\mathbf{V}_n(t)-\mathbbm{E}\mathbf{V}_n(t)\right)\quad t\in[0,1]
\end{equation}
and the law of a piecewise constant Gaussian process. Theorem \ref{theorem_continuous} estimates the distance between the law of (\ref{graph_ex}) and that of a continuous Gaussian process. These results extend the result of \cite{kasprzak3} bounding the distance between the distribution of the edge counts and a univariate Gaussian process. As a corollary to our results, we immediately obtain weak convergence of the law of (\ref{graph_ex}) in the Skorokhod and uniform topologies on the Skorokhod space to that of the continuous Gaussian process.
\end{enumerate}
\subsection{Stein's method in its generality}
Stein's method in its generality is a powerful technique used to obtain bounds on quantities of the form $|\mathbbm{E}_{\nu}h-\mathbbm{E}_{\mu}h|$, where $\mu$ is the target (known) distribution, $\nu$ is an approximating measure and $h$ is a real-valued test function chosen from a suitable class $\mathcal{H}$. The method is composed out of three main steps. First, one needs to find an operator $\mathcal{A}$ acting on a class of real-valued functions, such that \[\left(\forall f\in\text{Domain}(\mathcal{A})\quad\mathbbm{E}_{\pi}\mathcal{A}f=0\right)\quad \Longleftrightarrow \quad\pi=\mu.\]
Second, for a given function $h\in\mathcal{H}$, one solves the following Stein equation:
\[\mathcal{A}f=h-\mathbbm{E}_{\mu}h.\]
Finally, for $f=f_h$ solving the Stein equation, the following quantity:
\begin{equation}\label{quantity}
|\mathbbm{E}_{\nu}\mathcal{A}f_h|
\end{equation}
needs to be bounded. This is achieved using various mathematical tools (Taylor's expansions, Malliavin calculus, as described in \cite{nourdin}, coupling methods and others), applied in conjunction with smoothness properties of $f_h$.
For an accessible account of the method we recommend the surveys \cite{reinert} and \cite{ross} as well as the books \cite{janson} and \cite{normal_approx}, which treat the cases of Poisson and normal approximation, respectively, in detail. A database of information and publications connected to Stein's method can also be found in \cite{swan}.
\subsection{Stein's method of exchangeable pairs}
The exchangeable-pair approach to Stein's method was first developed in \cite{stein1}.  Therein, the author considered the setup in which, for a random variable $W$, one can construct another random variable $W'$ such that $(W,W')$ is an exchangeable pair and the following linear regression condition is satisfied
\begin{equation}\label{regression_condition}
\mathbbm{E}\left[W'-W|W\right]=-\lambda W
\end{equation}
for some $\lambda>0$. It follows from this assumption that
\begin{align*}
0=&\mathbbm{E}\left[(f(W)+f(W'))(W-W')\right]
=\mathbbm{E}\left[(f(W')-f(W))(W-W')\right]+2\lambda\mathbbm{E}[Wf(W)]
\end{align*}
and so
\[\mathbbm{E}[Wf(W)]=\frac{1}{2\lambda}\mathbbm{E}\left[(f(W)-f(W'))(W-W')\right].\]
Therefore, using Taylor's theorem, it can be proved that
\begin{align*}
&\left|\mathbbm{E}[f'(W)]-\mathbbm{E}[Wf(W)]\right|
\leq\frac{\|f'\|_{\infty}}{2\lambda}\sqrt{\text{Var}\left[\mathbbm{E}\left[(W-W')^2|W\right]\right]}+\frac{\|f''\|_{\infty}}{2\lambda}\mathbbm{E}|W-W'|^3,
\end{align*}
which provides a bound on the quantity \eqref{quantity} for $\nu=\mathcal{L}(W)$ and $\mathcal{A}$ being the canonical Stein operator corresponding to the standard normal law.

 A multivariate version of the method was first described in \cite{meckes} and then in \cite{reinert_roellin}. In  \cite{reinert_roellin}, for an exchangeable pair of $d$-dimensional vectors $(W,W')$, the following condition is used:
\begin{equation}\label{lambda_matrix}
\mathbbm{E}[W'-W|W]=-\Lambda W+R
\end{equation}
for some invertible matrix $\Lambda$ and a remainder term $R$.
The approach of \cite{reinert_roellin} was further reinterpreted and combined with the approach of \cite{meckes} in \cite{meckes09}. Extending this multivariate version of the exchangeable-pair method to multivariate functional approximations, with the linear regression condition taking form similar to \eqref{lambda_matrix}, is the subject of the current paper.

\subsection{Functional Stein's method}
Approximations by laws of stochastic processes using Stein's method have been studied in \cite{diffusion, functional_combinatorial, shih, Coutin, decreusefond_higher} and recently in  \cite{kasprzak1, kasprzak2, kasprzak3, decreusefond2, decreusefond_rough, campese}. These references can be divided into three groups. 

The ones belonging to the first group, containing \cite{diffusion, functional_combinatorial, kasprzak1, kasprzak2, kasprzak3}, all use, adapt and extend the setup of \cite{diffusion}. Therein, the author studied the rate of convergence in the celebrated functional central limit theorem, also called Donsker's theorem. Barbour considered test functions $g$ acting on the Skorokhod space $D \left( [0,1], \mathbbm{{R}} \right)$ of c\`adl\`ag real-valued maps on $[0,1]$, such that $g$ takes values in the reals, does not grow faster than a cubic, is twice Fr\'echet differentiable and its second derivative is Lipschitz. For each function $g$ belonging to this class he provided a bound on the absolute difference between the expectation of $g$ with respect to the law of a rescaled random walk and the expectation of $g$ with respect to the Wiener measure. Crucially, he also proved that this class of functions $g$ is so rich that his bounds imply weak convergence with respect to the Skorokhod topology of the considered rescaled random walk to Brownian Motion. This last property is vital for most applications of the limit theory for stochastic processes and may even be the main reason for the outstanding 
popularity of the Skorohod topology. Indeed, by means of the continuous mapping theorem, limit theorems for many natural, non-linear functionals such as the supremum over time, immediately follow from 
a weak limit theorem in the Skorokhod topology.

On the other hand, the results of the second group of references, containing \cite{Coutin, decreusefond_higher, decreusefond2, decreusefond_rough}, develop Stein's theory on a Hilbert space using a Besov-type topology. The bounds obtained therein, however, do not imply weak convergence in the Skorokhod topology. Therefore, the continuous mapping theorem does not apply in their setting. For instance, as opposed to the results of the first group of references, one cannot study convergence of the supremum of a process using the analysis of the second group of papers.

Finally, \cite{shih} develops approximations by abstract Wiener measures on a real separable Banach space and \cite{campese} proves bounds on measure-determining distances from Gaussian random variables valued in Hilbert spaces. As for the second group, despite the elegant abstract theory used and developed in these references, the results do not imply convergence in the Skorokhod topology on $D[0,1]$. 

In the current paper we shall follow the setup of the first group of references. We consider it more flexible than the one of the second group and more suited for applications to processes belonging to the widely-used (non-separable) Skorokhod space than the ones of the third group.

In the context of these three groups of references and the present paper, we also mention the recent paper \cite{DKP} which, although not relying on functional approximation by Stein's method, provides functional limit theorems for the class of (degenerate and non-degenerate) symmetric $U$-processes with a kernel that may depend on the sample size $n$. Since it implicitly relies on a 
multivariate Gaussian limit theorem derived by Stein's method from \cite{DP19}, it is also naturally related to Stein's method.

Moreover, since one main class of applications in the present paper involves weighted $U$-processes, it is worthwhile to compare our results and their applicability to those 
of \cite{DKP}. Firstly, as mentioned above, the paper \cite{DKP} focuses on Gaussian limit theorems for symmetric $U$-processes, which constitute a narrower class than the weighted $U$-processes considered in the present work. Moreover, thanks to the finite-dimensional convergence results from \cite{DP19}, the conditions for convergence from \cite{DKP} are phrased in term of $L^2$-norms of \textit{contraction kernels} and, as such, can be considered as fourth moment conditions. In contrast, as can be seen from the bounds and proofs of Section \ref{weighted}, the bounds and conditions in the present paper involve third moment quantities. This distinction is also clearly reflected in the respective applicability of the results proved in the present paper and those from \cite{DKP}. Indeed, whereas the symmetric $U$-processes considered in \cite{DKP} possess a global dependency structure, the results in Section \ref{weighted} are most useful whenever the dependence of the weighted $U$-process is local in the sense that the involved array of weighting coefficients $(a_J)_J$ is sparse in some sense. The runs example in Subsection \ref{runs} provides an instructive showcase for this observation. Moreover, the methods used in the proofs of the main results necessitate that the quantities in the bounds involve the absolute values of both the kernels and the coefficients. Hence, no cancellation effect, typically occuring under fourth moment conditions, may be relied on in this case. We therefore consider our theorems as rather complementary to the ones in \cite{DKP}. 

\subsection{Structure of the paper}
Section \ref{section_notation} includes some introductory remarks about notation and the spaces of test functions with respect to which bounds on distances between probability laws in this paper will be derived. Section \ref{section_setup_stein} gives a general form of the pre-limiting process to which all the processes of interest will be compared using Stein's method. It also presents the corresponding Stein equation, its solution and the smoothness properties of the solution. Section \ref{section_abstract} contains the main abstract result of this paper providing a bound on the distance between a process valued in the Skorokhod space $D([0,1],\mathbbm{R}^d)$ and the pre-limiting process described in the previous section. Section \ref{weighted} discusses the application of the abstract theorem to weighted, degenerate U-statistics and presents a bound on their distance from a continuous Gaussian process. It furthermore explains how the bound simplifies in the context of homogeneous sums and applies it to the example of $r$-runs on the line. Section \ref{section6} discusses the example concerning an Erd\H{o}s-Renyi random graph process and the bound on the distance between the number of its edges and two-stars and a continuous Gaussian process. Technical details of some of the proofs in this paper are postponed to Section \ref{appendix}.

\section{Notation and spaces $M$ and $M^0$}\label{section_notation}
The following notation, similar to the one of \cite{diffusion} and \cite{kasprzak2}, is used throughout the paper.  For a fixed positive integer $d$, let $D([0,1],\mathbbm{R}^d)$ be the Skorokhod space of  c\`adl\`ag $\mathbbm{R}^d$-valued functions on $[0,1]$. For $i=1,\cdots,d$, by $e_i$ we denote the $i$th unit vector of the canonical basis of $\mathbbm{R}^d$. The $i$th component of any $x\in\mathbbm{R}^d$ will be denoted by $x^{(i)}$, so that $x=\left(x^{(1)},\cdots,x^{(d)}\right)$.
For a function $w$ defined on  $[0,1]$ and taking values in a Euclidean space, we will also write \[\|w\|=\sup_{t\in[0,1]}|w(t)|,\] where $|\cdot|$ denotes the Euclidean norm. Moreover, the notation $\mathbbm{E}^W[\,\cdot\,]$ will be used to represent $\mathbbm{E}[\,\cdot\,|W]$.

Furthermore, we define
\[\|f\|_L:=\sup_{w\in D([0,1],\mathbbm{R}^d)}\frac{|f(w)|}{1+\|w\|^3}\text{,}\]
and let $L$ be the Banach space of continuous functions $f:D([0,1],\mathbbm{R}^d)\to\mathbbm{R}$ such that $\|f\|_L<\infty$.
By $D^kf$ we will always mean the $k$-th Fr\'echet derivative of $f$. The norm $\|\cdot\|$ of a $k$-linear form $B$ on $L$ will be taken to be
\[\|B\|=\sup_{\lbrace h:\|h_i\|\leq 1\,\forall i=1,\dots k\rbrace} |B[h_1,...,h_k]|,\]
where $B[h_1,\dots,h_k]$ denotes $B$ applied to arguments $h_1,\dots,h_k\in L$.

 As in \cite{diffusion}, we define $M\subset L$ as a subspace of $L$ consisting of the twice Fr\'echet differentiable functions $f$, such that:
\begin{equation}\label{space_m}
\|D^2f(w+h)-D^2f(w)\|\leq k_f\|h\|\text{,}
\end{equation}
for some constant $k_f$, uniformly in $w,h\in D([0,1],\mathbbm{R}^d)$. 
We have following lemma (whose proof we omit), which may be proved in an analogous way to that used to show (2.6) and (2.7) of \cite{diffusion}:
\begin{lemma}\label{first_der}
For every $f\in M$, let:
\begin{align*}
\|f\|_M:=&\sup_{w\in D([0,1],\mathbbm{R}^d)}\frac{|f(w)|}{1+\|w\|^3}+\sup_{w\in D([0,1],\mathbbm{R}^d)}\frac{\|Df(w)\|}{1+\|w\|^2}+\sup_{w\in D([0,1],\mathbbm{R}^d)}\frac{\|D^2f(w)\|}{1+\|w\|}\\
&+\sup_{w,h\in D([0,1],\mathbbm{R}^d)}\frac{\|D^2f(w+h)-D^2f(w)\|}{\|h\|}.
\end{align*}
Then, for all $f\in M$, we have $\|f\|_M<\infty$.
\end{lemma}

We, furthermore, let $M^0$ be the class of functionals $g\in M$ such that:
\begin{align}
\|g\|_{M^0}:=&\sup_{w\in D([0,1],\mathbbm{R}^d)}|g(w)|+\sup_{w\in D([0,1],\mathbbm{R}^d)}\|Dg(w)\|+\sup_{w\in D([0,1],\mathbbm{R}^d)}\|D^2g(w)\|\nonumber\\
&+\sup_{w,h\in D([0,1],\mathbbm{R}^d)}\frac{\|D^2g(w+h)-D^2g(w)\|}{\|h\|}<\infty\nonumber
\end{align}
and note that $M^0\subset M$. Below, we present a $d$-dimensional version of \cite[Proposition 3.1]{functional_combinatorial} providing conditions, under which weak convergence of the approximating measure to the target one may be deduced from convergence of the corresponding expectations of functions $g\in M^0$. Its proof can be found in the appendix of \cite{kasprzak2}.

\begin{proposition}\label{prop_m}
Suppose that, for each $n\geq 1$, the random element $\mathbf{Y}_n$ of $D([0,1],\mathbbm{R}^d)$ is piecewise constant with intervals of constancy of length at least $r_n$. Let $\left(\mathbf{Z}_n\right)_{n\geq 1}$ be random elements of $D^p$ converging in distribution in $D([0,1],\mathbbm{R}^d)$, with respect to the Skorokhod topology, to a random element $\mathbf{\mathbf{Z}}\in C\left([0,1],\mathbbm{R}^d\right)$. If:
\begin{equation}\label{assumption}
|\mathbbm{E}g(\mathbf{Y}_n)-\mathbbm{E}g(\mathbf{\mathbf{Z}}_n)|\leq C\mathscr{T}_n\|g\|_{M^0}
\end{equation}
for each $g\in M^0$ and if $\mathscr{T}_n\log^2(1/r_n)\xrightarrow{n\to\infty}0$, then the law of $\mathbf{Y}_n$ converges weakly to that of  $\mathbf{\mathbf{Z}}$ in $D([0,1],\mathbbm{R}^d)$, in both the uniform and the Skorokhod topologies.
\end{proposition}
\section{Setting up Stein's method for the pre-limiting approximation}\label{section_setup_stein}
We set up Stein's method in a fashion similar to \cite{diffusion} and \cite{kasprzak2}. First, we define the process $\mathbf{D}_n$ whose distribution will be treated as the target measure. We then construct a process $\left(\mathbf{W}_n(\cdot,u):u\geq 0\right)$ for which the target measure is stationary. We subsequently calculate its infinitesimal generator $\mathcal{A}_n$ and take it as our Stein operator. Next, we solve the Stein equation $\mathcal{A}_nf=g$, using the analysis of \cite{kasprzak}, and prove several smoothness properties of the solution $f_n=\phi_n(g)$.
\subsection{Target measure}
Let
\begin{equation}\label{d_n}
\mathbf{D}_n(t)=\sum_{i_1,\cdots,i_m=1}^{n}\left(\tilde{Z}^{(1)}_{i_1,\cdots,i_m}J^{(1)}_{i_1,\cdots,i_m}(t),\cdots,\tilde{Z}^{(d)}_{i_1,\cdots,i_m}J^{(d)}_{i_1,\cdots,i_m}(t)\right),\quad t\in[0,1],
\end{equation}
where $\tilde{Z}^{(k)}_{i_1,\cdots,i_m}$'s for $k=1,\cdots,d$ are centred Gaussian and:
\begin{enumerate}
\item the covariance matrix $\Sigma_n\in\mathbbm{R}^{(n^md)\times(n^md)}$ of $\tilde{Z}$ is positive definite, for $\tilde{Z}\in \mathbbm{R}^{(n^md)}$ built out of the $\tilde{Z}^{(k)}_{i_1,\cdots,i_m}$'s in such a way that they appear in the lexicographic order with $\tilde{Z}^{(k)}_{i_1,\cdots,i_m}$ appearing before $\tilde{Z}^{(k+1)}_{j_1,\cdots,j_m}$'s for any $k=1,\cdots,d-1$ and $i_1,\cdots,i_m,j_1,\cdots,j_m=1,\cdots,n$;
\item the collection of functions 
\[\left\lbrace J^{(k)}_{i_1,\cdots,i_m}\in D\left([0,1],\mathbbm{R}\right)\,:\, i_1,\cdots,i_m\in\lbrace 1,\cdots, n\rbrace, k\in\lbrace 1,\cdots,p\rbrace\right\rbrace\] is independent of the collection $\left\lbrace \tilde{Z}^{(k)}_{i_1,\cdots,i_m}\,:\, i_1,\cdots,i_m\in\lbrace 1,\cdots, n\rbrace, k\in\lbrace 1,\cdots,p\rbrace\right\rbrace$; a natural example of those would be $J^{(k)}_{i_1,\cdots,i_m}=\mathbbm{1}_{A^{(k)}_{i_1,\cdots,i_m}}$ for some measurable set $A^{(k)}_{i_1,\cdots,i_m}\subset[0,1]$.
\end{enumerate}
\begin{remark}
It is worth noting that processes $\mathbf{D}_n$ of the form (\ref{d_n}) are often approximations of interesting continuous Gaussian processes. An example is  $\mathbf{D}_n$ of (\ref{d_n}), where  all the $\tilde{Z}^{(k)}_{i_1,\cdots,i_m}$'s are standard normal and independent, $m=1$ and $J_i^{(k)}=\mathbbm{1}_{[i/n,1]}$ for all $k=1,\cdots,d$ and $i=1,\cdots,n$. By Donsker's theorem, it approximates the standard Brownian motion. By Proposition \ref{prop_m}, under several assumptions, if a piecewise constant process $\mathbf{Y}_n$ is close enough to process $\mathbf{D}_n$, then the law of $\mathbf{Y}_n$ converges weakly to that of the continuous process that $\mathbf{D}_n$ approximates.
\end{remark}

Now consider an array of i.i.d. Ornstein-Uhlenbeck processes with stationary law $\mathcal{N}(0,1)$, independent of the $J^{(k)}_{i_1,\cdots,i_m}$'s, given by $\lbrace (\mathscr{X}^{(k)}_{i_1,\cdots,i_m}(u),u\geq 0):i_1,\cdots,i_m=1,...,n,\,k=1,...,d\rbrace$. Let $\tilde{\mathscr{U}}(u)=\left(\Sigma_n\right)^{1/2}\mathscr{X}(u)$, where $\Sigma_n$ is the covariance matrix of $\tilde{Z}$, as above, and $\mathscr{X}(u)\in\mathbbm{R}^{n^md}$ is a vector composed out of the $\mathscr{X}^{(k)}_{i_1,\cdots,i_m}(u)$'s in such a way that they are ordered exactly as $\tilde{Z}^{(k)}_{i_1,\cdots,i_m}$'s are ordered in $\tilde{Z}$.
Write $\mathscr{U}_{i_1,\cdots,i_m}^{(k)}(u)=\left(\tilde{\mathscr{U}}(u)\right)_{I(k,i_1,\cdots,i_m)}$ using the bijection $I:\lbrace (k,i_1,\cdots,i_m):i_1,\cdots,i_m=1,\cdots, n, k=1,\cdots, d\rbrace\to\lbrace 1,\cdots,dn^m\rbrace$, given by:
\begin{equation}\label{i}
I(k,i_1,\cdots,i_m)=(k-1)n^m+(i_1-1)n^{m-1}+\cdots+(i_{m-1}-1)n+i_m.
\end{equation}
We will look at the process
\[\mathbf{W}_n(t,u)=\left(\mathbf{W}_n^{(1)}(t,u),\cdots,\mathbf{W}_n^{(d)}(t,u)\right),\quad t\in[0,1],u\geq 0,\]
where, for all $k=1,\cdots,d$:
\[\mathbf{W}_n^{(k)}(t,u)=\sum_{i_1,\cdots,i_m=1}^{n}\mathscr{U}^{(k)}_{i_1,\cdots,i_m}(u)J^{(k)}_{i_1,\cdots,i_m}(t),\quad t\in[0,1],u\geq 0.\]
It is easy to see that the stationary law of the process $\left(\mathbf{W}_n(\cdot,u)\right)_{u\geq 0}$ (which, for any fixed $u$, takes value in $D([0,1],\mathbbm{R}^d)$) is exactly the law of $\mathbf{D}_n$.
\subsection{Stein equation}
The following result follows immediately from \cite[Propositions 4.1 and 4.4]{kasprzak}:
\begin{proposition}\label{prop12.7}
The infinitesimal generator of the process $\left(\mathbf{W}_n(\cdot,u)\right)_{u\geq 0}$ acts on any $f\in M$ (for $M$ defined in Section \ref{section_notation}) in the following way:
\begin{align*}
&\mathcal{A}_nf(w)=-Df(w)[w]+\mathbbm{E}D^2f(w)\left[\mathbf{D}_n,\mathbf{D}_n\right].
\end{align*}
Moreover, for any $g\in M$ such that $\mathbbm{E}g(\mathbf{D}_n)=0$, the Stein equation $\mathcal{A}_nf_n=g$ is solved by:
\begin{equation}\label{phi}
f_n=\phi_n(g)=-\int_0^{\infty}T_{n,u}gdu,
\end{equation}
where $(T_{n,u}f)(w)=\mathbbm{E}\left[f(we^{-u}+\sqrt{1-e^{-2u}}\mathbf{D}_n(\cdot)\right]$. Furthermore, for $g\in M$:
\begin{align}
\text{A)} \quad &\|D\phi_n(g)(w)\|\leq \|g\|_{ M}\left(1+\frac{2}{3}\|w\|^2+\frac{4}{3}\mathbbm{E}\|\mathbf{D}_n\|^2\right)\text{,}\nonumber\\
\text{B)} \quad &\|D^2\phi_n(g)(w)\|\leq \|g\|_{ M}\left(\frac{1}{2}+\frac{\|w\|}{3}+\frac{\mathbbm{E}\|\mathbf{D}_n\|}{3}\right)\text{,}\nonumber\\
\text{C)}\quad&\frac{\left\|D^2\phi_n(g)(w+h)-D^2\phi_n(g)(w)\right\|}{\|h\|}\nonumber\\
\leq&\sup_{w,h\in D^p}\frac{\|D^2(g+c)(w+h)-D^2(g+c)(w)\|}{3\|h\|},\label{m_bound}
\end{align}
for any constant function $c:D([0,1],\mathbbm{R}^d)\to\mathbbm{R}$ and for all $w,h\in D([0,1],\mathbbm{R}^d)$.
\begin{remark}
The fact that the process $\left(\mathbf{W}_n(\cdot,u)\right)_{u\geq 0}$ is built using Ornstein-Uhlenbeck processes and that the corresponding semigroup $T_{n,u}$ takes the convenient form, coming from Mehler's formula, plays an important role in the proof of Proposition \ref{prop12.7}. It is not clear to us whether this result can easily be extended beyond this context.
\end{remark}
\end{proposition}
\section{An abstract approximation theorem}\label{section_abstract}
The following result provides an expression for a bound on the distance between a process $\mathbf{Y}_n$ and $\mathbf{D}_n$, defined by (\ref{d_n}). It assumes that we can find some $\mathbf{Y}_n'$ such that $(\mathbf{Y}_n,\mathbf{Y}_n')$ is an exchangeable pair satisfying an appropriate condition. We explain in Remark \ref{remark_th_1} how our condition is similar to that of  \cite[(1.7)]{reinert_roellin}.
\begin{theorem}\label{theorem1}
Assume that $(\mathbf{Y}_n,\mathbf{Y}_n')$ is an exchangeable pair of $D\left([0,1],\mathbbm{R}^d\right)$-valued random vectors such that:
\begin{equation}
Df(\mathbf{Y}_n)[\mathbf{Y}_n]=2\mathbbm{E}^{\mathbf{Y}_n}Df(\mathbf{Y}_n)\left[(\mathbf{Y}_n-\mathbf{Y}_n')\Lambda_n\right]+R_f,
\label{condition}
\end{equation}
where $\mathbbm{E}^{\mathbf{Y}_n}[\cdot]:=\mathbbm{E}\left[\cdot|\mathbf{Y}_n\right]$, for all $f\in M$, some $\Lambda_n\in\mathbbm{R}^{d\times d}$ and some random variable $R_f=R_f(\mathbf{Y}_n)$.  Let $\mathbf{D}_n$ be defined by (\ref{d_n}). Then, for any $g\in M$:
\begin{align*}
\left|\mathbbm{E}g(\mathbf{Y}_n)-\mathbbm{E}g(\mathbf{D}_n)\right|\leq \epsilon_1+\epsilon_2+\epsilon_3,
\end{align*}
where
\begin{align*}
\epsilon_1&=\frac{\|g\|_M}{6}\mathbbm{E}\left[\|(\mathbf{Y}_n-\mathbf{Y}_n')\Lambda_n\|\|\mathbf{Y}_n-\mathbf{Y}_n'\|^2\right],\\
\epsilon_2&=\left|\mathbbm{E}D^2f(\mathbf{Y}_n)\left[(\mathbf{Y}_n-\mathbf{Y}_n')\Lambda_n,\mathbf{Y}_n-\mathbf{Y}_n'\right]-\mathbbm{E}D^2f(\mathbf{Y}_n)\left[\mathbf{D}_n,\mathbf{D}_n\right]\right|,\\
\epsilon_3&=|\mathbbm{E}R_f|,
\end{align*}
and $f=\phi_n(g)$, as defined by (\ref{phi}).
\end{theorem}
\begin{remark}[Relevance of terms in the bound]
Term $\epsilon_1$ measures how close $\mathbf{Y}_n$ and $\mathbf{Y}_n'$ are and how \textit{small} (in a certain sense) $\Lambda_n$ is. Term $\epsilon_2$ quantifies the difference between the covariance structures of $\mathbf{Y}_n-\mathbf{Y}_n'$ and $\mathbf{D}_n$. This term may be estimated in several applications (see Theorems \ref{theorem_weighted_pre} and \ref{theorem_pre_limiting} below), yet this often requires some effort. Term $\epsilon_3$ measures the error in the exchangeable-pair linear regression condition (\ref{condition}).

\end{remark}
\begin{remark}
Condition (\ref{condition}) is always satisfied, for example with $\Lambda_n=0$ and $R_f=Df(\mathbf{Y}_n)[\mathbf{Y}_n]$ for all $f\in M$. However, for the bound in Theorem \ref{theorem1} to be small, we require the expectation of $R_f$ to be small in absolute value.
\end{remark}
\begin{remark}
The term
\[\left|\mathbbm{E}D^2f(\mathbf{Y}_n)\left[(\mathbf{Y}_n-\mathbf{Y}_n')\Lambda_n,\mathbf{Y}_n-\mathbf{Y}_n'\right]-\mathbbm{E}D^2f(\mathbf{Y}_n)\left[\mathbf{D}_n,\mathbf{D}_n\right]\right|\] in the bound obtained in Theorem \ref{theorem1} is an analogue of the second condition in \cite[Theorem 3]{meckes09}. The main result of that paper provides a bound on approximation by $\mathcal{N}(0,\Sigma)$ of a $d$-dimensional vector $X$. This is achieved by constructing an exchangeable pair $(X,X')$ satisfying:
\[\mathbbm{E}^X[X'-X]=\Lambda X+E\quad\text{and}\quad\mathbbm{E}^X[(X'-X)(X'-X)^T]=2\Lambda\Sigma+E'\]
for some invertible matrix $\Lambda$ and some remainder terms $E$ and $E'$. In the same spirit, Theorem \ref{theorem1} could be rewritten to assume (\ref{condition}) and:
\[ \mathbbm{E}^{\mathbf{Y}_n}D^2f(\mathbf{Y}_n)\left[(\mathbf{Y}_n-\mathbf{Y}_n')\Lambda_n,\mathbf{Y}_n-\mathbf{Y}_n'\right]=D^2f(\mathbf{Y}_n)\left[\mathbf{D}_n,\mathbf{D}_n\right]+R^1_f,\]
for all $f\in M$. The bound would then take the form:
\begin{align*}
\left|\mathbbm{E}g(\mathbf{Y}_n)-\mathbbm{E}g(\mathbf{D}_n)\right|\leq&\frac{\|g\|_M}{6}\mathbbm{E}\left[\|(\mathbf{Y}_n-\mathbf{Y}_n')\Lambda_n\|\|\mathbf{Y}_n-\mathbf{Y}_n'\|^2\right]+|\mathbbm{E}R_f|+|\mathbbm{E}R^1_f|,
\end{align*}
for $f=\phi_n(g)$.
\end{remark}
\begin{remark}\label{remark_th_1}
The role of $\Lambda_n$ in condition (\ref{condition}) is equivalent to that played by $\Lambda^{-1}$ in \cite{reinert_roellin} for $\Lambda$ defined by (1.7) therein.  In the functional setting, condition (\ref{condition}) is more appropriate than a straightforward adaptation of the setup of \cite{reinert_roellin}. This is because, for general processes $\mathbf{Y}_n$, the properties of the Fr\'echet derivative do not allow us to treat evaluating the derivative in the direction of $\mathbf{Y}_n-\mathbf{Y}_n'$ as matrix multiplication. Indeed, multiplying both sides of the hypothetical condition:
\[-Df(\mathbf{Y}_n)[\Lambda \mathbf{Y}_n]=\mathbbm{E}^{\mathbf{Y}_n}Df(\mathbf{Y}_n)[\mathbf{Y}_n-\mathbf{Y}_n']\]
by $\Lambda^{-1}$ does not yield:
\[-Df(\mathbf{Y}_n)[\mathbf{Y}_n]=\mathbbm{E}^{\mathbf{Y}_n}Df(\mathbf{Y}_n)[\Lambda^{-1}(\mathbf{Y}_n-\mathbf{Y}_n')].\]
\end{remark}
\begin{proof}[Proof of Theorem \ref{theorem1}]
We will bound $\left|\mathbbm{E}g(\mathbf{Y}_n)-\mathbbm{E}g(\mathbf{D}_n)\right|$ by bounding $\left|\mathbbm{E}\mathcal{A}_nf(\mathbf{Y}_n)\right|$, where $f$ is the solution to the Stein equation:
\[\mathcal{A}_nf=g-\mathbbm{E}g(\mathbf{D}_n),\]
for $\mathcal{A}_n$ defined in Proposition \ref{prop12.7}. Note that, by exchangeability of $(\mathbf{Y}_n,\mathbf{Y}_n')$ and (\ref{condition}):
\begin{align*}
0=&\mathbbm{E}\left(Df(\mathbf{Y}_n')+Df(\mathbf{Y}_n)\right)\left[(\mathbf{Y}_n-\mathbf{Y}_n')\Lambda_n\right]\\
=&\mathbbm{E}\left(Df(\mathbf{Y}_n')-Df(\mathbf{Y}_n)\right)\left[(\mathbf{Y}_n-\mathbf{Y}_n')\Lambda_n\right]+2\mathbbm{E}\left\lbrace\mathbbm{E}^{\mathbf{Y}_n}Df(\mathbf{Y}_n)\left[(\mathbf{Y}_n-\mathbf{Y}_n')\Lambda_n\right]\right\rbrace\\
=&\mathbbm{E}\left(Df(\mathbf{Y}_n')-Df(\mathbf{Y}_n)\right)\left[(\mathbf{Y}_n-\mathbf{Y}_n')\Lambda_n\right]+\mathbbm{E}Df(\mathbf{Y}_n)[\mathbf{Y}_n]-\mathbbm{E}R_f
\end{align*}
and so:
\begin{equation*}
\mathbbm{E}Df(\mathbf{Y}_n)[\mathbf{Y}_n]=\mathbbm{E}\left(Df(\mathbf{Y}_n)-Df(\mathbf{Y}_n')\right)\left[(\mathbf{Y}_n-\mathbf{Y}_n')\Lambda_n\right]+\mathbbm{E}R_f.
\end{equation*}
Therefore:
\begin{align}
&\left|\mathbbm{E}\mathcal{A}_nf(\mathbf{Y}_n)\right|\notag\\
=&\left|\mathbbm{E}Df(\mathbf{Y}_n)[\mathbf{Y}_n]-\mathbbm{E}D^2f(\mathbf{Y}_n)\left[\mathbf{D}_n,\mathbf{D}_n\right]\right|\nonumber\\
=&\left|\mathbbm{E}\left(Df(\mathbf{Y}_n)-Df(\mathbf{Y}_n')\right)\left[(\mathbf{Y}_n-\mathbf{Y}_n')\Lambda_n\right]-\mathbbm{E}D^2f(\mathbf{Y}_n)\left[\mathbf{D}_n,\mathbf{D}_n\right]+\mathbbm{E}R_f\right|\nonumber\\
\leq& \left|\mathbbm{E}\left(Df(\mathbf{Y}_n)-Df(\mathbf{Y}_n')\right)\left[(\mathbf{Y}_n-\mathbf{Y}_n')\Lambda_n\right]-\mathbbm{E}D^2f(\mathbf{Y}_n')\left[(\mathbf{Y}_n-\mathbf{Y}_n')\Lambda_n,\mathbf{Y}_n-\mathbf{Y}_n'\right]\right|\nonumber\\
&+\left|\mathbbm{E}D^2f(\mathbf{Y}_n)\left[(\mathbf{Y}_n-\mathbf{Y}_n')\Lambda_n,\mathbf{Y}_n-\mathbf{Y}_n'\right]-\mathbbm{E}D^2f(\mathbf{Y}_n)\left[\mathbf{D}_n,\mathbf{D}_n\right]\right|+|\mathbbm{E}R_f|\nonumber\\
\leq&\frac{\|g\|_M}{6}\mathbbm{E}\left[\|(\mathbf{Y}_n-\mathbf{Y}_n')\Lambda_n\|\|\mathbf{Y}_n-\mathbf{Y}_n'\|^2\right]+|\mathbbm{E}R_f|\nonumber\\
&+\left|\mathbbm{E}D^2f(\mathbf{Y}_n)\left[(\mathbf{Y}_n-\mathbf{Y}_n')\Lambda_n,\mathbf{Y}_n-\mathbf{Y}_n'\right]-\mathbbm{E}D^2f(\mathbf{Y}_n)\left[\mathbf{D}_n,\mathbf{D}_n\right]\right|,\nonumber
\end{align}
where the last inequality follows by Taylor's theorem and Proposition \ref{prop12.7}.
\end{proof}

\section{Weighted, degenerate $U$-statistics}\label{weighted}
In this Section we will apply Theorem \ref{theorem1} in order to prove bounds for the approximation of a vector of weighted, degenerate
$U$-processes by suitable Gaussian processes. 
\subsection{Introduction}\label{intro_weighted}
The setup will be the following. We fix positive integers $d,p_1,\dotsc,p_d$ and consider a sequence $(X_i)_{i\in\mathbbm{N}}$ of i.i.d. random variables with distribution $\mu$ on some measurable space $(E,\mathcal{E})$. Moreover, for $1\leq i\leq d$, we let $\psi(i)\in L^2(\mu^{p_i})$ be a symmetric kernel such that
$\mathbbm{E}[\psi(i)^2(X_1,\dotsc,X_{p_i})]>0$. We assume that $\psi(i)$ is \textbf{(completely) degenerate} with respect to $\mu$, i.e. that
\[
\mathbbm{E}[\left.\psi(i)(X_1,\dotsc,X_{p_i})\,\right|\,X_1,\dotsc,X_{p_i-1}]=0,\quad\text{a.s.}
\]
 We denote by
$\mathcal{D}_p(n)$ the collection of $p$-subsets of the set $[n]:=\{1,\dotsc,n\}$ (if $p>n$, we set $\mathcal{D}_p(n)=\emptyset$).

Furthermore, we fix an integer $n\geq \max(p_1,\dotsc,p_d)$ and let $\lbrace a_J(i):\, 1\leq i\leq d,\,
J\in\mathcal{D}_{p_i}(n)\rbrace$, be a (given) set of real numbers (weights).
We further let $\lbrace\sigma_n(i):\,1\leq i\leq d\rbrace$ be a set of positive real numbers and, for $t\in[0,1]$, define
\[\mathbf{Y}_n^{(i)}(t):=\frac{1}{\sigma_n(i)}\sum_{J\in\mathcal{D}_{p_i}(\lfloor nt\rfloor)}
 a_J(i)\psi(i)(X_j,j\in J)\,.
\]
In some applications it may be natural to take
\[\sigma_n(i)^2=\mathbbm{E}[\psi(i)^2(X_1,\dotsc,X_{p_i})]\sum_{J\in\mathcal{D}_{p_i}(n)}a_J(i)^2\,,\quad 1\leq i\leq d,\]
i.e. equal to the variance of the sum in the definition of $\mathbf{Y}_n^{(i)}(1)$. This is, however, not necessary for our results.
For fixed $t$ (in particular for $t=1$), the quantity $\mathbf{Y}_n^{(i)}(t)$ is customarily referred to as a \textbf{degenerate, weighted $U$-statistic} based on $X_1,\dotsc, X_{\lfloor nt\rfloor}$ and, thus, we coin the whole random function $\mathbf{Y}_n^{(i)}$ a \textbf{degenerate, weighted $U$-process} . Limit theorems (not necessarily central) for such weighted $U$-statistics have been derived in \cite{OnRe, RiRo97, Major, RiUt} and in the (somehow) more special case of incomplete $U$-statistics in \cite{Ja, Blom, BrKil}. However, we have not been able to find FCLTs for \textbf{degenerate, weighted $U$-process} in the literature.

With the above definitions, we let
\[\mathbf{Y}_n:=(\mathbf{Y}_n^{(1)},\dotsc,\mathbf{Y}_n^{(d)})\,,\]
which is, as one can easily observe, an element of $D([0,1],\mathbbm{R}^d)$.
We will write $X:=(X_1,\dotsc,X_n)$ and construct an $X':=(X_1',\dotsc,X_n')$ such that the pair $(X,X')$ is exchangeable. Specifically, we let $X_0$ be another random variable with distribution $\mu$ and let $I$ be uniformly distributed on $[n]$ in such a way that $I,X_0,(X_j)_{j\in\mathbbm{N}}$ are jointly independent. For $1\leq j\leq n$, we let
\begin{equation*}
 X_j':=\begin{cases}
        X_j\,,&\text{if }j\not=I\\
        X_0\,,&\text{if }j=I\,.
       \end{cases}
\end{equation*}
Then, for $t\in[0,1]$ and $1\leq i\leq d$, we define
\[(\mathbf{Y}_n^{(i)})'(t):=\frac{1}{\sigma_n(i)}\sum_{J\in\mathcal{D}_{p_i}(\lfloor nt\rfloor)}
 a_J(i)\psi(i)(X'_j,j\in J)\]
and
 \[\mathbf{Y}'_n:=((\mathbf{Y}_n^{(1)})',\dotsc,(\mathbf{Y}_n^{(d)})')\,.\]
The pair $(\mathbf{Y}_n,\mathbf{Y}'_n)$ is clearly exchangeable and, for $f\in M$, similarly as in the proof of \cite[Lemma 2.3]{DP16}, one can use degeneracy to show that
\begin{equation*}
Df(\mathbf{Y}_n)[\mathbf{Y}_n]=2\mathbbm{E}^{\mathbf{Y}_n}Df(\mathbf{Y}_n)\left[(\mathbf{Y}_n-\mathbf{Y}_n')\Lambda_n\right],
\end{equation*}
where
\begin{equation}\label{lambda_weighted}
\Lambda_n=\text{diag}\left(\frac{n}{2p_1},\dots,\frac{n}{2p_d}\right).
\end{equation}
Therefore condition (\ref{condition}) is satisfied for $\Lambda_n$ of (\ref{lambda_weighted}) and $R_f=0$. In what follows we will assume that
$1\leq p_1\leq p_2\leq\cdots\leq p_d$.

\subsection{A pre-limiting process}\label{pre_lim_weighted}
We will construct a pre-limiting Gaussian process $\mathbf{D}_n$ of the form (\ref{d_n}) which has the same covariance structure as $\mathbf{Y}_n$. We take $\mathbf{D}_n=\left(\mathbf{D}_n^{(1)},\dots,\mathbf{D}_n^{(d)}\right)$ for
\begin{equation*}
\mathbf{D}_n^{(i)}(t)=\frac{1}{\sigma_n(i)}\sum_{J\in\mathcal{D}_{p_i}(\lfloor nt\rfloor)}
 a_J(i)Z_J(i),
\end{equation*}
where, for $i=1,\dots,d$ and $J\in\mathcal{D}_{p_i}(n)$, $Z_J(i)$  are jointly Gaussian random variables that are independent of $X$ and satisfy
\begin{equation*}
\mathbbm{E}\left[Z_J(i)Z_K(l)\right]=\begin{cases}
\mathbbm{E}[\psi(i)(X_1,\dotsc,X_{p_i})\psi(l)(X_1,\dotsc,X_{p_l}) ],&\text{if }p_i=p_l\text{ and } K=J\\
0,&\text{otherwise,}\end{cases}
\end{equation*}
for $i,l=1,\dots,d$, $J\in\mathcal{D}_{p_i}(n)$ and $K\in\mathcal{D}_{p_l}(n)$.

\subsection{Distance from the pre-limiting process}
Having established the setup and defined the pre-limiting process above, we prove the following result:
\begin{theorem}\label{theorem_weighted_pre}
Let $\mathbf{Y}_n$ be defined as in Section \ref{intro_weighted} and $\mathbf{D}_n$ be defined as in Section \ref{pre_lim_weighted}. Then, for any $g\in M$,

\begin{align*}
&\Bigl|\mathbbm{E}[g(\mathbf{Y_n})]-\mathbbm{E}[g(\mathbf{D}_n)]\Bigr|
\leq\frac{2\sqrt{d}\|g\|_M}{3p_1}\sum_{i=1}^d\frac{\|\psi(i)\|_{L^3(\mu^{p_i})}^3 }{\sigma_n(i)^3}\sum_{l=1}^n
\left(\sum_{\substack{J\in\mathcal{D}_{p_i}(n):\\ l\in J}} |a_J(i)|\right)^3\notag\\
&\;+ \|g\|_M \sum_{i,j,k=1}^d \frac{\|\psi(i)\|_{L^3(\mu^{p_i})}\|\psi(j)\|_{L^3(\mu^{p_j})}\|\psi(k)\|_{L^3(\mu^{p_k})}}{\sigma_n(i)\sigma_n(j)\sigma_n(k)}\sum_{\substack{J\in\mathcal{D}_{p_i}(n),\\K\in\mathcal{D}_{p_j}(n),\\L\in \mathcal{D}_{p_k}(n):\\J\cap K\not=\emptyset,\\L\cap(J\cup K)\not=\emptyset}}|a_J(i)a_K(j)a_L(k)|.
\end{align*}
\end{theorem}

\begin{proof}~\\
\textbf{Step 1.} First note that, for $\epsilon_1$ in Theorem \ref{theorem1},
\begin{equation}\label{s1_1}
\left\|(\mathbf{Y}_n-\mathbf{Y}_n')\Lambda_n\right\|\left\|\mathbf{Y}_n-\mathbf{Y}_n'\right\|^2\leq\frac{n}{2p_1}\left\|\mathbf{Y}_n-\mathbf{Y}_n'\right\|^3,
\end{equation}
which follows directly from the definition of $\Lambda_n$ in (\ref{lambda_weighted}) and our assumption that $p_1\leq\dots\leq p_d$.
Now, note that
\begin{align}
\left\|\mathbf{Y}_n-\mathbf{Y}_n'\right\|^3=&\sup_{t\in[0,1]}\left[\left(\mathbf{Y}_n^{(1)}(t)-\left(\mathbf{Y}_n^{(1)}\right)'(t)\right)^2+\dots+\left(\mathbf{Y}_n^{(d)}(t)-\left(\mathbf{Y}_n^{(d)}\right)'(t)\right)^2\right]^{3/2}\nonumber\\
\leq&\sqrt{d}\sup_{t\in[0,1]}\left[\left|\mathbf{Y}_n^{(1)}(t)-\left(\mathbf{Y}_n^{(1)}\right)'(t)\right|^3+\dots+\left|\mathbf{Y}_n^{(d)}(t)-\left(\mathbf{Y}_n^{(d)}\right)'(t)\right|^3\right]\nonumber\\
\leq&\sqrt{d}\left[\left\|\mathbf{Y}_n^{(1)}-\left(\mathbf{Y}_n^{(1)}\right)'\right\|^3+\dots+\left\|\mathbf{Y}_n^{(d)}-\left(\mathbf{Y}_n^{(d)}\right)'\right\|^3\right].\label{s1_2}
\end{align}
 Furthermore, for $\max(J):=\max\lbrace j\,:\,j\in J\rbrace$ and for all $i=1,\dots,d$:
\begin{align}
&\mathbbm{E}\left\|\mathbf{Y}_n^{(i)}-\left(\mathbf{Y}_n^{(i)}\right)'\right\|^3\nonumber\\
&=\frac{1}{\sigma_n(i)^3}\notag\\
&\cdot\mathbbm{E}\left\{\sup_{t\in[0,1]}\left| \sum_{\substack{J\in\mathcal{D}_{p_i}(\lfloor nt\rfloor):\\ I\in J}} a_J(i)\bigl(\psi(i)(X_j,j\in J)- \psi(i)(X_0,X_j,j\in J\setminus\{I\})\bigr)
\mathbbm{1}_{[\frac{\max(J)}{n},1]}(t)\right|^3\right\}\notag\\
&\leq \frac{1}{\sigma_n(i)^3}\mathbbm{E}\left(\sum_{\substack{J\in\mathcal{D}_{p_i}(n):\\ I\in J}} |a_J(i)|\bigl|\psi(i)(X_j,j\in J)- \psi(i)(X_0,X_j,j\in J\setminus\{I\})\bigr|\right)^3\notag\\
&\leq\frac{1}{n\sigma_n(i)^3}\sum_{l=1}^n \sum_{\substack{J,K,L\in\mathcal{D}_{p_i}(n):\\ l\in J\cap K\cap L}} |a_J(i)a_K(i)a_L(i)| \mathbbm{E}\Biggl[\bigl|\psi(i)(X_j,j\in J)- \psi(i)(X_0,X_j,j\in J\setminus\{l\})\bigr|\notag\\
&\cdot\bigl|\psi(i)(X_j,j\in K)- \psi(i)(X_0,X_j,j\in K\setminus\{l\})\bigr|\bigl|\psi(i)(X_j,j\in L)- \psi(i)(X_0,X_j,j\in L\setminus\{l\})\bigr|\Biggr]\notag\\
&\leq \frac{\mathbbm{E}\bigl|\psi(i)(X_1,\dotsc,X_{p_i})-\psi(i)(X_2,\dotsc,X_{p_{i+1}})\bigr|^3 }{n\sigma_n(i)^3}\sum_{l=1}^n\sum_{\substack{J,K,L\in\mathcal{D}_{p_i}(n):\\ l\in J\cap K\cap L}} |a_J(i)a_K(i)a_L(i)|\label{s1_3a}\\
&\leq \frac{8\mathbbm{E}\bigl|\psi(i)(X_1,\dotsc,X_{p_i})\bigr|^3 }{n\sigma_n(i)^3}\sum_{l=1}^n
\left(\sum_{\substack{J\in\mathcal{D}_{p_i}(n):\\ l\in J}} |a_J(i)|\right)^3\label{s1_3b}\,.
\end{align}
Combining \eqref{s1_1} -\eqref{s1_3b} we obtain
\begin{align}
\epsilon_1&\leq \frac{\sqrt{d}\|g\|_M}{12p_1}\sum_{i=1}^d\frac{\mathbbm{E}\bigl|\psi(i)(X_1,\dotsc,X_{p_i})-\psi(i)(X_2,\dotsc,X_{p_{i+1}})\bigr|^3 }{\sigma_n(i)^3}\notag\\
&\hspace{1.7cm}\cdot\sum_{l=1}^n\sum_{\substack{J,K,L\in\mathcal{D}_{p_i}(n):\\ l\in J\cap K\cap L}} |a_J(i)a_K(i)a_L(i)|\label{eps_1_weighted_a}\\
&\leq\frac{2\sqrt{d}\|g\|_M}{3p_1}\sum_{i=1}^d\frac{\mathbbm{E}\bigl|\psi(i)(X_1,\dotsc,X_{p_i})\bigr|^3 }{\sigma_n(i)^3}\sum_{l=1}^n
\left(\sum_{\substack{J\in\mathcal{D}_{p_i}(n):\\ l\in J}} |a_J(i)|\right)^3.\label{eps_1_weighted_b}
\end{align}
\textbf{Step 2.}
We will now bound $\epsilon_2$ of Theorem \ref{theorem1}. Denoting by $e_i$ the $i$th element of the canonical basis of $\mathbbm{R}^d$, for $i=1,\dots,d$, for any $f\in M$, we have
\begin{align}\label{s2a}
&D^2f(\mathbf{Y}_n)\left[\left(\mathbf{Y}_n-\mathbf{Y}_n'\right)\Lambda_n,\mathbf{Y}_n-\mathbf{Y}_n'\right]\nonumber\\
=&D^2f(\mathbf{Y}_n)\left[\sum_{i=1}^d\frac{n}{2p_i}\left(\mathbf{Y}^{(i)}_n-\left(\mathbf{Y}^{(i)}_n\right)'\right)e_i,\sum_{i=1}^d\left(\mathbf{Y}^{(i)}_n-\left(\mathbf{Y}^{(i)}_n\right)'\right)e_i\right]\nonumber\\
=&\sum_{i,j=1}^d\frac{n}{2p_i}D^2f(\mathbf{Y}_n)\left[\left(\mathbf{Y}^{(i)}_n-\left(\mathbf{Y}^{(i)}_n\right)'\right)e_i,\left(\mathbf{Y}^{(j)}_n-\left(\mathbf{Y}^{(j)}_n\right)'\right)e_j\right].
\end{align}
We now let $f=\phi_n(g)$, as defined by \eqref{phi}, and fix some $i,j\in\lbrace 1,\dots,d\rbrace$. We have that
\begin{align}\label{s2b}
&\Bigg|\frac{n}{2p_i}\mathbbm{E}D^2f(\mathbf{Y}_n)\left[\left(\mathbf{Y}^{(i)}_n-\left(\mathbf{Y}^{(i)}_n\right)'\right)e_i,\left(\mathbf{Y}^{(j)}_n-\left(\mathbf{Y}^{(j)}_n\right)'\right)e_j\right]\notag\\
&\hspace{9cm}-\mathbbm{E}D^2f(\mathbf{Y}_n)\left[\mathbf{D}_n^{(i)}e_i,\mathbf{D}_n^{(j)}e_j\right]\Bigg|\nonumber\\
&=\frac{1}{\sigma_n(i)\sigma_n(j)}\Biggl|\frac{n}{2p_i}
\sum_{\substack{J\in\mathcal{D}_{p_i}(n),\\
K\in\mathcal{D}_{p_j}(n)}}a_J(i)a_K(j)
\mathbbm{E}\Bigl[
\bigl(\psi(i)(X_u,u\in J)-\psi(i)(X'_u,u\in J) \bigr)\notag\\
&\hspace{2cm}\cdot
\bigl(\psi(j)(X_u,u\in K)-\psi(j)(X'_u,u\in K) \bigr)D^2f(\mathbf{Y}_n)
\bigl[\mathbbm{1}_{[\frac{\max(J)}{n},1]}e_i,\mathbbm{1}_{[\frac{\max(K)}{n},1]}e_j\bigr]\Bigr]\notag\\
&\hspace{3cm}-\sum_{\substack{J\in\mathcal{D}_{p_i}(n),\\
K\in\mathcal{D}_{p_j}(n)}}a_J(i)a_K(j)\mathbbm{1}_{\{J=K\}}
\mathbbm{E}\bigl[\psi(i)(X_1,\dotsc,X_{p_i})\psi(j)(X_1,\dotsc,X_{p_j})\bigr]\notag\\
&\hspace{7cm}\cdot\mathbbm{E}\Bigl[D^2f(\mathbf{Y}_n)
\bigl[\mathbbm{1}_{[\frac{\max(J)}{n},1]}e_i,\mathbbm{1}_{[\frac{\max(K)}{n},1]}e_j\bigr]\Bigr]\Biggr|\notag\\
&=\frac{1}{2p_i\sigma_n(i)\sigma_n(j)}\Biggl|\sum_{l=1}^n\sum_{\substack{J\in\mathcal{D}_{p_i}(n),\\
K\in\mathcal{D}_{p_j}(n),\\l\in J\cap K  }}a_J(i)a_K(j)
\mathbbm{E}\Bigl[
\bigl(\psi(i)(X_u,u\in J)-\psi(i)(X_0,X_u,u\in J\setminus\{l\}) \bigr)\notag\\
&\hspace{0.5cm}\cdot
\bigl(\psi(j)(X_u,u\in K)-\psi(j)(X_0,X_u,u\in K\setminus\{l\}) \bigr)D^2f(\mathbf{Y}_n)
\bigl[\mathbbm{1}_{[\frac{\max(J)}{n},1]}e_i,\mathbbm{1}_{[\frac{\max(K)}{n},1]}e_j\bigr]\Bigr]\notag\\
&-2\sum_{l=1}^n\sum_{\substack{J\in\mathcal{D}_{p_i}(n),\\
K\in\mathcal{D}_{p_j}(n),\\ l\in J\cap K}}a_J(i)a_K(j)\mathbbm{1}_{\{J=K\}}
\mathbbm{E}\bigl[\psi(i)(X_1,\dotsc,X_{p_i})\psi(j)(X_1,\dotsc,X_{p_j})\bigr]\notag\\
&\hspace{7cm}\cdot\mathbbm{E}\Bigl[D^2f(\mathbf{Y}_n)
\bigl[\mathbbm{1}_{[\frac{\max(J)}{n},1]}e_i,\mathbbm{1}_{[\frac{\max(K)}{n},1]}e_j\bigr]\Bigr]\Biggr|\notag\\
&=\frac{1}{2p_i\sigma_n(i)\sigma_n(j)}\Biggl|\sum_{l=1}^n\sum_{\substack{J\in\mathcal{D}_{p_i}(n),\\
K\in\mathcal{D}_{p_j}(n),\\l\in J\cap K  }}\hspace{-3mm}a_J(i)a_K(j)
\mathbbm{E}\Biggl[\Biggl(
\biggl(\psi(i)(X_u,u\in J)\notag\\
&\hspace{1.5cm}-\psi(i)(X_0,X_u,u\in J\setminus\{l\}) \biggr)\cdot
\biggl(\psi(j)(X_u,u\in K)-\psi(j)(X_0,X_u,u\in K\setminus\{l\}) \biggr)\notag\\
&-2\mathbbm{1}_{\{J=K\}}\mathbbm{E}\bigl[\psi(i)(X_1,\dotsc,X_{p_i})\psi(j)(X_1,\dotsc,X_{p_j})\bigr]\Biggr)D^2f(\mathbf{Y}_n)
\bigl[\mathbbm{1}_{[\frac{\max(J)}{n},1]}e_i,\mathbbm{1}_{[\frac{\max(K)}{n},1]}e_j\bigr]\Biggr]\Biggr|.
\end{align}
Now, we define
\[\mathbf{Y}_n^{J,K}:=\Bigl(\bigl(\mathbf{Y}_n^{J,K}\bigr)^{(1)},\cdots,\bigl(\mathbf{Y}_n^{J,K}\bigr)^{(d)}\Bigr)\]
via
\begin{equation*}
 \bigl(\mathbf{Y}_n^{J,K}\bigr)^{(i)}(t):=\frac{1}{\sigma_n(i)}\sum_{\substack{L\in\mathcal{D}_{p_i}(\lfloor nt\rfloor):\\ L\cap(J\cup K)=\emptyset}}
 a_J(i)\psi(i)(X_j,j\in L)\,,\quad 1\leq i\leq d,\, t\in[0,1]\,.
\end{equation*}
Then, using independence, from \eqref{s2b} we obtain that
\begin{align}\label{s2c}
&\Bigg|\frac{n}{2p_i}\mathbbm{E}D^2f(\mathbf{Y}_n)\left[\left(\mathbf{Y}^{(i)}_n-\left(\mathbf{Y}^{(i)}_n\right)'\right)e_i,\left(\mathbf{Y}^{(j)}_n-\left(\mathbf{Y}^{(j)}_n\right)'\right)e_j\right]\notag\\
&\hspace{9cm}-\mathbbm{E}D^2f(\mathbf{Y}_n)\left[\mathbf{D}_n^{(i)}e_i,\mathbf{D}_n^{(j)}e_j\right]\Bigg|\nonumber\\
&=\frac{1}{2p_i\sigma_n(i)\sigma_n(j)}\Biggl|\sum_{l=1}^n\sum_{\substack{J\in\mathcal{D}_{p_i}(n),\\
K\in\mathcal{D}_{p_j}(n),\\l\in J\cap K  }}a_J(i)a_K(j)
\mathbbm{E}\Biggl[\Biggl(
\biggl(\psi(i)(X_u,u\in J)\notag\\
&\hspace{2cm}-\psi(i)(X_0,X_u,u\in J\setminus\{l\}) \biggr)
\biggl(\psi(j)(X_u,u\in K)-\psi(j)(X_0,X_u,u\in K\setminus\{l\}) \biggr)
\notag\\
&\hspace{5.5cm}-2\mathbbm{1}_{\{J=K\}}\mathbbm{E}\bigl[\psi(i)(X_1,\dotsc,X_{p_i})\psi(j)(X_1,\dotsc,X_{p_j})\bigr]\Biggr)\notag\\
&\hspace{5cm}\cdot\Bigl(D^2f(\mathbf{Y}_n)-D^2f(\mathbf{Y}_n^{J,K})\bigr)
\bigl[\mathbbm{1}_{[\frac{\max(J)}{n},1]}e_i,\mathbbm{1}_{[\frac{\max(K)}{n},1]}e_j\bigr]\Biggr]\Biggr|\notag\\
&\stackrel{\eqref{m_bound}C}\leq \frac{\|g\|_M}{6p_i\sigma_n(i)\sigma_n(j)}\sum_{l=1}^n\sum_{\substack{J\in\mathcal{D}_{p_i}(n),\\
K\in\mathcal{D}_{p_j}(n),\\l\in J\cap K  }}|a_J(i)a_K(j)|
\mathbbm{E}\Biggl[\Biggl|
\biggl(\psi(i)(X_u,u\in J)\notag\\
&\hspace{1.5cm}-\psi(i)(X_0,X_u,u\in J\setminus\{l\}) \biggr)
\biggl(\psi(j)(X_u,u\in K)-\psi(j)(X_0,X_u,u\in K\setminus\{l\}) \biggr)\notag\\
&\hspace{3cm}
-2\mathbbm{1}_{\{J=K\}}\mathbbm{E}\bigl[\psi(i)(X_1,\dotsc,X_{p_i})\psi(j)(X_1,\dotsc,X_{p_j})\bigr]\Biggr|\cdot\|\mathbf{Y}_n-\mathbf{Y}_n^{J,K}\|   \Biggr].
\end{align}
Now, we observe that
\begin{align*}
 \|\mathbf{Y}_n-\mathbf{Y}_n^{J,K}\|&\leq
 \sum_{k=1}^d\frac{1}{\sigma_n(k)}\|\mathbf{Y}_n^{(k)}-(\mathbf{Y}_n^{J,K})^{(k)}\|\notag\\
& \leq\sum_{k=1}^d\frac{1}{\sigma_n(k)}\sum_{\substack{L\in\mathcal{D}_{p_k}(n):\\ L\cap(J\cup K)\not=\emptyset}}
 |a_J(k)| |\psi(k)(X_u,u\in L)|.
\end{align*}
Hence, \eqref{s2c} yields
\begin{align}\label{s2d}
&\Bigg|\frac{n}{2p_i}\mathbbm{E}D^2f(\mathbf{Y}_n)\left[\left(\mathbf{Y}^{(i)}_n-\left(\mathbf{Y}^{(i)}_n\right)'\right)e_i,\left(\mathbf{Y}^{(j)}_n-\left(\mathbf{Y}^{(j)}_n\right)'\right)e_j\right]\notag\\
&\hspace{9cm}-\mathbbm{E}D^2f(\mathbf{Y}_n)\left[\mathbf{D}_n^{(i)}e_i,\mathbf{D}_n^{(j)}e_j\right]\Bigg|\nonumber\\
&\leq\sum_{k=1}^d \frac{\|g\|_M}{6p_i\sigma_n(i)\sigma_n(j)\sigma_n(k)}\sum_{l=1}^n\sum_{\substack{J\in\mathcal{D}_{p_i}(n),\\
K\in\mathcal{D}_{p_j}(n),\\l\in J\cap K  }}\sum_{\substack{L\in\mathcal{D}_{p_k}(n):\\ L\cap(J\cup K)\not=\emptyset}}|a_J(i)a_K(j)a_L(k)|\mathbbm{E}\Biggl[\biggl|
\Bigl(\psi(i)(X_u,u\in J)\notag\\
&\hspace{1.5cm}-\psi(i)(X_0,X_u,u\in J\setminus\{l\}) \Bigr)\cdot
\Bigl(\psi(j)(X_u,u\in K)-\psi(j)(X_0,X_u,u\in K\setminus\{l\}) \Bigr)\notag\\
&\hspace{2.5cm}-2\mathbbm{1}_{\{J=K\}}\mathbbm{E}\bigl[\psi(i)(X_1,\dotsc,X_{p_i})\psi(j)(X_1,\dotsc,X_{p_j})\bigr]\biggr| \cdot|\psi(k)(X_u,u\in L)|  \Biggr]\notag\\
&\leq \sum_{k=1}^d \frac{\|g\|_M}{p_i\sigma_n(i)\sigma_n(j)\sigma_n(k)}\notag\\
&\hspace{1cm}\cdot\sum_{l=1}^n\sum_{\substack{J\in\mathcal{D}_{p_i}(n),\\
K\in\mathcal{D}_{p_j}(n),\\l\in J\cap K  }}\sum_{\substack{L\in\mathcal{D}_{p_k}(n):\\ L\cap(J\cup K)\not=\emptyset}}|a_J(i)a_K(j)a_L(k)|
\|\psi(i)\|_{L^3(\mu^{p_i})}\|\psi(j)\|_{L^3(\mu^{p_j})}\|\psi(k)\|_{L^3(\mu^{p_k})}\notag\\
&\leq \sum_{k=1}^d \frac{\|g\|_M \|\psi(i)\|_{L^3(\mu^{p_i})}\|\psi(j)\|_{L^3(\mu^{p_j})}\|\psi(k)\|_{L^3(\mu^{p_k})}}{\sigma_n(i)\sigma_n(j)\sigma_n(k)}\sum_{\substack{J\in\mathcal{D}_{p_i}(n),\\K\in\mathcal{D}_{p_j}(n),\\L\in \mathcal{D}_{p_k}(n):\\J\cap K\not=\emptyset,\\L\cap(J\cup K)\not=\emptyset}}|a_J(i)a_K(j)a_L(k)|.
\end{align}
Finally, \eqref{s2a} and \eqref{s2d} imply that
\begin{align*}
\epsilon_2\leq& \sum_{i,j=1}^d\Bigg|\frac{n}{2p_i}\mathbbm{E}D^2f(\mathbf{Y}_n)\left[\left(\mathbf{Y}^{(i)}_n-\left(\mathbf{Y}^{(i)}_n\right)'\right)e_i,\left(\mathbf{Y}^{(j)}_n-\left(\mathbf{Y}^{(j)}_n\right)'\right)e_j\right]\\
&\hspace{8cm}-\mathbbm{E}D^2f(\mathbf{Y}_n)\left[\mathbf{D}_n^{(i)}e_i,\mathbf{D}_n^{(j)}e_j\right]\Bigg|\nonumber\\
\leq & \sum_{i,j,k=1}^d \frac{\|g\|_M \|\psi(i)\|_{L^3(\mu^{p_i})}\|\psi(j)\|_{L^3(\mu^{p_j})}\|\psi(k)\|_{L^3(\mu^{p_k})}}{\sigma_n(i)\sigma_n(j)\sigma_n(k)}\sum_{\substack{J\in\mathcal{D}_{p_i}(n),\\K\in\mathcal{D}_{p_j}(n),\\L\in \mathcal{D}_{p_k}(n):\\J\cap K\not=\emptyset,\\L\cap(J\cup K)\not=\emptyset}}|a_J(i)a_K(j)a_L(k)|.
\end{align*}
\end{proof}

\subsection{Distance from a continuous process}
We now prove the following theorem, which bounds the distance between the law of $\mathbf{Y}_n$ and that of a continuous Gaussian process. Let us introduce some notation first.

Let $\Sigma_n^{(m)}\in\mathbbm{R}^{d\times d}$ be given by
\[\left(\Sigma_n^{(m)}\right)_{i,l}=
\begin{cases}
\frac{n}{\sigma_n{(i)}\sigma_n{(l)}}\underset{m=\max(J)}{\underset{{J\in\mathcal{D}_{p_i}(m):}}{\sum}}a_{J}{(i)}a_J{(l)}\mathbbm{E}\left[\psi(i)(X_1,\dots,X_{p_i})\psi(l)(X_1,\dots,X_{p_l})\right],&\text{if }p_i=p_l\\
0,&\text{otherwise,}
\end{cases}\]
for $i,l=1,\dots,d$. For $i=1,\dots,d$, let
\[\delta_n^{(i)}=\frac{1}{\left(\sigma_n{(i)}\right)^2}\sup_{m\in[n]}\underset{m=\max(J)}{\underset{{J\in\mathcal{D}_{p_i}(m):}}{\sum}}a_J(i)^2\mathbbm{E}\left[\psi(i)^2(X_1,\dots,X_{p_i})\right],\]
where $[n]:=\lbrace 1,\dots, n\rbrace$, and
\[T_n^{(i)}=\frac{1}{\left(\sigma_n{(i)}\right)^2}\sum_{J\in\mathcal{D}_{p_i}(n)}a_J(i)^2\mathbbm{E}\left[\psi(i)^2(X_1,\dots,X_{p_i})\right].\]
Furthermore, let
\[\varphi_n(s)=\sum_{m=p_1}^n\left(\Sigma_n^{(m)}\right)^{1/2}\mathbbm{1}_{\left(\frac{m-1}{n},\frac{m}{n}\right]}(s),\quad s\in[0,1]\]
and suppose that $\varphi:[0,1]\to\mathbbm{R}^{d\times d}$ is a matrix of $L^2([0,1])$-functions such that, for all $i,j=1,\dots,d$,
\[\lim_{n\to\infty}\int_0^1\left|\left(\varphi_n(s)-\varphi(s)\right)_{i,j}\right|^2\, ds=0.\]
Let $\|\cdot\|_F$ denote the Frobenius norm. Suppose that $\mathbf{W}$ is a $d$-dimensional standard Brownian motion. 

Let
\[\mathbf{Z}(t)=\int_0^t\varphi(s)d\mathbf{W}(s)\]
and $\mathbf{Y}_n$ be defined as in Section \ref{intro_weighted}. 
\begin{theorem}\label{theorem_weighted_con}
Under the above setup, for any $g\in M$,
\[\left|\mathbbm{E}g(\mathbf{Y}_n)-\mathbbm{E}g(\mathbf{Z})\right|\leq \|g\|_{M}(\gamma_1+\gamma_2+\gamma_3+\gamma_4+\gamma_5),\]
and, for any $g\in M^0$,
\[\left|\mathbbm{E}g(\mathbf{Y}_n)-\mathbbm{E}g(\mathbf{Z})\right|\leq \|g\|_{M^0}(\gamma_1+\gamma_2+\gamma_3),\]
where
\begin{align*}
\gamma_1&=\frac{2\sqrt{d}}{3p_1}\sum_{i=1}^d\frac{\|\psi(i)\|_{L^3(\mu^{p_i})}^3 }{\sigma_n(i)^3}\sum_{l=1}^n
\left(\sum_{\substack{J\in\mathcal{D}_{p_i}(n):\\ l\in J}} |a_J(i)|\right)^3;\notag\\
\gamma_2&= \sum_{i,j,k=1}^d \frac{\|\psi(i)\|_{L^3(\mu^{p_i})}\|\psi(j)\|_{L^3(\mu^{p_j})}\|\psi(k)\|_{L^3(\mu^{p_k})}}{\sigma_n(i)\sigma_n(j)\sigma_n(k)}\sum_{\substack{J\in\mathcal{D}_{p_i}(n),\\K\in\mathcal{D}_{p_j}(n),\\L\in \mathcal{D}_{p_k}(n):\\J\cap K\not=\emptyset,\\L\cap(J\cup K)\not=\emptyset}}|a_J(i)a_K(j)a_L(k)|;\\
\gamma_3&=2\sqrt{\int_0^1\left\|\varphi_n(s)-\varphi(s)\right\|_F^2ds}+12\sqrt{\sum_{i=1}^d\delta_n^{(i)}\log\left(\frac{2T_n^{(i)}}{\delta_n^{(i)}}\right)};\\
\gamma_4&=\sqrt{d}\sum_{i=1}^d\left[8447\left(\delta_n^{(i)}\log\left(\frac{2T_n^{(i)}}{\delta_n^{(i)}}\right)\right)^{3/2}+44\left(\sum_{j=1}^d\int_0^1\left[\left(\varphi_n(s)-\varphi(s)\right)_{i,j}\right]^2ds\right)^{3/2}\right];\\
\gamma_5&=\sqrt{d}\left(\int_0^1\left\|\varphi(s)\right\|_F^2ds\right)\sum_{i=1}^d\Bigg[50\sqrt{\delta_n^{(i)}\log\left(\frac{2T_n^{(i)}}{\delta_n^{(i)}}\right)}\\
&\hspace{7cm}+19\sqrt{\sum_{j=1}^d\int_0^1\left[\left(\varphi_n(s)-\varphi(s)\right)_{i,j}\right]^2ds}\;\Bigg].
\end{align*}
\end{theorem}
\begin{proof}
Let us write $\mathbf{W}=\left(\mathbf{W}^{(1)},\dots,\mathbf{W}^{(d)}\right)$, where $\mathbf{W}^{(1)},\dots,\mathbf{W}^{(d)}$ are i.i.d. standard Brownian motions in $\mathbbm{R}$.

\textbf{Step 1.} Consider process $\mathbf{D}_n$ defined in Section \ref{pre_lim_weighted}. Note that, for $i=1,\dots,d$,
\begin{align*}
\mathbf{D}_n^{(i)}(t)=&\frac{1}{\sigma_n^{(i)}}\sum_{J\in \mathcal{D}_{p_i}(\lfloor nt\rfloor)}a_J(i)Z_J(i)\\
=&\frac{1}{\sigma_n^{(i)}}\sum_{m=p_i}^{\lfloor nt\rfloor}\underset{m=\max(J)}{\sum_{J\in\mathcal{D}_{p_i}([m]):}}a_J(i)Z_J(i)\\
=&\frac{1}{\sigma_n^{(i)}}\sum_{m=p_i}^{\lfloor nt\rfloor}\tilde{Z}_m(i),
\end{align*}
where $\lbrace \tilde{Z}_{m}(i):m\in [n], i\in [d]\rbrace$ is a jointly Gaussian collection of centred random variables with the following covariance structure:
\begin{align*}
&\mathbbm{E}\left[\tilde{Z}_{m_1}(i)\tilde{Z}_{m_2}(l)\right]\\
=&
\begin{cases}
\underset{m_1=\max(J)}{\underset{{J\in\mathcal{D}_{p_i}([m_1]):}}{\sum}}a_{J}^{(i)}a_J^{(l)}\mathbbm{E}\left[\psi(i)(X_1,\dots,X_{p_i})\psi(l)(X_1,\dots,X_{p_l})\right],&\text{if }p_i=p_l\text{ and }m_1=m_2\\
0,&\text{otherwise}.
\end{cases}
\end{align*}
Using this observation, note that $\mathbf{D}_n$ has the same distribution as $\tilde{\mathbf{Z}}_n$ given by
\begin{align*}
\tilde{\mathbf{Z}}_n(t):=\frac{1}{\sqrt{n}}\sum_{m=p_1}^n\int_0^{\lfloor nt\rfloor}\left(\Sigma_n^{(m)}\right)^{1/2}\mathbbm{1}_{(m-1,m]}(s)d\mathbf{W}(s),\quad t\in[0,1],
\end{align*}
whose distribution, by a simple change of variables, is equal to that of
\begin{align*}
\mathbf{Z}_n(t):=\sum_{m=p_1}^n\int_0^{\lfloor nt\rfloor/n}\left(\Sigma_n^{(m)}\right)^{1/2}\mathbbm{1}_{\left(\frac{m-1}{n},\frac{m}{n}\right]}(s)d\mathbf{W}(s)=\int_0^{\lfloor nt\rfloor/n}\varphi_n(s)d\mathbf{W}(s),\quad t\in[0,1].
\end{align*}

\textbf{Step 2.} By Doob's $L^2$ inequality and It\^o's isometry, we note that
\begin{align}
\mathbbm{E}\sup_{t\in[0,1]}\left|\int_0^{t}\left(\varphi_n(s)-\varphi(s)\right)d\mathbf{W}(s)\right|^2=&\mathbbm{E}\left[\sup_{t\in[0,1]}\sum_{i=1}^d\left(\sum_{j=1}^d\int_0^t\left(\varphi_n(s)-\varphi(s)\right)_{i,j}d\mathbf{W}^{(j)}(s)\right)^2\right]\notag\\
\leq &4\sum_{i=1}^d\mathbbm{E}\left[\left(\sum_{j=1}^d\int_0^1\left(\varphi_n(s)-\varphi(s)\right)_{i,j}d\mathbf{W}^{(j)}(s)\right)^2\right]\notag\\
=&4\sum_{i,j=1}^d\mathbbm{E}\left[\left(\int_0^1\left(\varphi_n(s)-\varphi(s)\right)_{i,j}d\mathbf{W}^{(j)}(s)\right)^2\right]\notag\\
=&4\int_0^1\left\|\varphi_n(s)-\varphi(s)\right\|_F^2ds.\label{second1}
\end{align}
Similarly, by Doob's $L^{3}$ inequality, the formula for Gaussian moments and It\^o's isometry,
\begin{align}
&\mathbbm{E}\sup_{t\in[0,1]}\left|\int_0^{t}\left(\varphi_n(s)-\varphi(s)\right)d\mathbf{W}(s)\right|^3\\
=&\mathbbm{E}\left[\sup_{t\in[0,1]}\left(\sum_{i=1}^d\left(\sum_{j=1}^d\int_0^t\left(\varphi_n(s)-\varphi(s)\right)_{i,j}d\mathbf{W}^{(j)}(s)\right)^2\right)^{3/2}\right]\notag\\
\leq &\frac{27\sqrt{d}}{8}\sum_{i=1}^d\mathbbm{E}\left[\left|\sum_{j=1}^d\int_0^1\left(\varphi_n(s)-\varphi(s)\right)_{i,j}d\mathbf{W}^{(j)}(s)\right|^3\right]\notag\\
=&\frac{27\sqrt{d}}{2\sqrt{2\pi}}\sum_{i=1}^d\left(\mathbbm{E}\left[\left(\sum_{j=1}^d\int_0^1\left(\varphi_n(s)-\varphi(s)\right)_{i,j}d\mathbf{W}^{(j)}(s)\right)^2\right]\right)^{3/2}\notag\\
=&\frac{27\sqrt{d}}{2\sqrt{2\pi}}\sum_{i=1}^d\left(\sum_{j=1}^d\int_0^1\left[\left(\varphi_n(s)-\varphi(s)\right)_{i,j}\right]^2ds\right)^{3/2}.\label{third1}
\end{align}

\textbf{Step 3.} We now apply an argument similar to that of \cite[Theorem 1]{ito_processes}. Note that
\begin{align*}
\mathbf{M}_n(t)=\int_0^{t\wedge 1}\varphi_n(s)d\mathbf{W}(s)+\left(\mathbf{W}(t)-\mathbf{W}(1)\right)\mathbbm{1}_{[t>1]}
\end{align*}
is a martingale vanishing at zero. In particular, so are the coordinate processes
\[\mathbf{M}_n^{(i)}(t)=\int_0^{t\wedge 1} \sum_{j=1}^d\left(\varphi_n\right)_{i,j}d\mathbf{W}^{(j)}(s) +\left(\mathbf{W}^{(i)}(t)-\mathbf{W}^{(i)}(1)\right)\mathbbm{1}_{[t>1]}.\]
Note that, by the Dambis-Dubins-Schwarz theorem, for each $i=1,\dots,d$, there exists a Wiener process $\tilde{\mathbf{W}}^{(i)}$, such that
\begin{align*}
\mathbf{M}_n^{(i)}(t)=\tilde{\mathbf{W}}^{(i)}\left(\left<\mathbf{M}_n^{(i)}\right>_t\right),\quad t\geq 0,
\end{align*}
where $\left<\mathbf{M}_n^{(i)}\right>_t$ is the quadratic variation of $\mathbf{M}_n^{(i)}$, i.e.
\begin{align*}
&\left<\mathbf{M}_n^{(i)}\right>_t=\sum_{j=1}^d\int_0^{t\wedge 1}\left((\varphi_n)_{i,j}\right)^2ds+(t-1)\vee 0.
\end{align*}
Note that
\begin{align*}
\left<\mathbf{M}_n^{(i)}\right>_1=&\sum_{m=p_1}^n\int_0^1\left(\Sigma_n^{(m)}\right)_{i,i}\mathbbm{1}_{\left(\frac{m-1}{n},\frac{m}{n}\right]}(s)ds=\frac{1}{n}\sum_{m=p_1}^n\left(\Sigma_n^{(m)}\right)_{i,i}=T_n^{(i)}
\end{align*}
and
\begin{align*}
\sup_{t\in[0,1]}\left(\left<\mathbf{M}_n^{(i)}\right>_t-\left<\mathbf{M}_n^{(i)}\right>_{\lfloor nt\rfloor/n}\right)=&\sup_{t\in[0,1]}\sum_{j=1}^d\int_{\lfloor nt\rfloor/n}^{t}\left((\varphi_n)_{i,j}(s)\right)^2ds\\
=&\sup_{t\in[0,1]}\sum_{j=1}^d\int_{\lfloor nt\rfloor/n}^{t}\left(\left(\Sigma_n^{\left((\lfloor nt\rfloor +1)\wedge n\right)}\right)^{1/2}\right)_{i,j}^2ds\\
=&\sup_{t\in[0,1]}\left(t-\frac{\lfloor nt\rfloor}{n}\right)\left(\Sigma_n^{\left((\lfloor nt\rfloor +1)\wedge n\right)}\right)_{i,i}\\
\leq&\frac{1}{\left(\sigma_n^{(i)}\right)^2}\sup_{m\in[n]}\underset{m=\max(J)}{\underset{{J\in\mathcal{D}_{p_i}(m):}}{\sum}}a_J(i)^2\mathbbm{E}\left[\psi(i)^2(X_1,\dots,X_{p_i})\right]\\
=&\delta_n^{(i)}.
\end{align*}
Therefore, using \cite[Lemma 3]{ito_processes}, we have that
\begin{align*}
&\mathbbm{E}\sup_{t\in[0,1]}\left|\left(\int_{\lfloor nt\rfloor/n}^t\varphi_n(s)d\mathbf{W}(s)\right)_i\right|^2\\
\leq &\mathbbm{E}\sup\Bigg\{ \left|\tilde{\mathbf{W}}^{(i)}(u)-\tilde{\mathbf{W}}^{(i)}(v)\right|^2:\\
&\hspace{3cm}u,v\in\left[0,\left<\mathbf{M}_n^{(i)}\right>_1\right],\,|u-v|\leq\sup_{t\in[0,1]}\left(\left<\mathbf{M}_n^{(i)}\right>_t-\left<\mathbf{M}_n^{(i)}\right>_{\lfloor nt\rfloor/n}\right)\Bigg\}\\
\leq &\mathbbm{E}\sup\left\lbrace \left|\tilde{\mathbf{W}}^{(i)}(u)-\tilde{\mathbf{W}}^{(i)}(v)\right|^2:\,u,v\in\left[0,T_n^{(i)}\right],\,|u-v|\leq\delta_n^{(i)} \right\rbrace\\
\leq& \frac{5\cdot 6^2}{2\log 2}\left(\delta_n^{(i)}\log\frac{2T_n^{(i)}}{\delta_n^{(i)}}\right)
\end{align*}
and 
\begin{align*}
&\mathbbm{E}\sup_{t\in[0,1]}\left|\left(\int_{\lfloor nt\rfloor/n}^t\varphi_n(s)d\mathbf{W}(s)\right)_i\right|^3\\
\leq &\mathbbm{E}\sup\Bigg\{ \left|\tilde{\mathbf{W}}^{(i)}(u)-\tilde{\mathbf{W}}^{(i)}(v)\right|^3:\\
&\hspace{3cm}u,v\in\left[0,\left<\mathbf{M}_n^{(i)}\right>_1\right],\,|u-v|\leq\sup_{t\in[0,1]}\left(\left<\mathbf{M}_n^{(i)}\right>_t-\left<\mathbf{M}_n^{(i)}\right>_{\lfloor nt\rfloor/n}\right)\Bigg\}\\
\leq &\mathbbm{E}\sup\left\lbrace \left|\tilde{\mathbf{W}}^{(i)}(u)-\tilde{\mathbf{W}}^{(i)}(v)\right|^3:\,u,v\in\left[0,T_n^{(i)}\right],\,|u-v|\leq\delta_n^{(i)} \right\rbrace\\
\leq& \frac{5\cdot 6^3}{\sqrt{\pi}(\log 2)^{3/2}}\left(\delta_n^{(i)}\log\left(\frac{2T_n^{(i)}}{\delta_n^{(i)}}\right)\right)^{3/2}.
\end{align*}
Finally, it follows that
\begin{align}
&\mathbbm{E}\sup_{t\in[0,1]}\left|\int_{\lfloor nt\rfloor/n}^t\varphi_n(s)d\mathbf{W}(s)\right|\leq \frac{6\sqrt{5}}{\sqrt{2\log 2}}\sqrt{\sum_{i=1}^d\delta_n^{(i)}\log\left(\frac{2T_n^{(i)}}{\delta_n^{(i)}}\right)};\label{second2}\\
&\mathbbm{E}\sup_{t\in[0,1]}\left|\int_{\lfloor nt\rfloor/n}^t\varphi_n(s)d\mathbf{W}(s)\right|^3\leq\sqrt{d}\sum_{i=1}^d\mathbbm{E}\sup_{t\in[0,1]}\left|\left(\int_{\lfloor nt\rfloor/n}^t\varphi_n(s)d\mathbf{W}(s)\right)_i\right|^3\notag\\
&\hspace{2cm}\leq  \frac{5\cdot 6^3\sqrt{d}}{\sqrt{\pi}(\log 2)^{3/2}}\sum_{i=1}^d\left(\delta_n^{(i)}\log\left(\frac{2T_n^{(i)}}{\delta_n^{(i)}}\right)\right)^{3/2}.\label{third2}
\end{align}
\textbf{Step 3.}
Using the calculations above, we note that
\begin{align*}
&\mathbbm{E}\|\textbf{Z}_n-\textbf{Z}\|\stackrel{\eqref{second1},\eqref{second2}}\leq 2\sqrt{\int_0^1\left\|\varphi_n(s)-\varphi(s)\right\|_F^2ds}+\frac{6\sqrt{5}}{\sqrt{2\log 2}}\sqrt{\sum_{i=1}^d\delta_n^{(i)}\log\left(\frac{2T_n^{(i)}}{\delta_n^{(i)}}\right)};\\
&\mathbbm{E}\|\textbf{Z}_n-\textbf{Z}\|^3\stackrel{\eqref{third1},\eqref{third2}}\leq\frac{20\cdot 6^3\sqrt{d}}{\sqrt{\pi}(\log 2)^{3/2}}\left(\delta_n^{(i)}\log\left(\frac{2T_n^{(i)}}{\delta_n^{(i)}}\right)\right)^{3/2}\\
&\hspace{6cm}+\frac{54\sqrt{d}}{\sqrt{2\pi}}\sum_{i=1}^d\left(\sum_{j=1}^d\int_0^1\left[\left(\varphi_n(s)-\varphi(s)\right)_{i,j}\right]^2ds\right)^{3/2}.
\end{align*}
We furthermore note that, using Doob's $L^3$ inequality, the formula for Gaussian moments and It\^o's isometry,
\begin{align*}
\mathbbm{E}\|\mathbf{Z}\|^3=&\mathbbm{E}\left[\sup_{t\in[0,1]}\left(\sum_{i=1}^d\left(\sum_{j=1}^d\int_0^t\left(\varphi(s)\right)_{i,j}d\mathbf{W}^{(j)}(s)\right)^2\right)^{3/2}\right]\notag\\
\leq &\frac{27\sqrt{d}}{8}\sum_{i=1}^d\mathbbm{E}\left[\left|\sum_{j=1}^d\int_0^1\left(\varphi(s)\right)_{i,j}d\mathbf{W}^{(j)}(s)\right|^3\right]\notag\\
=&\frac{27\sqrt{d}}{2\sqrt{2\pi}}\sum_{i=1}^d\left(\mathbbm{E}\left[\left(\sum_{j=1}^d\int_0^1\left(\varphi(s)\right)_{i,j}d\mathbf{W}^{(j)}(s)\right)^2\right]\right)^{3/2}\notag\\
=&\frac{27\sqrt{d}}{2\sqrt{2\pi}}\sum_{i=1}^d\left(\sum_{j=1}^d\int_0^1\left|\left(\varphi(s)\right)_{i,j}\right|^2ds\right)^{3/2}.
\end{align*}
Therefore, using the mean value theorem
\begin{align*}
\left|\mathbbm{E}g(\mathbf{D}_n)-\mathbbm{E}g(\mathbf{Z})\right|
\leq& \mathbbm{E}\left[\sup_{c\in[0,1]}\|Dg(\mathbf{Z}+c(\mathbf{Z}_n-\mathbf{Z})\|\|\mathbf{Z}-\mathbf{Z}_n\|\right]\\
\leq&\|g\|_{M}\mathbbm{E}\left[\sup_{c\in[0,1]}\left(1+\|\mathbf{Z}+c(\mathbf{Z}_n-\mathbf{Z})\|^2\right)\|\mathbf{Z}-\mathbf{Z}_n\|\right]\\
\stackrel{\text{H\"older}}\leq&\|g\|_{M}\left\lbrace \mathbbm{E}\|\mathbf{Z}-\mathbf{Z}_n\|+2\mathbbm{E}\|\mathbf{Z}-\mathbf{Z}_n\|^3+2\left(\mathbbm{E}\|\mathbf{Z}\|^3\right)^{2/3}\left(\mathbbm{E}\|\mathbf{Z}-\mathbf{Z}_n\|^3\right)^{1/3}\right\rbrace\\
\leq &\|g\|_M(\gamma_3+\gamma_4+\gamma_5)
\end{align*}
and
\begin{align*}
\left|\mathbbm{E}g(\mathbf{D}_n)-\mathbbm{E}g(\mathbf{Z})\right|\leq& \|g\|_{M^0}\mathbbm{E}\left\|\mathbf{Z}_n-\mathbf{Z}\right\|
\leq \|g\|_{M^0}\gamma_3.
\end{align*}
The result now follows by Theorem \ref{theorem_weighted_pre} and the triangle inequality.
\end{proof}

\begin{remark}\label{remark_weighted}
The approximation results in this Section are merely stated for vectors of \textit{degenerate} weighted $U$-processes. In many applications, however, the given weighted $U$-process might involve non-degenerate kernels. If
\[\mathbf{U}_n(t)=\sum_{J\in\D_p(\lfloor nt\rfloor)} a_J\psi(X_j,j\in J)\]
is such a non-degenerate, weighted $U$-process, then it can be written in its \textit{Hoeffding decompoition} as a sum of degenerate, weighted $U$-processes as follows:
\begin{align*}
\mathbf{U}_n(t)&=\int_{E^p}\psi d\mu^p\sum_{J\in\D_p(\lfloor nt\rfloor)} a_J+\sum_{q=1}^p \sum_{K\in\D_q(\lfloor nt\rfloor)}
\Bigl(\sum_{\substack{J\in \D_p(\lfloor nt\rfloor):\\ K\subseteq J}}a_J \Bigr)\psi_q(X_i,i\in K)\\
&=:\int_{E^p}\psi d\mu^p\sum_{J\in\D_p(\lfloor nt\rfloor)} a_J+\sum_{q=1}^p \mathbf{U}^{(q)}_n(t)\,,
\end{align*}
where the kernels $\psi_q$, $1\leq q\leq p$, are degenerate kernels which are expressible in terms of $\psi$. Hence, the results of this Section for the vector $(\mathbf{U}_n^{(1)},\dotsc,\mathbf{U}_n^{(p)})$ together with the application of a linear functional immediately yield bounds on the approximation of $\mathbf{U}_n$ by a suitable Gaussian process. For simplicity we do not state the resulting bounds explicitly but leave their derivation to the interested reader. In the very particular example of $d$-runs on the line, however, we will work out this procedure in full detail.

\end{remark}

\subsection{Homogeneous sum processes}\label{homsums}
In this subsection we consider an important subclass of weighted, degenerate $U$-processess, namely the processes given as so-called
\textbf{homogeneous sums} or \textbf{homogeneous sum processes}. In this case, the random variables $X_i,i\in\mathbbm{N}$, are real-valued such that $\mathbbm{E}|X_1|^3<\infty$,
$\E[X_1]=0$ and $\E[X_1^2]=1$. Moreover, for each $1\leq i\leq d$, the kernel $\psi(i)$ is given by
\begin{equation*}
\psi(i)(x_1,\dotsc,x_{p_i})=\prod_{j=1}^{p_i} x_j\,.
\end{equation*}
In particular, $\psi(i)$ does not depend on $n$. Hence, for $1\leq i\leq d$ and $t\in[0,1]$ we have that
\[
\mathbf{Y}_n^{(i)}(t)=\frac{1}{\sigma_n(i)}\sum_{J\in\mathcal{D}_{p_i}(\lfloor nt\rfloor)}
 a_J(i)\prod_{j\in J} X_j\,,\]
where the $\sigma_n(i)$ are positive reals and, in this special case, the random variables $Z_J(i)$ making up the processes $\mathbf{D}_n^{(i)}$, defined in Subsection \ref{pre_lim_weighted}, are standard normally distributed.
In this situation we have the following results, which are direct consequences of Theorems \ref{theorem_weighted_pre} and \ref{theorem_weighted_con}, respectively.

\begin{corollary}\label{corhumsums1}
With the above definitions and notation we have that
\begin{align*}
&\Bigl|\mathbbm{E}g(\mathbf{Y_n})-\mathbbm{E}g(\mathbf{D}_n)\Bigr|
\leq\frac{2\sqrt{d}\|g\|_M}{3p_1}\sum_{i=1}^d\frac{\bigl(\E|X_1|^3\bigr)^{p_i}}{\sigma_n(i)^3}\sum_{l=1}^n
\left(\sum_{\substack{J\in\mathcal{D}_{p_i}(n):\\ l\in J}} |a_J(i)|\right)^3\notag\\
&\;+ \|g\|_M \sum_{i,j,k=1}^d \frac{\bigl(\E|X_1|^3\bigr)^{(p_i+p_j+p_k)/3}}{\sigma_n(i)\sigma_n(j)\sigma_n(k)}\sum_{\substack{J\in\mathcal{D}_{p_i}(n),\\K\in\mathcal{D}_{p_j}(n),\\L\in \mathcal{D}_{p_k}(n):\\J\cap K\not=\emptyset,\\L\cap(J\cup K)\not=\emptyset}}|a_J(i)a_K(j)a_L(k)|.
\end{align*}
\end{corollary}

\begin{corollary}\label{corhumsums2}
Let $\Sigma_n^{(m)}\in\mathbbm{R}^{d\times d}$ be given by
\[\left(\Sigma_n^{(m)}\right)_{i,l}=
\begin{cases}
\frac{n}{\sigma_n(i)\sigma_n(l)}\underset{m=\max(J)}{\underset{{J\in\mathcal{D}_{p_i}(m):}}{\sum}}a_{J}(i)a_J(l),&\text{if }p_i=p_l\\
0,&\text{otherwise,}
\end{cases}\]
for $i,l=1,\dots,d$. 

For $i=1,\dots,d$, let
\[\delta_n^{(i)}=\frac{1}{\left(\sigma_n{(i)}\right)^2}\sup_{m\in[n]}\underset{m=\max(J)}{\underset{{J\in\mathcal{D}_{p_i}(m):}}{\sum}}a_J(i)^2,\]
where $[n]=\{1,\dots,n\}$, and
\[T_n^{(i)}=\frac{1}{\left(\sigma_n{(i)}\right)^2}\sum_{J\in\mathcal{D}_{p_i}(n)}a_J(i)^2.\]
Furthermore, let
\[\varphi_n(s)=\sum_{m=p_1}^n\left(\Sigma_n^{(m)}\right)^{1/2}\mathbbm{1}_{\left(\frac{m-1}{n},\frac{m}{n}\right]}(s),\quad s\in[0,1]\]
and suppose that $\varphi:[0,1]\to\mathbbm{R}^{d\times d}$ is matrix of $L^2([0,1])$-functions such that, for any $i,j=1,\dots,d$,
\[\lim_{n\to\infty}\int_0^1\left|\left(\varphi_n(s)-\varphi(s)\right)_{i,j}\right|^2\, ds=0,\]

Let $\mathbf{Y}_n$ be defined as in Section \ref{intro_weighted} and $\|\cdot\|_F$ denote the Frobenius norm. Suppose that $\mathbf{W}$ is a $d$-dimensional standard Brownian motion and
\[\mathbf{Z}(t)=\int_0^t\varphi(s)d\mathbf{W}(s).\]
Then, for any $g\in M$,
\[\left|\mathbbm{E}g(\mathbf{Y}_n)-\mathbbm{E}g(\mathbf{Z})\right|\leq \|g\|_{M}(\gamma_1+\gamma_2+\gamma_3+\gamma_4+\gamma_5)\]
and for any $g\in M^0$,
\[\left|\mathbbm{E}g(\mathbf{Y}_n)-\mathbbm{E}g(\mathbf{Z})\right|\leq \|g\|_{M^0}(\gamma_1+\gamma_2+\gamma_3),\]
where
\begin{align*}
\gamma_1&=\frac{2\sqrt{d}}{3p_1}\sum_{i=1}^d\frac{\left(\mathbbm{E}|X_1|^3\right)^{p_i} }{\sigma_n(i)^3}\sum_{l=1}^n
\left(\sum_{\substack{J\in\mathcal{D}_{p_i}(n):\\ l\in J}} |a_J(i)|\right)^3;\notag\\
\gamma_2&= \sum_{i,j,k=1}^d \frac{\left(\mathbbm{E}|X_1|^3\right)^{(p_i+p_j+p_k)/3}}{\sigma_n(i)\sigma_n(j)\sigma_n(k)}\sum_{\substack{J\in\mathcal{D}_{p_i}(n),\\K\in\mathcal{D}_{p_j}(n),\\L\in \mathcal{D}_{p_k}(n):\\J\cap K\not=\emptyset,\\L\cap(J\cup K)\not=\emptyset}}|a_J(i)a_K(j)a_L(k)|;\\
\gamma_3&=2\sqrt{\int_0^1\left\|\varphi_n(s)-\varphi(s)\right\|_F^2ds}+12\sqrt{\sum_{i=1}^d\delta_n^{(i)}\log\left(\frac{2T_n^{(i)}}{\delta_n^{(i)}}\right)};\\
\gamma_4&=\sqrt{d}\sum_{i=1}^d\left[8447\left(\delta_n^{(i)}\log\left(\frac{2T_n^{(i)}}{\delta_n^{(i)}}\right)\right)^{3/2}+44\left(\sum_{j=1}^d\int_0^1\left[\left(\varphi_n(s)-\varphi(s)\right)_{i,j}\right]^2ds\right)^{3/2}\right];\\
\gamma_5&=\sqrt{d}\Big(\int_0^1\left\|\varphi(s)\right\|_F^2ds\Big)\hspace{-1mm}\sum_{i=1}^d\left[50\sqrt{\delta_n^{(i)}\log\left(\frac{2T_n^{(i)}}{\delta_n^{(i)}}\right)} +19\sqrt{\sum_{j=1}^d\int_0^1\left[\left(\varphi_n(s)-\varphi(s)\right)_{i,j}\right]^2ds}\right]\hspace{-1mm}.
\end{align*}    

\end{corollary}

\begin{remark}~\\
\begin{enumerate}
\item In the case $p=2$ the array $(a_J)_{J\in\D_2(n)}:=(a_J(1))_{J\in\D_2(n)}$ may be identified with the (symmetric) matrix $A=(a_{i,j})_{1\leq i,j\leq n}$, where 
$a_{i,i}=0$ and $a_{i,j}=a_{j,i}$ for all $1\leq i,j\leq n$. Many papers \cite{deJo87, GT, Mik, Rot1, NPR} have established sufficient conditions for the (univariate) CLT to hold for $Y_n:=\mathbf{Y}_n(1)$ in this case (with the choice of $\sigma_n^2=\sum_{1\leq i\neq j\leq n} a_{i,j}^2$). Remarkably, in \cite{NPR} the authors prove a universality principle for homogeneous sums of any order $p\geq 1$. In other words, they find necessary and sufficient conditions on the coefficient functions for the asymptotic normality of $Y_n$ to hold in the case when the $X_j$'s are i.i.d. standard Gaussian. They also show that these conditions imply asymptotic normality of $Y_n$ for any possible choice of the distribution of the $X_j$'s, as long as the $X_j$'s are independent and the usual moment assumptions hold.

Now concentrating on $p=2$ and letting 
\[\lambda_n^*:=\max\{|\lambda|\,:\,\lambda \text{ eigenvalue of } A\}\,,\] 
for the matrix $A$ introduced above, a well-known sufficent condition (see, e.g. \cite[Theorem 1.1]{Mik}) for $Y_n$, $n\in\N$, to be asymptotically normal is that $\lim_{n\to\infty} \lambda_n^*/\sigma_n=0$ (under our standing assumption that $\E|X_1|^3<\infty$). The well-known inequalities (see e.g. \cite{GT}) 
\[\rho_n:=\sqrt{\max_{1\leq i\leq n}\sum_{j:j\not=i}a_{i,j}^2}\leq \lambda_n^*\leq\Gamma_n:=\max_{1\leq i\leq n}\sum_{j:j\not=i}|a_{i,j}|\]
imply that this condition in particular implies the \textbf{Lindeberg type condition} $\lim_{n\to\infty}\rho_n^2/\sigma_n^2=0$, which roughly says that the \textbf{asymptotic influence} of every individual $X_i$ vanishes. On the other hand, it is implied by the stronger (and maybe easier to verify) condition that $\lim_{n\to\infty}\Gamma_n/\sigma_n=0$. We remark that the sufficient condition provided by \cite{NPR} for $d=2$ reduces to $\lim_{n\to\infty}\Tr(A^4)/\sigma_n^4=0$, which is easily seen to be equivalent to $\lim_{n\to\infty} \lambda_n^*/\sigma_n=0$. Here $\Tr(B)=\sum_{i=1}^n b_{i,i}$ denotes the trace of a matrix $B=(b_{i,j})_{1\leq i,j\leq n}$.

From the easy to derive inequality 
\[\sigma_n^{-3}\sum_{i=1}^n\biggl(\sum_{j:j\not=i}|a_{i,j}|\Bigr)^3\geq \Bigl(\sigma_n^{-2}\rho_n^2\Bigr)^{3/2}\]
we conclude that, \textbf{in the univariate case}, the condition $\gamma_1\to 0$ as $n\to\infty$, which follows from our bound in Corollary \ref{corhumsums2}, is also stronger than the Lindeberg condition. The Lindeberg condition is, however, neither necessary (consider e.g. $Y_n=(n-1)^{-1/2} X_1\sum_{j=2}^n X_j$ where the $X_j$ are i.i.d symmetric Rademacher random variables) nor sufficient for the asymptotic normality of the $Y_n$. Hence, by the above inequality, also the sufficient condition $\lim_{n\to\infty} \lambda_n^*/\sigma_n=0$ is not necessary for asymptotic normality to hold.
We now provide upper bounds on the quantities $\gamma_1$ and $\gamma_2$ from our bound in this special case. First note that 
\begin{align*}
&\sigma_n^{-3}\sum_{i=1}^n\biggl(\sum_{j:j\not=i}|a_{i,j}|\Bigr)^3=\sigma_n^{-3}\Biggl(\sum_{i=1}^n\sum_{j,k,l\not=i}|a_{i,j}||a_{i,k}|a_{i,l}|
\Biggr)\\
&=\sum_{i=1}^n\sum_{j:j\not=i}|a_{i,j}|^3 +3\sum_{(i,j,k)\in[n]^3_{\not=}}|a_{i,j}||a_{i,k}|^2+\sum_{(i,j,k,l)\in[n]^4_{\not=}}|a_{i,j}||a_{i,k}|a_{i,l}|\\
&=:\sigma_n^{-3}\bigl(S_1+3S_2+S_3),
\end{align*}
 where $[n]^p_{\not=}$ denotes the collection of all $(i_1,\dotsc,i_p)\in[n]^p$ such that $i_k\not=i_l$ whenever $k\not=l$.
We have 
\begin{align*}
S_1&\leq \left(\max_{k\not=l}|a_{k,l}|\right)\sum_{i=1}^n\sum_{j:j\not=i}|a_{i,j}|^2\leq \rho_n\sigma_n^2\,,\\
S_2&=\sum_{i\not=k}|a_{i,k}|^2\sum_{j:j\not=i,k} |a_{i,j}|\leq \Gamma_n\sum_{1\leq i\not=k\leq n}|a_{i,k}|^2=\Gamma_n\sigma_n^2,\\
S_3&=\sum_{i\not=j}|a_{i,j}|\sum_{k:k\not=i,j}|a_{i,k}|\sum_{l:l\not=i,k,j}|a_{i,l}|\leq\Gamma_n^2\sum_{1\leq i\not=j\leq n}|a_{i,j}|.
\end{align*}
Hence, there is an absolute constant $C_1$ such that 
\[\gamma_1\leq C_1\Biggl(\frac{\rho_n}{\sigma_n}+\frac{\Gamma_n}{\sigma_n}+\frac{\Gamma_n^2}{\sigma_n^2}\frac{\sum_{ i\not=j}|a_{i,j}|}{\sigma_n}\Biggr).\]
The second term $\gamma_2$ in our bound in this case is of the same order as 
\begin{align*}
\sigma_n^{-3}\sum_{\substack{J,K,L\in\mathcal{D}_{2}(n):\\J\cap K\not=\emptyset,\\L\cap(J\cup K)\not=\emptyset}}|a_Ja_Ka_L|
&  \asymp  \sigma_n^{-3}\Biggl(\sum_{i\not=j}|a_{i,j}|^3+\sum_{(i,j,k)\in[n]^3_{\not=}}|a_{i,j}|^2|a_{j,k}|   \\
&\;+ \sum_{(i,j,k)\in[n]^3_{\not=}}|a_{i,j}||a_{j,k}||a_{k,i}|+ \sum_{(i,j,k,l)\in[n]^4_{\not=}}|a_{i,j}||a_{i,k}||a_{k,l}|\Biggr)\,,
\end{align*}
where, for positive sequences, we write $b_n\asymp d_n$ if there are $0<c<C<\infty$ such that $c b_n<d_n < Cb_n$ for all sufficiently large $n$. Note that we have 
\begin{align*}
S_4&:=\sum_{(i,j,k)\in[n]^3_{\not=}}|a_{i,j}||a_{j,k}||a_{k,i}|=\sum_{i\not=j}|a_{i,j}|\sum_{k:k\not=i,j}|a_{j,k}||a_{k,i}|\\
&\leq \sum_{i\not=j}|a_{i,j}|\Bigl(\sum_{k:k\not=i,j}|a_{j,k}|^2\Bigr)^{1/2}\Bigl(\sum_{k:k\not=i,j}|a_{k,i}|^2\Bigr)^{1/2}
\leq \rho_n^2\sum_{i\not=j}|a_{i,j}|\,,\\
S_5&:=\sum_{(i,j,k,l)\in[n]^4_{\not=}}|a_{i,j}||a_{i,k}||a_{k,l}|=\sum_{i\not=j}|a_{i,j}|\sum_{k:k\not=i,j}|a_{i,k}|\sum_{l:l\not=i,j,k}|a_{k,l}|\leq \Gamma_n^2\sum_{i\not=j}|a_{i,j}|\,.
\end{align*}
Thus, there is another absolute constant $C_2$ such that 
\[\gamma_2\leq C_2\Biggl(\frac{\rho_n}{\sigma_n}+\frac{\Gamma_n}{\sigma_n}+\frac{\Gamma_n^2}{\sigma_n^2}\frac{\sum_{ i\not=j}|a_{i,j}|}{\sigma_n}\Biggr).\]
In particular, we obtain the asymptotic normality of $Y_n=\mathbf{Y}_n(1)$ under the assumption that 
\[\Gamma_n=o(\sigma_n) \quad\text{and}\quad \frac{\Gamma_n^2}{\sigma_n^2}=o\Biggl(\frac{\sigma_n}{\sum_{i\not=j}|a_{i,j}|}\Biggr)\,,\]
which is a stronger condition than $\lambda_n^*=o(\sigma_n)$. However, if additionally the terms $\gamma_3, \gamma_4$ and $\gamma_5$ in Corollary \ref{corhumsums2} converge to zero, we can conclude the much stronger result that the whole process $\mathbf{Y}_n$ converges to a continuous Gaussian process on $[0,1]$.

\item The literature around FCLTs for homogeneous sum processes is non-void but nevertheless extremely scarce. Indeed, the only references we have found, whose results might compare to ours (in the one-dimensional case) are \cite{Mik} and \cite{Basa}, of which \cite{Mik} only considers quadratic forms, i.e. the case $p=2$. It turns out that comparing our results to those in \cite{Mik} (for $p=2$) and to those in \cite{Basa} is complicated. Indeed, \cite[Theorem 1.6]{Mik} states the FCLT for the quadratic from $\mathbf{Y}_n$ under the (additional) assumption that $\|\tilde{A}\|^{-2}\|  \tilde{A}^T\tilde{A}\|\to0$ as $n\to\infty$, where $\|\cdot\|$ denotes the Frobenius norm of a matrix 
and where $\tilde{A}=(\tilde{a}_{i,j})_{1\leq i,j\leq n}$ has entries $\tilde{a}_{i,j}=a_{i,j} 1_{\{i>j\}}$. Thus, the matrix $C:=\tilde{A}^T\tilde{A}$ has entries $c_{i,j}=\sum_{k=(i\vee j)+1}^{n} a_{i,k}a_{k,j}$ and, hence, its Frobenius norm is given by a quite complicated expression.

Moreover, we have found that the argument leading to \cite[Theorem 1.1]{Basa} is flawed. Indeed, on page 187 therein, in the display below 
(2.9), one cannot simply drop the quantity $\tau_n^4$ (not even at the price of an enlarged absolute constant $C$) because the claimed inequality must hold for all fixed values of $n\in\N$ (sufficiently large) and $t_1,t_2\in[0,1]$. Moreover, the application of \cite[Theorem 15.6]{Bill} on page 188 seems to be a bit rushed, since the almost sure left-continuity of the limiting Gaussian process $\xi_k$ is not verified. Moreover, the claimed limiting process $\xi_k$ appearing in \cite[Theorem 1.1]{Basa} is not even completely determined, since equation (1.4) theorof only specifies the one-dimensional distributions of $\xi_k$ but not its covariance function. 
\end{enumerate}
\end{remark}

\subsection{Example: runs on the line}\label{runs}
Let $\xi_1,\dots,\xi_n$ be i.i.d. random variables, such that $\mathbbm{P}[\xi_1=1]=p=1-\mathbbm{P}[\xi_1=0]$, for $p\in(0,1)$. For any $1\leq r<n$ let $\sigma_n(r)=\sqrt{np^r(1-p)}$ and $V_r$ be the rescaled centred number of $r$-runs given by
\[\mathbf{V}_n^{(r)}(t):=\frac{1}{\sigma_n(r)}\sum_{m=1}^{\lfloor nt\rfloor}\left(\xi_{m}\cdot\xi_{m+1}\cdot\ldots\cdot\xi_{m+r-1}-p^r\right),\quad t\in[0,1],\]
where we adopt the torus convention, i.e. that $\xi_{n+1}=\xi_1, \xi_{n+2}=\xi_2$ and so on.

A similar setup was considered in \cite{reinert_roellin}, where the authors studied the rate of the (finite-dimensional) weak convergence of the law of $\textbf{V}_n^{(r)}(1)$ to the normal distribution. The authors of \cite{reinert_roellin} note that the standard exchangeable-pair construction of \cite{RiRo97} does not lead to a bound going to zero as $n\to\infty$. In order to solve this problem, they apply their \textit{embedding method} and study the joint convergence of $\left(\textbf{V}_n^{(1)}(1),\dots,\textbf{V}_n^{(r)}(1)\right)$ to a multivariate normal law, using a slightly unusual construction of the exchangeable pair. Our propositions in this subsection provide bounds on the rate of the \textit{functional} convergence of $\left(\mathbf{V}_n^{(r_1)},\dots,\mathbf{V}_n^{(r_d)}\right)$ to a Gaussian process for \textit{any} collection $\lbrace r_1,\dots,r_d\rbrace$. They implicitly use the standard exchangeable-pair construction of Subsection \ref{intro_weighted}. Our bounds are of the same order as the bound on the rate of the (finite-dimensional) convergence provided in \cite{reinert_roellin}.

We start with the following result on the pre-limiting approximation:

\begin{proposition}\label{prop1runs}
Adopt the notation from above. Let $d\geq 1$ and $\frac{n}{2}> r_1\geq r_2\geq\dots\geq r_d\geq 1$. Let
\[\mathbf{V}_n=\left(\mathbf{V}_n^{(r_1)},\dots,\mathbf{V}_n^{(r_d)}\right).\]
Let $\lbrace Z_J: J\in \mathcal{D}_j(n), \,j=1,\dots,r_1\rbrace$ be a collection of i.i.d. standard normal random variables. For $i=1,\dots,d$, let furthermore
\[\mathbf{D}_n^{(r_i)}(t)=\frac{1}{\sigma_n(r_i)}\sum_{m=1}^{\lfloor nt\rfloor}\sum_{j=1}^{r_i}\sum_{0\leq i_1<\dots<i_j\leq r_i-1}p^{r-j}Z_{m+i_1,\dots,m+i_j},\quad t\in[0,1].\]
and
\[\mathbf{D}_n=\left(\mathbf{D}_n^{(r_1)},\dots,\mathbf{D}_n^{(r_d)}\right).\]
Then, for any $g\in M^0$,
\[\left|\mathbbm{E}g(\mathbf{V}_n)-\mathbbm{E}g(\mathbf{D}_n)\right|\leq \|g\|_{M^0}\left(\gamma_1+\gamma_2\right)n^{-1/2},\]
where
\begin{align*}
\gamma_1=&\frac{2\sqrt{dr_1}\left(\sum_{i=1}^dr_i\right)^{3/2}}{3r_d}\sum_{i=1}^d\sum_{j=1}^{r_i}\frac{(1+p^3-2p^4)^jp^{3r_i/2-3j}}{(1-p)^{3/2}}{r_i-1\choose j-1}^3;\\
\gamma_2=&2\sqrt{dr_1}\left(\sum_{i=1}^dr_i\right)\sum_{u,v,w=1}^d\sum_{j_1=1}^{r_u} \sum_{j_2=1}^{r_v}\sum_{j_3=1}^{r_w}\frac{\bigl(1+p^3-2p^4\bigr)^{(j_1+j_2+j_3)/3}p^{(r_u+r_v+r_w)/2-j_1-j_2-j_3}}{(1-p)^{3/2}}\\
&\phantom{....................................................................}\cdot r_w(r_u\vee r_v)^2{r_u-1\choose j_1-1}{r_v-1\choose j_2-1}{r_w-1\choose j_3-1}.
\end{align*}
\end{proposition}

\begin{proof}~\\
 \textbf{Step 1.} For $i=1,2,\dots$, let $X_i=\xi_i-p$. It is easy to prove, by induction on $r$, that
\begin{align}\label{decomposition}
\mathbf{V}_n^{(r)}(t)=&\frac{1}{\sigma_n(r)}\sum_{m=1}^{\lfloor nt\rfloor}\sum_{j=1}^{r}\sum_{0\leq i_1<\dots<i_j\leq r-1}p^{r-j}X_{m+i_1}\dots X_{m+i_j},\quad t\in[0,1]
\end{align}

Indeed, for any $m=1,\dots,n$,
\begin{align*}
\xi_m-p=X_m
\end{align*}
and, assuming that
\begin{align}\label{induc_hyp}
\xi_m\xi_{m+1}\dots\xi_{m+r-1}-p^r=\sum_{j=1}^{r}\sum_{0\leq i_1<\dots<i_j\leq r-1}p^{r-j}X_{m+i_1}\dots X_{m+i_j},
\end{align}
we have
\begin{align*}
&\xi_m\xi_{m+1}\dots\xi_{m+r}-p^{r+1}\\
=&\left(\xi_m\xi_{m+1}\dots\xi_{m+r-1}-p^r\right)\left(\xi_{m+r}-p\right)+p\left(\xi_m\xi_{m+1}\dots\xi_{m+r-1}-p^r\right)+p^r\left(\xi_{m+r}-p\right)\\
\stackrel{(\ref{induc_hyp})}{=}&\sum_{j=1}^{r}\sum_{0\leq i_1<\dots<i_j\leq r-1}p^{r-j}X_{m+i_1}\dots X_{m+i_j}X_{m+r}\\
&+\sum_{j=1}^{r}\sum_{0\leq i_1<\dots<i_j\leq r-1}p^{r+1-j}X_{m+i_1}\dots X_{m+i_j}+p^rX_{m+r}\\
=&\sum_{j=2}^{r+1}\sum_{0\leq i_1<\dots<i_j=r} p^{r+1-j}X_{m+i_1}\dots X_{m+i_j}\\
&+\sum_{j=1}^{r}\sum_{0\leq i_1<\dots<i_j\leq r-1}p^{r+1-j}X_{m+i_1}\dots X_{m+i_j}+p^rX_{m+r}\\
=&\sum_{j=1}^{r+1}\sum_{0\leq i_1<\dots<i_j\leq r}p^{r+1-j}X_{m+i_1}\dots X_{m+i_j},
\end{align*}
as required.

\textbf{Step 2.} Now, for any $r=1,2,\dots,r_1$ and $j=1,\dots,r$, note that
\begin{align*}
&\frac{p^{r-j}}{\sigma_n(r)}\sum_{m=1}^{\lfloor nt\rfloor}\sum_{0\leq i_1<\dots<i_j\leq r-1}X_{m+i_1}\dots X_{m+i_j}\\
=&\frac{p^{r-j}}{\sigma_n(r)}\sum_{m=1}^{\lfloor nt\rfloor}\sum_{m\leq i_1<\dots<i_j\leq m+r-1}X_{i_1}\dots X_{i_j}\\
=&\frac{p^{r-j}}{\sigma_n(r)}\sum_{1\leq i_1<\dots<i_j\leq \lfloor nt\rfloor+r-1}\left(\left(r-i_j+i_1\right)\vee 0\right)X_{i_1}\dots X_{i_j}\\
=&\frac{p^{r-j}}{\sigma_n(r)}\sum_{J\in\mathcal{D}_j \left((\lfloor nt\rfloor+r-1)\wedge n\right)}a_J(r)X_{i_1}\dots X_{i_j},
\end{align*}
for
\begin{align*}
a_J(r):=&p^{r-j}\max\left(r-\max(J)+\min(J), 0\right)\\
&\hspace{-1cm}+p^{r-j}\max\left(r+\min(J\cap(n/2,n])-\max(J\cap[1,n/2))-n, 0\right)\mathbbm{1}_{\lbrace J\cap[1,n/2)\neq\emptyset\neq J\cap(n/2,n]\rbrace}.
\end{align*}
Furthermore, let
\begin{align*}
\mathbf{U}_n^{(r,j)}(t)=\frac{1}{\sigma_n(r)}\sum_{J\in D_{j}(\lfloor nt\rfloor)}a_J(r)\prod_{i\in J}X_i,\quad t\in\left[0,1\right]
\end{align*}
and define function $f:\left(D\left([0,1],\mathbbm{R}\right)\right)^{r_1+\dots+r_d}\to D\left([0,1],\mathbbm{R}^d\right)$, given by
\begin{align*}
&f\left(x_{1,1},\dots,x_{1,r_1},x_{2,2},\dots,x_{2,r_2},\dots,x_{d,1},\dots,x_{d,r_d}\right)\\
=&\left(\left(\sum_{j=1}^{r_1}x_{1,j}\left(\left(t+\frac{r_1-1}{n}\right)\wedge 1\right),\dots,\sum_{j=1}^{r_d}x_{d,j}\left(\left(t+\frac{r_d-1}{n}\right)\wedge 1\right)\right),t\in[0,1]\right).
\end{align*}

Hence, note that, by \eqref{decomposition},
\begin{align*}
g\left(\mathbf{V}_n\right)=g\circ f\left(\mathbf{U}_n^{(r_1,1)},\dots,\mathbf{U}_n^{(r_1,r_1)},\dots,\mathbf{U}_n^{(r_d,1)},\dots,\mathbf{U}_n^{(r_d,r_d)}\right).
\end{align*}
It is proved in Lemma \ref{lem_prop1runs} in Section \ref{sec_prop1runs} of the Appendix that
\begin{align}\label{m0_bound}
\|g\circ f\|_{M^0}\leq \|g\|_{M^0}\sqrt{dr_1}\sum_{i=1}^dr_i.
\end{align}

\textbf{Step 3.}
Now, note that, for $r,r_u,r_v,r_w\in\{1,2,\dots,r_1\}$,
\begin{align}
1)\quad &\sum_{l=1}^n\left(\underset{l\in J}{\sum_{J\in\mathcal{D}_{j}(n):}}|a_J(r)|\right)^3\leq p^{3r-3j}\sum_{l=1}^n\left(\sum_{m=l-r+1}^{l}{r-1\choose j-1}\right)^3=p^{3r-3j}r^3{r-1\choose j-1}^3n\notag\\
2)\quad & \sum_{\substack{J\in\mathcal{D}_{j_1}(n),\\K\in\mathcal{D}_{j_2}(n),\\L\in \mathcal{D}_{j_3}(n):\\J\cap K\not=\emptyset,\\L\cap(J\cup K)\not=\emptyset}}|a_J(r_u)a_K(r_v)a_L(r_w)|\notag\\
\leq&\frac{p^{r_u+r_v+r_w-j_1-j_2-j_3}}{r_u\wedge r_v}\notag\\
&\hspace{1cm}\cdot\sum_{l=1}^n\sum_{m_1=l-r_u+1}^l\sum_{m_2=l-r_v+1}^l\sum_{k=l-r_u\vee r_v+1}^{l+r_u\vee r_v-1}\sum_{m_3=k-r_w+1}^k{r_u-1\choose j_1-1}{r_v-1\choose j_2-1}{r_w-1\choose j_3-1}\notag\\
\leq&2p^{r_u+r_v+r_w-j_1-j_2-j_3}r_w(r_u\vee r_v)^2{r_u-1\choose j_1-1}{r_v-1\choose j_2-1}{r_w-1\choose j_3-1}n
\label{final_bounds}
\end{align}
and so, using \eqref{final_bounds} and \eqref{m0_bound}, for any $g\in M^0$,
\begin{align*}
A)\quad&\|g\circ f\|_{M^0}\frac{2\sqrt{\sum_{i=1}^dr_i}}{3r_d}\sum_{i=1}^d\sum_{j=1}^{r_i}\frac{\bigl(\E|X_1|^3\bigr)^{j}}{\sigma_n(r_i)^3}\sum_{l=1}^n
\left(\sum_{\substack{J\in\mathcal{D}_{j}(n):\\ l\in J}} |a_J(r_i)|\right)^3\\
\leq&\|g\|_{M^0}\frac{2\sqrt{dr_1}\left(\sum_{i=1}^dr_i\right)^{3/2}}{3r_d}\sum_{i=1}^d\sum_{j=1}^{r_i}\frac{(1+p^3-2p^4)^jp^{3r_i/2-3j}}{(1-p)^{3/2}}{r_i-1\choose j-1}^3n^{-1/2};\\
\leq&\|g\|_{M^0}\frac{2\sqrt{dr_1}\left(\sum_{i=1}^dr_i\right)^{3/2}}{3r_d}\sum_{i=1}^d\sum_{j=1}^{r_i}\frac{(1+p^3-2p^4)^jp^{3r_i/2-3j}}{(1-p)^{3/2}}{r_i-1\choose j-1}^3n^{-1/2};\\
B)\quad &\|g\circ f\|_{M^0}\sum_{u,v,w=1}^d\sum_{j_1=1}^{r_u} \sum_{j_2=1}^{r_v}\sum_{j_3=1}^{r_w}\frac{\bigl(\E|X_1|^3\bigr)^{(j_1+j_2+j_3)/3}}{\sigma_n(r_u)\sigma_n(r_v)\sigma_n(r_w)}\sum_{\substack{J\in\mathcal{D}_{j_1}(n),\\K\in\mathcal{D}_{j_2}(n),\\L\in \mathcal{D}_{j_3}(n):\\J\cap K\not=\emptyset,\\L\cap(J\cup K)\not=\emptyset}}|a_J(r_u)a_K(r_v)a_L(r_w)|\\
\leq &2\sqrt{dr_1}\left(\sum_{i=1}^dr_i\right)\sum_{u,v,w=1}^d\sum_{j_1=1}^{r_u} \sum_{j_2=1}^{r_v}\sum_{j_3=1}^{r_w}\frac{\bigl(1+p^3-2p^4\bigr)^{(j_1+j_2+j_3)/3}p^{(r_u+r_v+r_w)/2-j_1-j_2-j_3}}{(1-p)^{3/2}}\\
&\phantom{...........................................................}\cdot r_w(r_u\vee r_v)^2 {r_u-1\choose j_1-1}{r_v-1\choose j_2-1}{r_w-1\choose j_3-1}n^{-1/2}.
\end{align*}
The result now follows by Corollary \ref{corhumsums1}.
\end{proof}

Next, we deal with the continuous process approximation as given in Corollary \ref{corhumsums2}. For this, we need to either compute or estimate the quantities $\delta_n^{(i)}$, $T_n^{(i)}$ and 
$\Sigma_n^{(m)}$. After rearranging the entries of the random vector according to their order as homogeneous sums, we can write $\Sigma_n^{(m)}$ as a block diagonal matrix. More precisely, 
for $1\leq q\leq r_1$ letting 
\begin{align}\label{n}
N(q):=\max\{1\leq j\leq d\,:\, r_j\geq q\}\,,
\end{align}
we can write $\Sigma_n^{(m)}$ as a block diagonal matrix with blocks $\Sigma_n^{(m)}(1),\dotsc,\Sigma_n^{(m)}(r_1)$, where, for fixed $q=1,\dotsc,r_1$, $\Sigma_n^{(m)}(q)$ is an $N(q)\times N(q)$ matrix, namely the 
covariance matrix of the random vector
\[\Bigl(\hspace{-1mm}\sqrt{n}\mathbf{U}_n^{(r_1,q)}(m/n)-\sqrt{n}\mathbf{U}_n^{(r_1,q)}((m-1)/n),\dotsc,\sqrt{n}\mathbf{U}_n^{(r_{N(q)},q)}(m/n)-\sqrt{n}\mathbf{U}_n^{(r_{N(q)},q)}((m-1)/n)\hspace{-1mm}\Bigr)^T\hspace{-1mm}.\]
A simple computation shows that, for $q>1$ and $r_i\wedge r_l\leq m\leq n+1-r_i\wedge r_l$,
\begin{align*}
\Sigma_n^{(m)}(q)(i,l)&=\frac{n}{\sigma_n(r_i)\sigma_n(r_l)}\sum_{\substack{J\in\D_q(n):\\ \max (J)=m}} a_J(r_i)a_J(r_l)\\
&=\frac{p^{\frac{r_i+r_l}{2}-q}}{1-p}
\sum_{k=q-1}^{r_i\wedge r_l -1}{k-1\choose q-2}(r_i-k)(r_l-k).
\end{align*}
Otherwise, for $q>1$ and $m\geq n+2-r_i\wedge r_l$,
\begin{align*}
\Sigma_n^{(m)}(q)(i,l)=&\frac{p^{\frac{r_i+r_l}{2}-q}}{1-p}\Bigg[\sum_{k=q-1}^{r_i\wedge r_l -1}{k-1\choose q-2}(r_i-k)(r_l-k)\\
&\hspace{3cm}+\sum_{u=n+2-r_i\wedge r_l}^m\sum_{k=(q-1)\vee (n-u-1)}^{r_i\wedge r_l-1}{k-1\choose q-2}(r_i-k)(r_l-k)\Bigg].
\end{align*}
Moreover, for $q>1$ and $m\leq r_i\wedge r_l-1$,
\begin{align*}
\Sigma_n^{(m)}(q)(i,l)=\frac{p^{\frac{r_i+r_l}{2}-q}}{1-p}\sum_{k=q-1}^{m-1}{k-1\choose q-2}(r_i-k)(r_l-k),
\end{align*}
and, for all $1\leq m\leq n$
\[\Sigma_n^{(m)}(1)(i,l)=\frac{p^{\frac{r_i+r_l}{2}-1}}{1-p}r_ir_l.\]
Hence, we let $\Sigma$ be a block diagonal matrix with blocks $\Sigma(1)\in\mathbbm{R}^{N(1)\times N(1)},\dots,\Sigma(r_1)\in\mathbbm{R}^{N(r_1)\times N(r_1)}$, where
\begin{align}\label{sigma1}
\Sigma(1)(i,l)=\frac{p^{\frac{r_i+r_l}{2}-1}}{1-p}r_ir_l
\end{align}
 and for any $q=2,\dots,r_1$ and $i,l=1,\dots,N(q)$,
\begin{align}\label{sigma2}
\Sigma(q)(i,l)&=\frac{p^{\frac{r_i+r_l}{2}-q}}{1-p}
\sum_{k=q-1}^{r_i\wedge r_l -1}{k-1\choose q-2}(r_i-k)(r_l-k).
\end{align}
Note that, for $\varphi(s)\equiv\Sigma^{1/2}$ and $\varphi_n(s)=\sum_{m=1}^n\left(\Sigma_n^{(m)}\right)^{1/2}\mathbbm{1}_{\left[\frac{m-1}{n},\frac{m}{n}\right]}(s)$, $s\in[0,1]$,
\begin{align*}
&\int_0^1\|\varphi_n(s)-\varphi(s)\|_F^2ds\\
\leq& \frac{2(r_1-1)}{n}\left[\sum_{m=1}^{r_1-1}\left\|\left(\Sigma_n^{(m)}\right)^{1/2}-\Sigma^{1/2}\right\|_F^2+\sum_{m=n+2-r_1}^{n}\left\|\left(\Sigma_n^{(m)}\right)^{1/2}-\Sigma^{1/2}\right\|_F^2\right]\\
\leq &\frac{4(r_1-1)}{n}\sum_{k=1}^d\sum_{i=1}^{r_k}\left[\sum_{m=1}^{r_1-1}\left(\left|\left(\Sigma_n^{(m)}\right)_{i,i}\right|+\left|\Sigma_{i,i}\right|\right)+\sum_{m=n+2-r_1}^{n}\left(\left|\left(\Sigma_n^{(m)}\right)_{i,i}\right|+\left|\Sigma_{i,i}\right|\right)\right]\\
\leq &\frac{24(r_1)^3}{n}\sum_{q=1}^{r_1}\sum_{i=1}^{N(q)}\sum_{k=q-1}^{r_i-1}\left({k-1\choose q-2}\mathbbm{1}_{[q>1]}+\mathbbm{1}_{[q=1]}\right)\frac{p^{r_i-q}}{1-p}(r_i-k)^2.
\end{align*}
Moreover, with obvious notation,
\begin{align*}
T_n^{(i)}(q)=\frac{1}{\left(\sigma_n{(r_i)}\right)^2}\sum_{J\in\mathcal{D}_{q}(n)}a_J(r_i)^2
=&\frac{1}{n}\sum_{m=1}^n\Sigma_n^{(m)}(q)(i,i)\\
=&\begin{cases}
\frac{p^{r_i-1}}{1-p}r_i^2,\quad&\text{if }q=1,\\
\frac{p^{r_i-q}}{1-p}\sum_{k=q-1}^{r_i -1}{k-1\choose q-2}(r_i-k)^2,\quad &\text{if }q>1.
\end{cases}
\end{align*}
Furthermore, for $q>1$,
\begin{align*}
&\delta_n^{(i)}(q)\\
=&\frac{1}{\left(\sigma_n{(r_i)}\right)^2}\sup_{m\in[n]}\underset{m=\max(J)}{\underset{{J\in\mathcal{D}_{q}(m):}}{\sum}}a_J(i)^2\\
=&\frac{p^{r_i-q}}{n(1-p)}\sum_{k=q-1}^{r_i -1}{k-1\choose q-2}(r_i-k)^2+\frac{p^{r_i-q}}{n(1-p)}\sum_{u=n+2-r_i}^n\sum_{k=(q-1)\vee (n-u-1)}^{r_i-1}{k-1\choose q-2}(r_i-k)^2
\end{align*}
and
\[\delta_n^{(i)}(1)=\frac{p^{r_i-1}}{n(1-p)}r_i^2.\]
Therefore, for all $q=1,\dots,r_1$,
\[\frac{1}{n}T_n^{(i)}(q)\leq\delta_n^{(i)}(q)\leq \frac{r_i}{n}T_n^{(i)}(q).\]
Thus, taking \eqref{m0_bound} into account, we note that
\begin{align*}
&\|g\circ f\|_{M^0}\left[12\sqrt{\sum_{q=1}^{r_1}\sum_{i=1}^{N(q)}\delta_n^{(i)}(q)\log\left(\frac{2T_n^{(i)}(q)}{\delta_n^{(i)}(q)}\right)}+2\sqrt{\int_0^1\|\varphi_n(s)-\varphi(s)\|_F^2ds}\right]\\
\leq&\|g\|_{M^0}\sqrt{dr_1}\sum_{j=1}^dr_j\left[12\frac{\sqrt{\log n}}{\sqrt{n}}\Biggl(\sum_{q=2}^{r_1}\sum_{i=1}^{N(q)}\sum_{k=q-1}^{r_i -1}{k-1\choose q-2}\frac{p^{r_i-q}r_i}{(1-p)}(r_i-k)^2+\sum_{i=1}^d\frac{p^{r_i-1}}{1-p}r_i^3\Biggr)^{1/2}\right.\\
&\left.\phantom{...........}+\frac{4\sqrt{6}(r_1)^{3/2}}{\sqrt{n}}\left(\sum_{q=2}^{r_1}\sum_{i=1}^{N(q)}\sum_{k=q-1}^{r_i-1}{k-1\choose q-2}\frac{p^{r_i-q}}{1-p}(r_i-k)^2+\sum_{i=1}^d\frac{p^{r_i-1}}{1-p}r_i^2\right)^{1/2}\right]\\
\leq&\|g\|_{M^0}4\sqrt{d}r_1^2\left(\sum_{j=1}^dr_j\right)\left(\sum_{q=2}^{r_1}\sum_{i=1}^{N(q)}\sum_{k=q-1}^{r_i-1}{k-1\choose q-2}\frac{p^{r_i-q}}{1-p}(r_i-k)^2+\sum_{i=1}^d\frac{p^{r_i-1}}{1-p}r_i^2\right)^{1/2}\\
&\hspace{11cm}\cdot\frac{3\sqrt{\log n}+\sqrt{6}}{\sqrt{n}}.
\end{align*}
Hence, using Corollary \ref{corhumsums2} and Proposition \ref{prop1runs} (and noting that reordering the arguments of function $f$ does not change the bound on $\|g\circ f\|_{M^0}$ obtained in Lemma \ref{lem_prop1runs}), we obtain the following result:

\begin{proposition}\label{prop2runs}
Adopt the notation form above. In particular, let $N$ be as in \eqref{n},  $\mathbf{V}_n$ be defined as in Proposition \ref{prop1runs} and $\Sigma$ be the block diagonal matrix with blocks $\Sigma(1)\in\mathbbm{R}^{N(1)\times N(1)},\dots,\Sigma(r_1)\in\mathbbm{R}^{N(r_1)\times N(r_1)}$ defined by \eqref{sigma1} and \eqref{sigma2}. Let $\mathbf{Z}'=\Sigma^{1/2}\mathbf{W}$, where $\mathbf{W}$ is a $(\sum_{i=1}^d r_i)$-dimensional standard Brownian motion and write $\mathbf{Z}'=\left(\left(\mathbf{Z}'\right)^{(1)},\left(\mathbf{Z}'\right)^{(2)},\dots\right)$. Set $N(0)=0$.
 For $i=1,\dots,d$ and $t\in[0,1]$, define
\begin{align*}
\mathbf{Z}^{(i)}(t)=&\left(\left(\mathbf{Z}'\right)^{(i)}+\left(\mathbf{Z}'\right)^{(N(1)+i)}+\left(\mathbf{Z}'\right)^{(N(1)+N(2)+i)}+\dots+\left(\mathbf{Z}'\right)^{(N(1)+N(2)+\dots,N(r_i-1)+i)}\right)\\
&\hspace{9cm}\cdot\left(\left(t+\frac{r_i-1}{n}\right)\wedge 1\right)
\end{align*}
and let
\[\mathbf{Z}=\left(\mathbf{Z}^{(1)},\dots,\mathbf{Z}^{(d)}\right).\]
 Then, for any $g\in M^0$, we have 
\[\left|\mathbbm{E}g(\mathbf{V}_n)-\mathbbm{E}g(\mathbf{Z})\right|\leq n^{-1/2}\|g\|_{M^0}\left(\gamma_1+\gamma_2+\gamma_3\sqrt{\log n}\right),\]
where $\gamma_1$ and $\gamma_2$ are as in Proposition \ref{prop1runs} and 
\[\gamma_3=22\sqrt{d}r_1^2\left(\sum_{j=1}^dr_j\right)\left(\sum_{q=2}^{r_1}\sum_{i=1}^{N(q)}\sum_{k=q-1}^{r_i-1}{k-1\choose q-2}\frac{p^{r_i-q}}{1-p}(r_i-k)^2+\sum_{i=1}^d\frac{p^{r_i-1}}{1-p}r_i^2\right)^{1/2}\,.\]
\end{proposition}
\begin{remark}
Assuming that $d, r_1,\dots,r_d$ are all fixed and do not depend on $n$, the bound in Proposition \ref{prop2runs} is of order $\sqrt{\frac{\log n}{n}}$. Therefore, by Proposition \ref{prop_m}, weak convergence of the law of $\mathbf{V}_n$ to that of $\mathbf{Z}$, in both the Skorokhod and the uniform topologies on the Skorokhod space, follows immediately from Proposition \ref{prop2runs} as a corollary.
\end{remark}
\begin{remark}
It is possible to obtain bounds similar to those in Propositions \ref{prop1runs} and \ref{prop2runs} for the larger class of test functions $M$. It would, however, require some more involved computations, which would make the discussion of this example rather long.
\end{remark}

\section{Edge and two-star counts in Erd\H{o}s-Renyi random graphs}\label{section6}
In this section we study an Erd\H{o}s-Renyi random graph with a fixed edge probability $p$ and $\lfloor nt\rfloor$ edges for $t\in[0,1]$. We analyse the asymptotic behaviour of the joint law of its (rescaled) number of edges and its (rescaled) number of two-stars (i.e. subgraphs which are trees with one internal node and $2$ leaves). Hence, we extend the result of \cite{kasprzak3}, where the univariate process convergence of the rescaled number of edges is studied. 
We also extend the analysis of \cite{reinert_roellin1}, whose authors provide a bound on the distance between the (three-dimensional) joint law of the (rescaled) number of edges, two-stars and triangles in a $G(n,p)$ graph and a Gaussian vector. In Theorem \ref{theorem_pre_limiting}, we establish a bound on the distance between our process and a pre-limiting Gaussian processes with paths in $D([0,1],\mathbbm{R}^2)$. Then, in Theorem \ref{theorem_continuous}, a bound on the quality of a continuous Gaussian process approximation is provided.

It is worth noting that the analysis of a three-dimensional process representing the number of edges, triangles and two-stars in a $G(\lfloor nt\rfloor, p)$ graph does not pose any additional challenges except that it makes the algebraic computations more involved. The only reason we do not do it here is that it would make this section rather lengthy.
\subsection{Introduction}\label{section1}
Consider an Erd\H{o}s-Renyi random graph $G(\lfloor nt\rfloor,p)$ on $\lfloor nt\rfloor$ vertices, for $t\in[0,1]$, with a fixed edge probability $p$. Let $I_{i,j}=I_{j,i}$'s be i.i.d. Bernoulli$(p)$ random variables indicating that edge $(i,j)$ is present in this graph. We consider the following process, representing the re-scaled total number of edges
\begin{align}\label{t_n}
\mathbf{T}_n(t)=\frac{\lfloor nt\rfloor-2}{2n^2}\sum_{1\leq i\neq j\leq\lfloor nt\rfloor} I_{i,j}=\frac{\lfloor nt\rfloor-2}{n^2}\sum_{1\leq i<j\leq \lfloor nt\rfloor}I_{i,j},
\end{align}
and a re-scaled statistic related to the number of two-stars
\begin{align}\label{v_n}
\mathbf{V}_n(t)=\frac{1}{6n^2}\underset{i,j,k\text{ distinct}}{\sum_{1\leq i,j,k\leq \lfloor nt\rfloor}}I_{ij}I_{jk}=\frac{1}{n^2}\sum_{1\leq i<j<k\leq \lfloor nt\rfloor}\left(I_{i,j}I_{j,k}+I_{i,j}I_{i,k}+I_{j,k}I_{i,k}\right).
\end{align}
Furthermore, let $\mathbf{Y}_n(t)=\left(\mathbf{T}_n(t)-\mathbbm{E}\mathbf{T}_n(t),\mathbf{V}_n(t)-\mathbbm{E}\mathbf{V}_n(t)\right)$ for $t\in\left[0,1\right]$.
\begin{remark}
Note that, for all $t\in[0,1]$, $\mathbbm{E}\mathbf{T}_n(t)=\frac{\lfloor nt\rfloor -2}{n^2}{\lfloor nt\rfloor \choose 2}p$ and $\mathbbm{E}\mathbf{V}_n(t)=\frac{3}{n^2}{\lfloor nt\rfloor \choose 3}p^2$ and, by an argument similar to that of \cite[Section 5]{reinert_roellin1}, the covariance matrix of $\left(\mathbf{T}_n(t)-\mathbbm{E}\mathbf{T}_n(t),\mathbf{V}_n(t)-\mathbbm{E}\mathbf{V}_n(t)\right)$ is given by
\[3\frac{{\lfloor nt\rfloor \choose 3}}{n^4}p(1-p)\left(\begin{array}{cc}
(\lfloor nt\rfloor -2)&2p(\lfloor nt\rfloor -2)\\
2p(\lfloor nt\rfloor -2)&4p^2(\lfloor nt\rfloor -2)+p(1-p) \end{array}\right).\]
The scaling therefore ensures that the covariances are of the same order in $n$.
\end{remark}
\subsection{Exchangeable pair setup}\label{section2}
In order to construct a suitable exchangeable pair, following \cite{reinert_roellin1}, we pick $(I,J)$ according to $\mathbbm{P}[I=i,J=j]=\frac{1}{{n\choose 2}}$ for $1\leq i<j\leq n$. If $I=i,J=j$, we replace $I_{i,j}=I_{j,i}$ by an independent copy $I_{i,j}'=I_{j,i}'$ and set:
\begin{align*}
\mathbf{T}_n'(t)&=\mathbf{T}_n(t)-\frac{\lfloor nt\rfloor-2}{n^2}\left(I_{I,J}-I_{I,J}'\right)\mathbbm{1}_{[I/n,1]\cap[J/n,1]}(t)\\
\mathbf{V}_n'(t)&=\mathbf{V}_n(t)-\frac{1}{n^2}\sum_{k:k\neq I,J}\left(I_{I,J}-I_{I,J}'\right)\left(I_{J,k}+I_{I,k}\right)\mathbbm{1}_{[I/n,1]\cap[J/n,1]\cap[k/n,1]}(t).
\end{align*}
We, similarly, let $\mathbf{Y}_n'(t)=\left(\mathbf{T}_n'(t)-\mathbbm{E}\mathbf{T}_n(t),\mathbf{V}_n'(t)-\mathbbm{E}\mathbf{V}_n(t)\right)$ and note that, for $\mathbf{Y}_n=\left(\mathbf{Y}_n(t),t\in[0,1]\right)$ and $\mathbf{Y}_n'=\left(\mathbf{Y}_n'(t),t\in[0,1]\right)$, $(\mathbf{Y}_n,\mathbf{Y}_n')$ forms an exchangeable pair.
Note that, for any $m=1,2$, any $f\in M$, as defined in Section \ref{section_notation}, and $e_1,e_2$ denoting the canonical basis vectors $(1,0)$ and $(0,1)$, respectively, we have
\begin{align*}
&\mathbbm{E}^{\mathbf{Y}_n}\left\lbrace Df(\mathbf{Y}_n)\left[\left(\mathbf{T}_n'-\mathbf{T}_n\right)e_m\right]\right\rbrace\\
=&\mathbbm{E}^{\mathbf{Y}_n}\left\lbrace Df(\mathbf{Y}_n)\left[\frac{\lfloor n\cdot\rfloor-2}{n^2}\left(I_{I,J}'-I_{I,J}\right)\mathbbm{1}_{[I/n,1]\cap[J/n,1]}e_m\right]\right\rbrace\\
=&\frac{2}{n^3(n-1)}\sum_{i<j}\mathbbm{E}^{\mathbf{Y}_n}\left\lbrace Df(\mathbf{Y}_n)\left[(\lfloor n\cdot\rfloor-2)\left(I_{i,j}'-I_{i,j}\right)\mathbbm{1}_{[i/n,1]\cap[j/n,1]}e_m\right]|I=i,J=j\right\rbrace\\
=&-\frac{1}{{n\choose 2}}Df(\mathbf{Y}_n)[\mathbf{T}_ne_m]+\frac{2}{n^3(n-1)}p\sum_{i<j}Df(\mathbf{Y}_n)\left[(\lfloor n\cdot\rfloor-2)\mathbbm{1}_{[i/n,1]\cap[j/n,1]}e_m\right]\\
=&-\frac{1}{{n\choose 2}}Df(\mathbf{Y}_n)[\left(\mathbf{T}_n(\cdot)-\mathbbm{E}\mathbf{T}_n(\cdot)\right)e_m].\\
\end{align*}
Also:
\begin{align*}
&\mathbbm{E}^{\mathbf{Y}_n}Df(\mathbf{Y}_n)[(\mathbf{V}_n-\mathbf{V}_n')e_m]\\
=&\frac{1}{n^2{n\choose 2}}\sum_{i<j}\mathbbm{E}^{\mathbf{Y}_n}\Bigg\{\sum_{k:k\neq i,j}Df(\mathbf{Y}_n)\bigg[\left(I_{i,j}-I_{i,j}'\right)\left(I_{j,k}+I_{i,k}\right)\\
&\hspace{7cm}\cdot\mathbbm{1}_{[i/n,1]\cap[j/n,1]\cap[k/n,1]}e_m\bigg]\Bigg|\,I=i,J=j\bigg\}\\
=&\frac{2}{{n\choose 2}}Df(\mathbf{Y}_n)[\mathbf{V}_ne_m]\\
&\hspace{2cm}-\frac{p}{n^2{n\choose 2}}\sum_{i<j}\sum_{k:k\neq i,j}\mathbbm{E}^{\mathbf{Y}_n}Df(\mathbf{Y}_n)\left[\left(I_{j,k}+I_{i,k}\right)\mathbbm{1}_{[i/n,1]\cap[j/n,1]\cap[k/n,1]}e_m\right]\\
=&\frac{2}{{n\choose 2}}Df(\mathbf{Y}_n)[\mathbf{V}_ne_m]-\frac{p}{n^2{n\choose 2}}\underset{i,j,k\text{ distinct}}{\sum_{1\leq i,j,k\leq n}}\mathbbm{E}^{\mathbf{Y}_n}Df(\mathbf{Y}_n)\left[I_{i,j}\mathbbm{1}_{[i/n,1]\cap[j/n,1]\cap[k/n,1]}e_m\right]\\
=&\frac{2}{{n\choose 2}}Df(\mathbf{Y}_n)[\left(\mathbf{V}_n-\mathbbm{E}\mathbf{V}_n(\cdot)\right)e_m]\\
&\hspace{3cm}-\frac{p}{n^2{n\choose 2}}\underset{i,j,k\text{ distinct}}{\sum_{1\leq i,j,k\leq n}}\mathbbm{E}^{\mathbf{Y}_n}Df(\mathbf{Y}_n)\left[(I_{i,j}-p)\mathbbm{1}_{[i/n,1]\cap[j/n,1]\cap[k/n,1]}e_m\right]\\
=&\frac{2}{{n\choose 2}}Df(\mathbf{Y}_n)[\left(\mathbf{V}_n-\mathbbm{E}\mathbf{V}_n(\cdot)\right)e_m]\\
&\hspace{3cm}-\frac{2p}{{n\choose 2}}Df(\mathbf{Y}_n)\left[\frac{1}{\lfloor n\cdot\rfloor-2}\left(\mathbf{T}_n-\mathbbm{E}\mathbf{T}_n(\cdot)\right)e_m\left(\sum_{k=1}^n\mathbbm{1}_{[k/n,1]}-2\right)\right]\\
=&\frac{2}{{n\choose 2}}Df(\mathbf{Y}_n)[\left(\mathbf{V}_n-\mathbbm{E}\mathbf{V}_n(\cdot)\right)e_m]-\frac{2p}{{n\choose 2}}Df(\mathbf{Y}_n)\left[(\mathbf{T}_n-\mathbbm{E}\mathbf{T}_n(\cdot))e_m\right].
\end{align*}
Therefore, for any $m=1,2$:
\begin{align*}
\text{A)}\quad Df(\mathbf{Y}_n)\left[\left(\mathbf{T}_n-\mathbbm{E}\mathbf{T}_n\right)e_m\right]=&\frac{n(n-1)}{2}\mathbbm{E}^{\mathbf{Y}_n}\left\lbrace Df(\mathbf{Y}_n)\left[(\mathbf{T}_n-\mathbf{T}_n')e_m\right]\right\rbrace\\
\text{B)}\quad Df(\mathbf{Y}_n)\left[\left(\mathbf{V}_n-\mathbbm{E}\mathbf{V}_n\right)e_m\right]
=&\frac{n(n-1)}{4}\mathbbm{E}^{\mathbf{Y}_n}\bigg\{ Df(\mathbf{Y}_n)\left[(\mathbf{V}_n-\mathbf{V}_n')e_m\right]\\
&+pDf(\mathbf{Y}_n)\left[\left(\mathbf{T}_n-\mathbbm{E}\mathbf{T}_n\right)e_m\right]\bigg\}\\
=&\frac{n(n-1)}{4}\mathbbm{E}^{\mathbf{Y}_n}\left\lbrace Df(\mathbf{Y}_n)\left[\left(2p(\mathbf{T}_n-\mathbf{T}_n')+\mathbf{V}_n-\mathbf{V}_n'\right)e_m\right]\right\rbrace
\end{align*}
and so
\[Df(\mathbf{Y}_n)[\mathbf{Y}_n]=2\mathbbm{E}^{\mathbf{Y}_n}Df(\mathbf{Y}_n)\left[(\mathbf{Y}_n-\mathbf{Y}_n')\Lambda_n\right],\]
where
\begin{equation}\label{eq_lambda}
\Lambda_n=\frac{n(n-1)}{8}\left(\begin{array}{ccc}
2&2p\\
0&1
\end{array}\right).
\end{equation}
Therefore, condition (\ref{condition}) is satisfied with $\Lambda_n$ of (\ref{eq_lambda}) and $R_f=0$.
\subsection{A pre-limiting process}\label{section3}
Suppose that the collection $\lbrace Z^{(1)}_{i,j}: i,j\in [n], i<j\rbrace\cup\lbrace Z^{(2)}_{i,j,k}: i,j,k\in [n], i<j<k\rbrace$ is jointly centred Gaussian with the following covariance structure:
\begin{align*}
&\mathbbm{E}Z_{ij}^{(1)}Z_{kl}^{(1)}=\begin{cases}
\frac{p(1-p)}{n^4},&\text{if }(i,j)=(k,l),\\
0,&\text{otherwise,}\end{cases}\\
&\mathbbm{E}Z_{i,j,k}^{(2)}Z_{l,m}^{(1)}=\begin{cases}
\frac{2p^2(1-p)}{n^4},&\text{if }\lbrace l,m\rbrace\subset\lbrace i,j,k\rbrace,\\
0,&\text{otherwise,}
\end{cases}\\
&\mathbbm{E}Z_{i,j,k}^{(2)}Z_{r,s,t}^{(2)}=\begin{cases}
\frac{3p^2(1+2p-3p^2)}{n^4},&\text{if }(i,j,k)=(r,s,t),\\
\frac{4p^3(1-p)}{n^4},& \text{if }\left|\lbrace i,j,k\rbrace\cap\lbrace r,s,t\rbrace\right|=2.\\
0,&\text{otherwise.}
\end{cases}
\end{align*}
Let $\mathbf{D}_n=(\mathbf{D}_n^{(1)},\mathbf{D}_n^{(2)})$ be defined in the following way:
\begin{align*}
&\mathbf{D}_n^{(1)}(t)=\left(\lfloor nt\rfloor-2\right)\sum_{1\leq i<j\leq\lfloor nt\rfloor}Z_{i,j}^{(1)},\quad t\in[0,1]\\
&\mathbf{D}_n^{(2)}(t)=\sum_{1\leq i<j<k\leq\lfloor nt\rfloor}Z_{i,j,k}^{(2)},\quad t\in[0,1].
\end{align*}

Note that the covariance structure of the collection $\lbrace Z^{(1)}_{i,j}: i,j\in [n], i<j\rbrace\cup\lbrace Z^{(2)}_{i,j,k}: i,j,k\in [n], i<j<k\rbrace$ is the same as the covariance structure of the summands in the formulas \eqref{t_n} and \eqref{v_n}.
\subsection{Distance from the pre-limiting process}
We provide an estimate of the distance between $\mathbf{Y}_n$ and the pre-limiting piecewise constant Gaussian process.
\begin{theorem}\label{theorem_pre_limiting}
Let $\mathbf{Y}_n$ be defined as in Section \ref{section1} and $\mathbf{D}_n$ be defined as in Section \ref{section3}. Then, for any $g\in M$,
\[\left|\mathbbm{E}g(\mathbf{Y}_n)-\mathbbm{E}g(\mathbf{D}_n)\right|\leq 23\|g\|_{M}n^{-1}.\]
\end{theorem}
\begin{remark}
Our bound in Theorem \ref{theorem_pre_limiting} is of the same order as the analogous bound obtained in \cite{reinert_roellin1} on the distance between the (finite-dimensional) distributions of $\textbf{Y}_n(1)$ and $\textbf{D}_n(1)$.
\end{remark}
The proof is based on Theorem \ref{theorem1}. In \textbf{Step 1} we estimate term $\epsilon_1$, which involves bounding $\|\Lambda_n\|_2$ of (\ref{eq_lambda}) and the third moment of $\|\mathbf{Y}_n-\mathbf{Y}_n'\|$. In \textbf{Step 2} we treat $\epsilon_2$, using involved computations, which are, to a large extent, postponed to the appendix. Term $\epsilon_3$ is equal to zero as $R_f$ of Section \ref{section2} is equal to zero.
\begin{proof}[Proof of Theorem \ref{theorem_pre_limiting}]
We adopt the notation of sections \ref{section1}, \ref{section2}, \ref{section3} and apply Theorem \ref{theorem1}.

\textbf{Step 1.}
 First note that, for $\epsilon_1$ in Theorem \ref{theorem1},
\begin{align*}
|(\mathbf{Y}_n-\mathbf{Y}_n')\Lambda_n|\leq \|\Lambda_n\|_2|\mathbf{Y}_n-\mathbf{Y}_n'|,
\end{align*}
where $|\cdot|$ denotes the Euclidean norm in $\mathbbm{R}^2$ and $\|\cdot\|_2$ is the induced operator $2$-norm. Furthermore,
\begin{align*}
\|\Lambda_n\|_2\leq \|\Lambda_n\|_F=\frac{n(n-1)}{8}\sqrt{2^2+(2p)^2+0^2+1^2}\leq \frac{3n(n-1)}{8},
\end{align*}
for $\|\cdot\|_F$ denoting the Frobenius norm (which, for $\Theta\in\mathbbm{R}^{d_1\times d_2}$ is defined by $\|\Theta\|_F=\sqrt{\sum_{i=1}^{d_1}\sum_{j=1}^{d_2}|\Theta_{i,j}|}$).
Therefore:
\begin{align}
&\mathbbm{E}\left[\|(\mathbf{Y}_n-\mathbf{Y}_n')\Lambda_n\|\|\mathbf{Y}_n-\mathbf{Y}_n'\|^2\right]\nonumber\\
\leq&\frac{3n(n-1)}{8} \mathbbm{E}\|\mathbf{Y}_n-\mathbf{Y}_n'\|^3\nonumber\\
\leq& \frac{3n(n-1)}{8}\mathbbm{E}\left[\frac{(n-2)^2}{n^4}\left(I_{I,J}-I'_{I,J}\right)^2+\frac{1}{n^4}\left(\sum_{k:k\neq I,J}(I_{I,J}-I_{I,J}')\left(I_{J,k}+I_{I,k}\right)\right)^2\right]^{3/2}\nonumber\\
\leq&\frac{3n(n-1)}{8}\left[\frac{(n-2)^2}{n^4}+\frac{\left(2(n-2)\right)^2}{n^4}\right]^{3/2}\nonumber\\
\leq&\frac{5}{n},\nonumber
\end{align}
where the third inequality follows because $|I_{I,J}-I_{I,J}'|\leq 1$ and $|I_{J,k}+I_{I,k}|\leq 2$ for all $k$ and
\begin{align}
\epsilon_1\leq \frac{5\|g\|_{M}}{6n}.\label{4.1.1.1}
\end{align}

\textbf{Step 2.}
In order to deal with $\epsilon_2$ in Theorem \ref{theorem1}, we need to bound
\begin{align}
&\left|\mathbbm{E}D^2f(\mathbf{Y}_n)\left[\left(\mathbf{Y}_n-\mathbf{Y}_n'\right)\Lambda_n,\mathbf{Y}_n-\mathbf{Y}_n'\right]-\mathbbm{E}D^2f(\mathbf{Y}_n)\left[\mathbf{D}_n,\mathbf{D}_n\right]\right|\nonumber\\
=&\left|\frac{n(n-1)}{8}\mathbbm{E}D^2f(\mathbf{Y}_n)\left[\left(2(\mathbf{T}_n-\mathbf{T}_n'),2p(\mathbf{T}_n-\mathbf{T}_n')+(\mathbf{V}_n-\mathbf{V}_n')\right),\left(\mathbf{T}_n-\mathbf{T}_n',\mathbf{V}_n-\mathbf{V}_n'\right)\right]\right.\nonumber\\
&\left.-\mathbbm{E}D^2f(\mathbf{Y}_n)\left[\mathbf{D}_n,\mathbf{D}_n\right]\right|\nonumber\\
\leq&S_1+S_2+S_3+S_4,\label{4_int}
\end{align}
where:
\begin{align}
S_1&=\Bigg|\frac{n(n-1)}{8}\mathbbm{E}D^2f(\mathbf{Y}_n)\left[\left(2(\mathbf{T}_n-\mathbf{T}_n'),0\right),(\mathbf{T}_n-\mathbf{T}_n',0)\right]\notag\\
&\hspace{7cm}-\mathbbm{E}D^2f(\mathbf{Y}_n)\left[\left(\mathbf{D}_n^{(1)},0\right),\left(\mathbf{D}_n^{(1)},0\right)\right]\Bigg|\nonumber\\
S_2&=\left|\frac{n(n-1)}{8}\mathbbm{E}D^2f(\mathbf{Y}_n)\left[\left(0,2p(\mathbf{T}_n-\mathbf{T}_n')+\mathbf{V}_n-\mathbf{V}_n'\right),(\mathbf{T}_n-\mathbf{T}_n',0)\right]\right.\nonumber\\
&\left.\hspace{7cm}-\mathbbm{E}D^2f(\mathbf{Y}_n)\left[\left(0,\mathbf{D}_n^{(2)}\right),\left(\mathbf{D}_n^{(1)},0\right)\right]\right|\nonumber\\
S_3&=\Bigg|\frac{n(n-1)}{8}\mathbbm{E}D^2f(\mathbf{Y}_n)\left[\left(2(\mathbf{T}_n-\mathbf{T}_n'),0\right),(0,\mathbf{V}_n-\mathbf{V}_n')\right]\notag\\
&\hspace{7cm}-\mathbbm{E}D^2f(\mathbf{Y}_n)\left[\left(\mathbf{D}_n^{(1)},0\right),\left(0,\mathbf{D}_n^{(2)}\right)\right]\Bigg|\nonumber\\
S_4&=\left|\frac{n(n-1)}{8}\mathbbm{E}D^2f(\mathbf{Y}_n)\left[\left(0,2p(\mathbf{T}_n-\mathbf{T}_n')+\mathbf{V}_n-\mathbf{V}_n'\right),(0,\mathbf{V}_n-\mathbf{V}_n')\right]\right.\nonumber\\
&\left.\phantom{............................................................................}-\mathbbm{E}D^2f(\mathbf{Y}_n)\left[\left(0,\mathbf{D}_n^{(2)}\right),\left(0,\mathbf{D}_n^{(2)}\right)\right]\right|\label{s1s2}.
\end{align}
In Lemma \ref{lemma8_app}, in the appendix, we obtain the following estimates:
\begin{align}
&S_1\leq \frac{\sqrt{5}\|g\|_{M}}{12n}, \quad S_2\leq\frac{\sqrt{178}\|g\|_{M}}{6n}, \quad S_3\leq \frac{\sqrt{178}\|g\|_{M}}{6n},\quad S_4\leq\frac{(\sqrt{612}+\sqrt{178})\|g\|_{M}}{3n}. \label{s_est}
\end{align}

Note that, therefore, by (\ref{4_int}) and (\ref{s_est}),
\begin{align}
&\epsilon_2=\left|\mathbbm{E}D^2f(\mathbf{Y}_n)\left[(\mathbf{Y}_n-\mathbf{Y}_n')\Lambda_n,\mathbf{Y}_n-\mathbf{Y}_n'\right]-\mathbbm{E}D^2f(\mathbf{Y}_n)\left[\mathbf{D}_n,\mathbf{D}_n\right]\right|
\leq 18\|g\|_{M}n^{-1}.\label{4.8}
\end{align}

Using Theorem \ref{theorem1} together with (\ref{4.8}) and (\ref{4.1.1.1}) gives the desired result.
\end{proof}

\subsection{Distance from the continuous process}
We now study the approximation of $\mathbf{Y}_n$ by a continuous Gaussian process with covariance equal to the limit of the covariance of $\mathbf{D}_n$. We obtain a bound on the quality of this approximation. This is achieved by applying  Theorem \ref{theorem_pre_limiting} and by bounding the distance between $\mathbf{D}_n$ and the continuous process via the Brownian modulus of continuity.
 \begin{theorem}\label{theorem_continuous}
Let $\mathbf{Y}_n$ be defined as in Subsection \ref{section1} and let $\mathbf{Z}=(\mathbf{Z}^{(1)},\mathbf{Z}^{(2)})$ be defined by:
\[\begin{cases}
\mathbf{Z}^{(1)}(t)=\frac{\sqrt{p(1-p)}}{\sqrt{2+8p^2}}t\mathbf{B}_1(t^2)+\frac{p\sqrt{2p(1-p)}}{\sqrt{1+4p^2}}t\mathbf{B}_2(t^2),\\
\mathbf{Z}^{(2)}(t)=\frac{p\sqrt{2p(1-p)}}{\sqrt{1+4p^2}}t\mathbf{B}_1(t^2)+\frac{2p^2\sqrt{2p(1-p)}}{\sqrt{1+4p^2}}t\mathbf{B}_2(t^2)
\end{cases},\]
where $\mathbf{B}_1,\mathbf{B}_2$ are independent standard Brownian Motions. Then, for any $g\in M$:
\[\left|\mathbbm{E}g(\mathbf{Y}_n)-\mathbbm{E}g(\mathbf{Z})\right|\leq \|g\|_{M}\left(16422n^{-1/2}\sqrt{\log n}+138n^{-1/2}\right).\]
\end{theorem}
\begin{remark}
\item Theorem \ref{theorem_continuous}, together with Proposition \ref{prop_m}, implies that $\mathbf{Y}_n$ converges to $\mathbf{Z}$ in distribution with respect to the Skorokhod and uniform topologies.
\end{remark}
\begin{remark}
Theorem \ref{theorem_continuous} can be adapted to situations in which $p=p_n$ varies with $n$. More precisely, as indicated by the necessary and sufficient conditions for approximate normality 
of the marginal distributions given in \cite{Ruc}, Theorem \ref{theorem_continuous} can be modified to yield a quantitative functional CLT in the case that $n^3p_n^2\to\infty$ and $n^2(1-p_n)\to\infty$. 
\end{remark}
In \textbf{Step 1} of the proof of Theorem \ref{theorem_continuous}, we use i.i.d standard Brownian Motions to construct a process $\mathbf{Z}_n$ having the same distribution as $\mathbf{D}_n$. In \textbf{Step 2} we couple $\mathbf{Z}_n$ and $\mathbf{Z}$ and use the Brownian modulus of continuity to bound  moments of the supremum distance between them. In \textbf{Step 3} we combine those bounds with the mean value theorem to obtain the desired final estimate.
\begin{proof}[Proof of Theorem \ref{theorem_continuous}]~\\
\textbf{Step 1.}
Let $\textbf{B}_3$ be another standard Brownian Motion, mutually independent with $\textbf{B}_1$ and $\textbf{B}_2$. Let $\mathbf{Z}_n=\left(\mathbf{Z}_n^{(1)},\mathbf{Z}_n^{(2)}\right)$ be defined by:
\begin{align*}
\text{A)}\quad \mathbf{Z}_n^{(1)}(t)=&\frac{(\lfloor nt\rfloor -2)\sqrt{p(1-p)}}{n^2\sqrt{2+8p^2}}\mathbf{B}_1\left(\lfloor nt\rfloor(\lfloor nt\rfloor -1)\right)\\
&+\frac{(\lfloor nt\rfloor -2)p\sqrt{2p(1-p)}}{n^2\sqrt{1+4p^2}}\mathbf{B}_2\left(\lfloor nt\rfloor(\lfloor nt\rfloor -1)\right);\\
\text{B)}\quad \mathbf{Z}_n^{(2)}(t)=&\frac{(\lfloor nt\rfloor -2)p\sqrt{2p(1-p)}}{n^2\sqrt{1+4p^2}}\mathbf{B}_1\left(\lfloor nt\rfloor(\lfloor nt\rfloor -1)\right)\\
&+\frac{(\lfloor nt\rfloor -2)2p^2\sqrt{2p(1-p)}}{n^2\sqrt{1+4p^2}}\mathbf{B}_2\left(\lfloor nt\rfloor(\lfloor nt\rfloor -1)\right)\\
&+\frac{p(1-p)}{n^{2}\sqrt{2}}\mathbf{B}_3\left(\lfloor nt\rfloor(\lfloor nt\rfloor-1)(\lfloor nt\rfloor-2)\right).
\end{align*}
Now, note that $\left(\mathbf{D}_n^{(1)},\mathbf{D}_n^{(2)}\right)\stackrel{\mathcal{D}}=\left(\mathbf{Z}_n^{(1)},\mathbf{Z}_n^{(2)}\right)$. To see this, observe that for all $u,t\in[0,1]$,
\begin{align}
\text{A)}\quad&\mathbbm{E}\mathbf{D}_n^{(1)}(t)\mathbf{D}_n^{(1)}(u)\notag\\
=&(\lfloor nt\rfloor -2)(\lfloor nu\rfloor -2)\lfloor n(t\wedge u)\rfloor(\lfloor n(t\wedge u)\rfloor-1)\frac{p(1-p)}{2n^4}\nonumber\\
=&\mathbbm{E}\mathbf{Z}_n^{(1)}(t)\mathbf{Z}_n^{(1)}(u)\nonumber;\\
\text{B)}\quad&
\mathbbm{E}\textbf{D}_n^{(2)}(t)\textbf{D}_n^{(2)}(u)\notag\\
=&{\lfloor n(t\wedge u)\rfloor\choose 3}\frac{3p^2(1+2p-3p^2)}{n^4}\notag\\
&+{\lfloor n(t\wedge u)\rfloor\choose 2}\left[(\lfloor nt\rfloor -2)(\lfloor nu\rfloor-2)-\left(\lfloor n(t\wedge u)\rfloor -2\right)\right]\frac{4p^3(1-p)}{n^4}\notag\\
=&\lfloor n(t\wedge u)\rfloor(\lfloor n(t\wedge u)\rfloor-1)\frac{4p^3(1-p)(\lfloor nt\rfloor -2)(\lfloor nu\rfloor-2)+(\lfloor n(t\wedge u)\rfloor-2)p^2(1-p)^2}{2n^4}\notag\\
=&\mathbbm{E}\textbf{Z}^{(2)}(t)\textbf{Z}_{n}^{(2)};\notag\\
\text{C)}\quad&\mathbbm{E}\mathbf{D}_n^{(1)}(t)\mathbf{D}_n^{(2)}(u)\notag\\
=&(\lfloor nt\rfloor -2)(\lfloor nu\rfloor -2)\lfloor n(t\wedge u)\rfloor(\lfloor n(t\wedge u)\rfloor-1)\frac{p^2(1-p)}{n^4}\nonumber\\
=&\mathbbm{E}\mathbf{Z}_n^{(1)}(t)\mathbf{Z}_n^{(2)}(u).\label{cov_structure}
\end{align}

\textbf{Step 2.} We now let $\mathbf{Z}$ be constructed as in Theorem \ref{theorem_continuous}, using the same Brownian Motions $\mathbf{B}_1,\mathbf{B}_2$, as the ones used in the construction of $\mathbf{Z}_n$. In Lemma \ref{lemma10_app}, proved in the appendix, we obtain the following bounds:
\begin{align}
&\mathbbm{E}\left\|\mathbf{Z}_n-\mathbf{Z}\right\|\leq\frac{8}{\sqrt{n}}+\frac{39\sqrt{\log n}}{\sqrt{n}}\notag\\
&\mathbbm{E}\left\|\mathbf{Z}_n-\mathbf{Z}\right\|^3\leq\frac{49}{n^{3/2}}+\frac{8167(\log n)^{3/2}}{n^{3/2}}\notag\\
&\mathbbm{E}\|\mathbf{Z}\|^2\leq \frac{4}{3}.\label{4.1.2.1}
\end{align}

\textbf{Step 3.} We note that, by (\ref{4.1.2.1}):
\begin{align*}
\left|\mathbbm{E}g(\mathbf{Z})-\mathbbm{E}g(\mathbf{D}_n)\right|\stackrel{\text{MVT}}\leq&\mathbbm{E}\left[\sup_{c\in[0,1]}\left\|Dg(\mathbf{Z}+c(\mathbf{Z}_n-\mathbf{Z}))\right\|\|\mathbf{Z}-\mathbf{Z}_n\|\right]\\
\leq&\|g\|_{M}\mathbbm{E}\left[\sup_{c\in[0,1]}\left(1+\|\mathbf{Z}+c(\mathbf{Z}_n-\mathbf{Z})\|^2\right)\|\mathbf{Z}-\mathbf{Z}_n\|\right]\\
\leq&\|g\|_{M}\mathbbm{E}\left[\|\mathbf{Z}-\mathbf{Z}_n\|+\|\mathbf{Z}\|\|\mathbf{Z}-\mathbf{Z}_n\|+\|\mathbf{Z}-\mathbf{Z}_n\|^2\right]\\
\leq&\|g\|_{M}\left[\mathbbm{E}\|\mathbf{Z}-\mathbf{Z}_n\|+2\mathbbm{E}\|\mathbf{Z}-\mathbf{Z}_n\|^3+2\left(\mathbbm{E}\|\mathbf{Z}\|^3\right)^{2/3}\left(\mathbbm{E}\|\mathbf{Z}-\mathbf{Z}_n\|^3\right)^{1/3}\right]\\
\leq&\|g\|_{M}\left(\frac{115}{\sqrt{n}}+\frac{16422\sqrt{\log n}}{\sqrt{n}}\right),
\end{align*}
which, together with Theorem \ref{theorem_pre_limiting} gives the desired estimate.
\end{proof}
\begin{remark}
The representation of $\mathbf{Z}$ in terms of two independent Brownian Motions comes from a careful analysis of the limiting covariance of $\mathbf{D}_n$, which may be derived using (\ref{cov_structure}).
\end{remark}

\section{Appendix - technical details of the proofs of Proposition \ref{prop1runs} and Theorems  \ref{theorem_pre_limiting} and \ref{theorem_continuous}}\label{appendix}

\subsection{Technical details of the proof of Proposition \ref{prop1runs} }\label{sec_prop1runs}
\begin{lemma}\label{lem_prop1runs}
Let $n,d\in\mathbbm{N}$ and $r_1\geq r_2\geq\dots\geq r_d\geq 1$.
Define function \\ $f:\left(D\left([0,1],\mathbbm{R}\right)\right)^{r_1+\dots+r_d}\to D\left([0,1],\mathbbm{R}^d\right)$, given by
\begin{align*}
&f\left(x_{1,1},\dots,x_{1,r_1},x_{2,2},\dots,x_{2,r_2},\dots,x_{d,1},\dots,x_{d,r_d}\right)\\
=&\left(\left(\sum_{j=1}^{r_1}x_{1,j}\left(\left(t+\frac{r_1-1}{n}\right)\wedge 1\right),\dots,\sum_{j=1}^{r_d}x_{d,j}\left(\left(t+\frac{r_d-1}{n}\right)\wedge 1\right)\right),t\in[0,1]\right).
\end{align*}
Then , for any $g\in M^0$,
\[\|g\circ f\|_{M^0}\leq \|g\|_{M^0}\sqrt{dr_1}\sum_{i=1}^dr_i.\]
\end{lemma}
\begin{proof}
Note that function $f$ is twice Fr\'echet differentiable with
\begin{align*}
(A)\quad &Df(w)\left[\left(x_{1,1},\dots,x_{1,r_1},x_{2,1},\dots,x_{2,r_2},\dots,x_{d,1},\dots,x_{d,r_d}\right)\right]\\
=&\left(\left(\sum_{j=1}^{r_1}x_{1,j}\left(\left(t+\frac{r_1-1}{n}\right)\wedge 1\right),\dots,\sum_{j=1}^{r_d}x_{d,j}\left(\left(t+\frac{r_d-1}{n}\right)\wedge 1\right)\right),t\in[0,1]\right)\\
(B)\quad &D^2f(w)[x^{(1)},x^{(2)}]=0
\end{align*}
for all $w,x^{(1)},x^{(2)},\left(x_{1,1},\dots,x_{1,r_1},x_{2,1},\dots,x_{2,r_2},\dots,x_{d,1},\dots,x_{d,r_d}\right)\in \left(D\left([0,1],\mathbbm{R}\right)\right)^{r_1+\dots+r_d}$.
Furthermore, for any $w\in\left(D\left([0,1],\mathbbm{R}\right)\right)^{r_1+\dots+r_d}$,
\begin{align*}
a)&\quad \|f(w)\|\\
\leq&\sqrt{\sup_{t\in[0,1]}\left(\sum_{j=1}^{r_1}w_{1,j}\left(\left(t+\frac{r_1-1}{n}\right)\wedge 1\right)\right)^2+\dots+\sup_{t\in[0,1]}\left(\sum_{j=1}^{r_d}w_{d,j}\left(\left(t+\frac{r_d-1}{n}\right)\wedge 1\right)\right)^2}\\
\leq &\sqrt{\sum_{i=1}^d\sup_{t\in[0,1]}\left|\sum_{j=1}^{r_i}w_{i,j}(t)\right|^2}\\
b)&\quad \|Df(w)\|\leq \sqrt{\sum_{i=1}^dr_i}. 
\end{align*}
Therefore, for any $w,h\in\left(D\left([0,1],\mathbbm{R}\right)\right)^{r_1+\dots+r_d}$,

\begin{align}
A)\quad &|g\circ f(w)|\leq\|g\|_{M^0};\notag\\
B)\quad&\left\|D(g\circ f)(w)\right\|=\left\|Dg(f(w))[Df(w)[\cdot]]\right\|\leq \|g\|_{M^0}\|Df(w)\|
\leq\|g\|_{M^0}\sqrt{\sum_{i=1}^dr_i};\notag\\
C)\quad&\left\|D^2(g\circ f)(w)\right\|=\left\|D^2g(f(w))\left[Df(w),Df(w)\right]\right\|\leq \|g\|_{M^0}\left\|Df(w)\right\|^2\leq  \|g\|_{M^0}\sum_{i=1}^dr_i;\notag\\
D)\quad&\left\|D^2(g\circ f)(w+h)-D^2(g\circ f)(w)\right\|\notag\\
=&\left\|D^2g(f(w+h))\left[Df(w+h),Df(w+h)\right]-D^2g(f(w))[Df(w),Df(w)]\right\|\notag\\
\leq &\left\|D^2g(f(w+h))\left[Df(w+h),Df(w+h)\right]-D^2g(f(w))\left[Df(w+h),Df(w+h)\right]\right\|\notag\\
&+\left\|D^2g(f(w))\left[Df(w+h),Df(w+h)\right]-D^2g(f(w))[Df(w),Df(w)]\right\|\notag\\
\leq &\|g\|_{M^0}\|f(w+h)-f(w)\|\|Df(w+h)\|^2\notag\\
\leq &\|g\|_{M^0}\left(\sum_{i=1}^d\sup_{t\in[0,1]}\left|\sum_{j=1}^{r_i}h_{i,j}(t)\right|^2\right)^{1/2}\sum_{i=1}^d r_i,\label{m01}
\end{align}
where $D)$ follows from the fact that $Df(w)=Df(w+h)$. Moreover,
\begin{align}
\frac{\left(\sum_{i=1}^d\sup_{t\in[0,1]}\left|\sum_{j=1}^{r_i}h_{i,j}(t)\right|^2\right)^{1/2}}{\sup_{t\in[0,1]}\left(\sum_{i=1}^d\sum_{j=1}^{r_i}h_{i,j}^2(t)\right)^{1/2}}\leq \frac{\sup_{t\in[0,1]}\left(dr_1\sum_{i=1}^d\sum_{j=1}^{r_i}h_{i,j}^2(t)\right)^{1/2}}{\sup_{t\in[0,1]}\left(\sum_{i=1}^d\sum_{j=1}^{r_i}h_{i,j}^2(t)\right)^{1/2}}=\sqrt{dr_1}.\label{m02}
\end{align}

Therefore,  using \eqref{m01} and \eqref{m02},
\begin{align*}
\|g\circ f\|_{M^0}\leq \|g\|_{M^0}\sqrt{dr_1}\sum_{i=1}^dr_i.
\end{align*}
\end{proof}

\subsection{Technical details of the proof of Theorem \ref{theorem_pre_limiting}}
\begin{lemma}\label{lemma8_app}
For $S_i,i=1,2,3,4$ of (\ref{s1s2}), we have the following estimates:
\begin{align*}
&S_1\leq \frac{\sqrt{5}\|g\|_{M}}{12n}, \quad S_2\leq\frac{\sqrt{178}\|g\|_{M}}{6n}, \quad S_3\leq \frac{\sqrt{178}\|g\|_{M}}{6n},\quad S_4\leq\frac{(\sqrt{612}+\sqrt{178})\|g\|_{M}}{3n}.
\end{align*}
\end{lemma}

\begin{proof}
For $S_1$, for fixed $i,j\in\lbrace 1,\cdots, n\rbrace$, let $\mathbf{Y}_n^{ij}$ be equal to $\mathbf{Y}_n$ except for the fact that $I_{ij}$ is replaced by an independent copy, i.e. for all $t\in [0,1]$ let:
\begin{align*}
\mathbf{T}_n^{ij}(t)&=\mathbf{T}_n(t)-\frac{\lfloor nt\rfloor-2}{n^2}\left(I_{ij}-I_{ij}'\right)\mathbbm{1}_{[i/n,1]\cap[j/n,1]}(t)\\
\mathbf{V}_n^{ij}(t)&=\mathbf{V}_n(t)-\frac{1}{n^2}\sum_{k:k\neq i,j}\left(I_{ij}-I_{ij}'\right)\left(I_{jk}+I_{ik}\right)\mathbbm{1}_{[i/n,1]\cap[j/n,1]\cap[k/n,1]}(t)
\end{align*}
and let $\mathbf{Y}_n^{ij}(t)=\left(\mathbf{T}_n^{ij}(t)-\mathbbm{E}\mathbf{T}_n(t),\mathbf{V}_n^{ij}(t)-\mathbbm{E}\mathbf{V}_n(t)\right)$.

By noting that the mean zero $Z_{i,k}^{(1)}$ and $Z_{i',j}^{(1)}$ are independent for $i\neq i'$, we obtain:
\begin{align}
S_1=&\Bigg|\vphantom{\sum_1^1}\frac{n(n-1)}{8}\mathbbm{E}D^2f(\mathbf{Y}_n)\left[(\mathbf{T}_n-\mathbf{T}_n')(2,0),(\mathbf{T}_n-\mathbf{T}_n')(1,0)\right]\nonumber\\
&-\mathbbm{E}D^2f(\mathbf{Y}_n)\bigg[\sum_{1\leq i< j\leq n}Z_{i,j}^{(1)}(\lfloor n\cdot \rfloor -2)(1,0)\mathbbm{1}_{[i/n,1]\cap[j/n,1]},\notag\\
&\hspace{5cm}\sum_{1\leq i< j\leq n}Z_{i,j}^{(1)}(\lfloor n\cdot \rfloor -2)(1,0)\mathbbm{1}_{[i/n,1]\cap[j/n,1]}\bigg]\Bigg|\nonumber\\
=&\left|\frac{1}{2n^4}\sum_{1\leq i< j\leq n}\mathbbm{E}\bigg\{\left(I_{i,j}-2pI_{i,j}+p\right)\right.\nonumber\\
&\hspace{1cm}\cdot D^2f(\mathbf{Y}_n)\left[(\lfloor n\cdot \rfloor -2)(1,0)\mathbbm{1}_{[i/n,1]\cap[j/n,1]},(\lfloor n\cdot \rfloor -2)(1,0)\mathbbm{1}_{[i/n,1]\cap[j/n,1]}\right]\bigg\}\nonumber\\
&-\sum_{1\leq i<j\leq n}\Bigg\{\mathbbm{E}\left(Z_{i,j}^{(1)}\right)^2\nonumber\\
&\left.\cdot\mathbbm{E}D^2f(\mathbf{Y}_n)\left[(\lfloor n\cdot \rfloor -2)(1,0)\mathbbm{1}_{[i/n,1]\cap[j/n,1]},(\lfloor n\cdot \rfloor -2)(1,0)\mathbbm{1}_{[i/n,1]\cap[j/n,1]}\right]\vphantom{\left(Z^{1}_2\right)^2}\bigg\}\right|\nonumber\\
=&\left|\sum_{1\leq i< j\leq n}\mathbbm{E}\left\lbrace \left(\frac{1}{2n^4}(I_{i,j}-2pI_{i,j}+p)-\mathbbm{E}\left(Z_{i,j}^{(1)}\right)^2\right)\right.\right.\nonumber\\
&\left.\left.\hphantom{\sum_{1\leq i\neq j\leq n}\mathbbm{E}}\cdot D^2f(\mathbf{Y}_n)\left[(\lfloor n\cdot \rfloor -2)(1,0)\mathbbm{1}_{[i/n,1]\cap[j/n,1]},(\lfloor n\cdot \rfloor -2)(1,0)\mathbbm{1}_{[i/n,1]\cap[j/n,1]}\right]\vphantom{\frac{1}{4n^2}}\right\rbrace\vphantom{\sum_1^2}\right|\nonumber\\
=&\Bigg|\sum_{1\leq i< j\leq n}\mathbbm{E}\bigg\{ \frac{1}{2n^4}(I_{i,j}-2pI_{i,j}+p)\left(D^2f(\mathbf{Y}_n)-D^2f(\mathbf{Y}_n^{ij})\right)\nonumber\\
&\hspace{3.5cm}\Big[(\lfloor n\cdot \rfloor -2)(1,0)\mathbbm{1}_{[i/n,1]\cap[j/n,1]},(\lfloor n\cdot \rfloor -2)(1,0)\mathbbm{1}_{[i/n,1]\cap[j/n,1]}\Big]\vphantom{\frac{1}{n^2}}\bigg\}\Bigg|\nonumber\\
\stackrel{\eqref{m_bound}}\leq&\frac{\|g\|_{M}}{6n^2}\sum_{1\leq i< j\leq n}\mathbbm{E}\left|(I_{i,j}-2pI_{i,j}+p)\right|\left\|\mathbf{Y}_n-\mathbf{Y}_n^{ij}\right\|\label{eq_star},
\end{align}
where (\ref{eq_star}) follows from Proposition \ref{prop12.7}. Now,
\[\left\|\mathbf{Y}_n-\mathbf{Y}_n^{ij}\right\|\leq\frac{1}{n^2}\sqrt{(\lfloor n\cdot\rfloor -2)^2(I_{ij}-I_{ij}')^2+\left(\sum_{k:k\neq i,j}|I_{ij}-I_{ij}'|(I_{jk}+I_{ik})\right)^2}\]
and so, by (\ref{eq_star}),
\begin{align}
S_1\leq&\frac{\|g\|_{M}}{6n^4}\sum_{1\leq i< j\leq n}\mathbbm{E}\Bigg\{\left|I_{i,j}-2pI_{i,j}+p\right|\notag\\
&\hspace{4cm}\cdot\sqrt{(n-2)^2(I_{ij}-I_{ij}')^2+\bigg(\sum_{k\neq i,j}|I_{ij}-I_{ij}'|(I_{jk}+I_{ik})\bigg)^2}\Bigg\}\nonumber\\
\leq&\frac{\|g\|_{M}}{6n^3}\sum_{1\leq i<j\leq n}\mathbbm{E}\left\lbrace\left|I_{i,j}-2pI_{i,j}+p\right|\cdot\sqrt{(I_{ij}-I_{ij}')^2+\left(|I_{ij}-I_{ij}'|(I_{jk}+I_{ik})\right)^2}\right\rbrace\nonumber\\
\leq&\frac{\sqrt{5}\|g\|_{M}}{12n},\label{4.1}
\end{align}
where the last inequality holds because $|I_{ij}-2pI_{ij}+p|\leq 1$, $|I_{ij}-I_{ij}'|\leq 1$ and $I_{jk}+I_{ik}\leq 2$ for all $k\in\lbrace 1,\cdots,n\rbrace$.

For $S_2$, let $\mathbf{Y}_n^{ijk}$ equal to $\mathbf{Y}_n$ except that $I_{ij},I_{jk},I_{ik}$ are replaced by $I_{ij}'$, $I_{jk}'$, $I_{ik}'$, i.e. for all $t\in[0,1]$ let
\begin{align}
\mathbf{T}_n^{ijk}(t)=&\mathbf{T}_n(t)-\frac{\lfloor nt\rfloor -2}{n^2}\left[(I_{ij}-I_{ij}')\mathbbm{1}_{[i/n,1]\cap[j/n,1]}(t)\right.\nonumber\\
&\left.+(I_{jk}-I_{jk}')\mathbbm{1}_{[j/n,1]\cap[k/n,1]}(t)+(I_{ik}-I_{ik}')\mathbbm{1}_{[i/n,1]\cap[k/n,1]}(t)\right]\nonumber\\
\mathbf{V}_n^{ijk}(t)=&\mathbf{V}_n(t)-\frac{1}{n^2}\sum_{l:l\neq i,j,k}\bigg[\left(I_{ij}-I_{ij}'\right)\left(I_{jl}+I_{il}\right)\mathbbm{1}_{[i/n,1]\cap[j/n,1]\cap[l/n,1]}(t)\nonumber\\
&\hspace{5cm}+\left(I_{jk}-I_{jk}'\right)\left(I_{jl}+I_{kl}\right)\mathbbm{1}_{[k/n,1]\cap[j/n,1]\cap[l/n,1]}(t)\notag\\
&\hspace{5cm}+\left(I_{ik}-I_{ik}'\right)\left(I_{jl}+I_{il}\right)\mathbbm{1}_{[i/n,1]\cap[k/n,1]\cap[l/n,1]}(t)\bigg]\nonumber\\
&\hspace{-1cm}-\frac{1}{n^2}\left[(I_{ij}I_{jk}-I_{ij}'I_{jk}')+(I_{ij}I_{ik}-I_{ij}'I_{ik}')+(I_{ik}I_{jk}-I_{ik}'I_{jk}')\right]\mathbbm{1}_{[i/n,1]\cap[j/n,1]\cap[k/n,1]}(t).\label{eq_ijk}
\end{align}
Let $\mathbf{Y}_n^{ijk}(t)=\left(\mathbf{T}_n^{ijk}(t)-\mathbbm{E}\mathbf{T}_n(t),\mathbf{V}_n^{ijk}(t)-\mathbbm{E}\mathbf{V}_n(t)\right)$ for all $t\in[0,1]$. Then
\begin{align}
S_2=&\Bigg|\vphantom{\sum_{1\leq i<j<k\leq n}}\frac{n(n-1)}{8}\mathbbm{E}\left\lbrace D^2f(\mathbf{Y}_n)\left[(\mathbf{T}_n-\mathbf{T}_n')(0,2p)+(\mathbf{V}_n-\mathbf{V}_n')(0,1),(\mathbf{T}_n-\mathbf{T}_n')(1,0)\right]\right\rbrace\nonumber\\
&-\mathbbm{E}D^2f(\mathbf{Y}_n)\bigg[\sum_{1\leq i<j<k\leq n}Z_{i,j,k}^{(2)}(0,1)\mathbbm{1}_{[i/n,1]\cap[j/n,1]\cap[k/n,1]},\notag\\
&\hspace{6cm}\sum_{1\leq i< j\leq n}Z_{i,j}^{(1)}(\lfloor n\cdot \rfloor -2)(1,0)\mathbbm{1}_{[i/n,1]\cap[j/n,1]}\bigg]\Bigg|\nonumber\\
=&\left|\frac{1}{4n^4}\hspace{-1mm}\sum_{1\leq i<j\leq n}\underset{k\not\in\lbrace i,j\rbrace}{\sum_{1\leq k\leq n}}\mathbbm{E}\left\lbrace \left[2p\left(I_{ij}-2pI_{ij}+p\right)+\left(I_{ij}-2pI_{ij}+p\right)(I_{jk}+I_{ik})-8p^2(1-p)\right]\right.\right.\notag\\
&\left.\left.\cdot \left(D^2f(\mathbf{Y}_n)-D^2f(\mathbf{Y}^{ijk}_n)\right)\left[(0,1)\mathbbm{1}_{[i/n,1]\cap[j/n,1]\cap[k/n,1]},(\lfloor n\cdot\rfloor-2)(1,0)\mathbbm{1}_{[i/n,1]\cap[j/n,1]}\right]\right\rbrace\vphantom{\underset{j}{\sum_{1\leq i<j<k\leq n}}}\right|\nonumber\\
\stackrel{(\ref{m_bound})}\leq &\frac{\|g\|_{M}}{12n^3}\sum_{1\leq i<j\leq n}\underset{k\not\in\lbrace i,j\rbrace}{\sum_{1\leq k\leq n}}\mathbbm{E}\bigg\{\Big|2p\left(I_{ij}-2pI_{ij}+p\right)\notag\\
&\hspace{5cm}+\left(I_{ij}-2pI_{ij}+p\right)(I_{jk}+I_{ik})\Big|\|\mathbf{Y}_n-\mathbf{Y}_n^{ijk}\|\bigg\}\notag\\
\leq &\frac{\|g\|_{M}}{3n^3}\sum_{1\leq i<j\leq n}\underset{k\not\in\lbrace i,j\rbrace}{\sum_{1\leq k\leq n}}\mathbbm{E}\|\mathbf{Y}_n-\mathbf{Y}_n^{ijk}\|.
\label{4.2}
\end{align}
Now, by (\ref{eq_ijk}), we note that:
\begin{align}
&\|\mathbf{Y}_n-\mathbf{Y}_n^{ijk}\|\notag\\
\leq&\frac{1}{n^2}\left\lbrace(n-2)^2(|I_{ij}-I_{ij}'|+|I_{jk}-I_{jk}'|+|I_{ik}-I_{ik}'|)^2\vphantom{\left[\sum_{k\neq l}^l\right]^2}\right.\notag\\
&+\left[\sum_{l:l\neq i,j,k}\Big(|I_{ij}-I'_{ij}|(I_{jl}+I_{il})+|I_{jk}-I_{jk}'|(I_{jl}+I_{kl})+|I_{ik}-I_{ik}'|(I_{jl}+I_{il})\right.\notag\\
&\left.\left.+|I_{ik}-I_{ik}'|(I_{jl}+I_{il})\Big)+|I_{ij}I_{jk}-I_{ij}'I_{jk}'|+|I_{ij}I_{ik}-I_{ij}'I_{ik}'|+|I_{ij}I_{jk}-I_{ij}'I_{jk}'| \vphantom{\sum_{k\neq i}l}\right]^2\right\rbrace^{1/2}\notag\\
\leq&\frac{1}{n^2}\sqrt{ 9(n-2)^2+(8(n-3)+3)^2}\notag\\
=&\frac{\sqrt{73n^2-372n+477}}{n^2},\label{4.55}
\end{align}
where the second inequality follows from the fact that for all $a,b,c\in\lbrace 1,\cdots, n\rbrace$, $|I_{ab}-I_{ab}'|\leq 1$, $(I_{ab}+I_{bc})\leq 2$ and $|I_{ab}I_{bc}-I_{ab}'I_{bc}'|\leq 1$.
Therefore, by (\ref{4.2}):
\begin{align}
S_2\leq&\frac{\|g\|_Mn(n-1)(n-2)\sqrt{73n^2-372n+477}}{6n^5}\leq \frac{\sqrt{178}\|g\|_{M}}{6n}.\label{4.3}
\end{align}

Similarly, for $S_3$,
\begin{align}
S_3=&\left|\vphantom{\sum_{1\leq i<j<k\leq n}}\frac{n(n-1)}{8}\mathbbm{E}\left\lbrace D^2f(\mathbf{Y}_n)\left[(\mathbf{T}_n-\mathbf{T}_n')(2,0),(\mathbf{V}_n-\mathbf{V}_n')(0,1)\right]\right\rbrace\right.\nonumber\\
&-\mathbbm{E}D^2f(\mathbf{Y}_n)\bigg[\sum_{1\leq i<j<k\leq n}Z_{i,j,k}^{(2)}(0,1)\mathbbm{1}_{[i/n,1]\cap[j/n,1]\cap[k/n,1]},\notag\\
&\hspace{6.5cm}\left.\sum_{1\leq i< j\leq n}Z_{i,j}^{(1)}(\lfloor n\cdot \rfloor -2)(1,0)\mathbbm{1}_{[i/n,1]\cap[j/n,1]}\bigg]\right|\nonumber\\
=&\left|\frac{1}{4n^4}\sum_{1\leq i<j\leq n}\underset{k\not\in\lbrace i,j\rbrace}{\sum_{1\leq k\leq n}}\mathbbm{E}\left\lbrace \left[2\left(I_{ij}-2pI_{ij}+p\right)(I_{jk}+I_{ik})-8p^2(1-p)\right]\right.\right.\notag\\
&\left.\left.\cdot \left(D^2f(\mathbf{Y}_n)-D^2f(\mathbf{Y}^{ijk}_n)\right)\left[(0,1)\mathbbm{1}_{[i/n,1]\cap[j/n,1]\cap[k/n,1]},(\lfloor n\cdot\rfloor-2)(1,0)\mathbbm{1}_{[i/n,1]\cap[j/n,1]}\right]\right\rbrace\vphantom{\underset{j}{\sum_{1\leq i<j<k\leq n}}}\right|\nonumber\\
\stackrel{(\ref{m_bound})}\leq &\frac{\|g\|_{M}}{12n^3}\sum_{1\leq i<j\leq n}\underset{k\not\in\lbrace i,j\rbrace}{\sum_{1\leq k\leq n}}\mathbbm{E}\left\lbrace\left|2\left(I_{ij}-2pI_{ij}+p\right)(I_{jk}+I_{ik})\right|\|\mathbf{Y}_n-\mathbf{Y}_n^{ijk}\|\right\rbrace\notag\\
\leq &\frac{\|g\|_{M}}{3n^3}\sum_{1\leq i<j\leq n}\underset{k\not\in\lbrace i,j\rbrace}{\sum_{1\leq k\leq n}}\mathbbm{E}\|\mathbf{Y}_n-\mathbf{Y}_n^{ijk}\|\notag\\
\stackrel{\eqref{4.55}}\leq &\frac{\sqrt{178}\|g\|_{M}}{6n}.
\label{s_3}
\end{align}

Now, for $S_4$, let $\mathbf{Y}_n^{ijkl}$ be equal to $\mathbf{Y}_n$ except that $I_{ij},I_{ik},I_{il},I_{jk},I_{jl},I_{kl}$ are replaced with independent copies $I_{ij}',I_{ik}',I_{il}',I_{jk}',I_{jl}',I_{kl}'$, i.e. for all $t\in[0,1]$ let
\begin{align}
\mathbf{T}_n^{ijkl}(t)=&\mathbf{T}_n(t)-\frac{\lfloor nt\rfloor -2}{n^2}\left[(I_{ij}-I_{ij}')\mathbbm{1}_{[i/n,1]\cap[j/n,1]}(t)+(I_{ik}-I_{ik}')\mathbbm{1}_{[i/n,1]\cap[k/n,1]}(t)\right.\nonumber\\
&+(I_{il}-I_{il}')\mathbbm{1}_{[i/n,1]\cap[l/n,1]}(t)+(I_{jk}-I_{jk}')\mathbbm{1}_{[j/n,1]\cap[k/n,1]}(t)\nonumber\\
&\left.+(I_{jl}-I_{jl}')\mathbbm{1}_{[j/n,1]\cap[l/n,1]}(t)+(I_{kl}-I_{kl}')\mathbbm{1}_{[k/n,1]\cap[l/n,1]}(t)\right]\nonumber\\
\mathbf{V}_n^{ijkl}(t)=&\mathbf{V}_n(t)-\frac{1}{n^2}\sum_{m:m\neq i,j,k,l}\left[\left(I_{ij}-I_{ij}'\right)\left(I_{im}+I_{jm}\right)\mathbbm{1}_{[i/n,1]\cap[j/n,1]\cap[m/n,1]}(t)\right.\nonumber\\
&\hspace{-1cm}+\left(I_{ik}-I_{ik}'\right)\left(I_{im}+I_{km}\right)\mathbbm{1}_{[i/n,1]\cap[k/n,1]\cap[m/n,1]}(t)\nonumber\\
&\hspace{-1cm}+\left(I_{il}-I_{il}'\right)\left(I_{im}+I_{lm}\right)\mathbbm{1}_{[i/n,1]\cap[l/n,1]\cap[m/n,1]}(t)\nonumber\\
&\hspace{-1cm}+\left(I_{jk}-I_{jk}'\right)\left(I_{jm}+I_{km}\right)\mathbbm{1}_{[j/n,1]\cap[k/n,1]\cap[m/n,1]}(t)\nonumber\\
&\hspace{-1cm}+\left(I_{jl}-I_{jl}'\right)\left(I_{jm}+I_{lm}\right)\mathbbm{1}_{[j/n,1]\cap[l/n,1]\cap[m/n,1]}(t)\nonumber\\
&\left.\hspace{-1cm}+\left(I_{kl}-I_{ll}'\right)\left(I_{km}+I_{lm}\right)\mathbbm{1}_{[k/n,1]\cap[l/n,1]\cap[m/n,1]}(t)\right]\nonumber\\
&\hspace{-1cm}-\frac{1}{n^2}\left[(I_{ij}I_{jk}-I_{ij}'I_{jk}')+(I_{ij}I_{ik}-I_{ij}'I_{ik}')+(I_{ik}I_{jk}-I_{ij}'I_{jk}')\right]\mathbbm{1}_{[i/n,1]\cap[j/n,1]\cap[k/n,1]}(t)\nonumber\\
&\hspace{-1cm}-\frac{1}{n^2}\left[(I_{ij}I_{jl}-I_{ij}'I_{jl}')+(I_{ij}I_{il}-I_{ij}'I_{il}')+(I_{il}I_{jl}-I_{ij}'I_{jl}')\right]\mathbbm{1}_{[i/n,1]\cap[j/n,1]\cap[l/n,1]}(t)\nonumber\\
&\hspace{-1cm}-\frac{1}{n^2}\left[(I_{ik}I_{kl}-I_{ik}'I_{kl}')+(I_{ik}I_{il}-I_{ik}'I_{il}')+(I_{il}I_{kl}-I_{ik}'I_{kl}')\right]\mathbbm{1}_{[i/n,1]\cap[k/n,1]\cap[l/n,1]}(t)\nonumber\\
&\hspace{-1cm}-\frac{1}{n^2}\left[(I_{jk}I_{jl}-I_{jk}'I_{jl}')+(I_{jl}I_{kl}-I_{jl}'I_{kl}')+(I_{kl}I_{jk}-I_{kl}'I_{jk}')\right]\mathbbm{1}_{[j/n,1]\cap[k/n,1]\cap[l/n,1]}(t)
\label{eq_ijkl}
\end{align}
and for all $t\in[0,1]$ let
$\mathbf{Y}_n^{ijkl}(t)=\left(\mathbf{T}_n^{ijkl}(t)-\mathbbm{E}\mathbf{T}_n,\mathbf{V}_n^{ijkl}(t)-\mathbbm{E}\mathbf{V}_n(t)\right).$
Note that:
\begin{align}
S_4=&\Bigg|\frac{n(n-1)}{8}\mathbbm{E}\left\lbrace D^2f(\mathbf{Y}_n)\left[(\mathbf{T}_n-\mathbf{T}_n')(0,2p)+(\mathbf{V}_n-\mathbf{V}_n')(0,1),(\mathbf{V}_n-\mathbf{V}_n')(0,1)\right]\right\rbrace\nonumber\\
&-\mathbbm{E}D^2f(\mathbf{Y}_n)\bigg[\sum_{1\leq i<j<k\leq n}Z_{i,j,k}^{(2)}(0,1)\mathbbm{1}_{[i/n,1]\cap[j/n,1]\cap[k/n,1]},\notag\\
&\hspace{6.5cm}\sum_{1\leq i<j<k\leq n}Z_{i,j,k}^{(2)}(0,1)\mathbbm{1}_{[i/n,1]\cap[j/n,1]\cap[k/n,1]}\bigg]\Bigg|\nonumber\\
\leq&\Bigg|\frac{1}{4n^4}\sum_{1\leq i<j\leq n}\underset{\lbrace k,l\rbrace\cap\lbrace i,j\rbrace=\emptyset}{\sum_{1\leq k\neq l\leq n}}\mathbbm{E}\bigg\{ \Big[2p\left(I_{ij}-2pI_{ij}+p\right)(I_{jk}+I_{ik})\notag\\
&\hspace{3cm}+(I_{ij}-2pI_{ij}+p)\left(I_{ik}I_{il}+I_{ik}I_{jl}+I_{jk}I_{il}+I_{jk}I_{jl}\right)-16p^3(1-p)\big]\notag\\
&\hspace{0.5cm}\cdot \left(D^2f(\mathbf{Y}_n)-D^2f(\mathbf{Y}^{ijkl}_n)\right)\left[(0,1)\mathbbm{1}_{[i/n,1]\cap[j/n,1]\cap[l/n,1]},(0,1)\mathbbm{1}_{[i/n,1]\cap[j/n,1]\cap[k/n,1]}\right]\bigg\}\Bigg|\nonumber\\
&+\Bigg|\frac{1}{4n^4}\sum_{1\leq i<j\leq n}\underset{k\not\in\lbrace i,j\rbrace}{\sum_{1\leq k\leq n}}\mathbbm{E}\bigg\{ \Big[2p\left(I_{ij}-2pI_{ij}+p\right)(I_{jk}+I_{ik})\notag\\
&\hspace{4cm}+(I_{ij}-2pI_{ij}+p)\left(I_{ik}+2I_{ik}I_{jk}+I_{jk}\right)-4p^2(1+2p-3p^2)\Big]\notag\\
&\cdot \left(D^2f(\mathbf{Y}_n)-D^2f(\mathbf{Y}^{ijk}_n)\right)\left[(0,1)\mathbbm{1}_{[i/n,1]\cap[j/n,1]\cap[k/n,1]},(0,1)\mathbbm{1}_{[i/n,1]\cap[j/n,1]\cap[k/n,1]}\right]\Big\}\Bigg|\nonumber\\
\leq &\frac{\|g\|_{M}}{12n^4}\sum_{1\leq i<j\leq n}\underset{\lbrace k,l\rbrace\cap\lbrace i,j\rbrace=\emptyset}{\sum_{1\leq k\neq l\leq n}}\mathbbm{E}\bigg\{\Big|\left(I_{ij}-2pI_{ij}+p\right)\notag\\
&\hspace{2.5cm}\cdot(2pI_{jk}+2pI_{ik}+I_{ik}I_{il}+I_{ik}I_{jk}+I_{jk}I_{il}+I_{jk}I_{jl})\Big|\cdot\|\mathbf{Y}_n-\mathbf{Y}_n^{ijkl}\|\bigg\}\notag\\
&+\frac{\|g\|_{M}}{12n^4}\sum_{1\leq i<j\leq n}\underset{k\not\in\lbrace i,j\rbrace}{\sum_{1\leq k\leq n}}\mathbbm{E}\Big\{\left|\left(I_{ij}-2pI_{ij}+p\right)(2pI_{jk}+2pI_{ik}+I_{ik}+2I_{ik}I_{jk}+I_{jk})\right|\notag\\
&\hspace{10cm}\cdot\|\mathbf{Y}_n-\mathbf{Y}_n^{ijk}\|\Big\}
\notag\\
\leq &\frac{2\|g\|_{M}}{3n^4}\sum_{1\leq i<j\leq n}\underset{\lbrace k,l\rbrace\cap\lbrace i,j\rbrace=\emptyset}{\sum_{1\leq k\neq l\leq n}}\mathbbm{E}\|\mathbf{Y}_n-\mathbf{Y}_n^{ijkl}\|+\frac{2\|g\|_{M}}{3n^4}\sum_{1\leq i<j\leq n}\underset{k\not\in\lbrace i,j\rbrace}{\sum_{1\leq k\leq n}}\mathbbm{E}\|\mathbf{Y}_n-\mathbf{Y}_n^{ijk}\|.
\label{4.5}
\end{align}
Now, by (\ref{eq_ijkl}), note that:
\begin{align*}
&\|\mathbf{Y}_n-\mathbf{Y}_n^{ijkl}\|\\
\leq&\frac{1}{n^2}\Bigg\{\vphantom{\sum_{m:m\neq i,j,k,l}}(n-2)^2\left(|I_{ij}-I_{ij}'|+|I_{ik}-I_{ik}'|+|I_{il}-I_{i}'|+|I_{jk}-I_{jk}'|+|I_{jl}-I_{jl}'|+|I_{kl}-I_{kl}'|\right)^2\\
&+\Bigg[\sum_{m:m\neq i,j,k,l}\left[\left|I_{ij}-I_{ij}'\right|\left(I_{im}+I_{jm}\right)+\left|I_{ik}-I_{ik}'\right|\left(I_{im}+I_{km}\right)+\left|I_{il}-I_{il}'\right|\left(I_{im}+I_{lm}\right)\right.\\
&\left.\hphantom{\sum_{m:m\neq i,j,k,l}}+\left|I_{jk}-I_{jk}'\right|\left(I_{jm}+I_{km}\right)+\left|I_{jl}-I_{jl}'\right|\left(I_{jm}+I_{lm}\right)+\left|I_{kl}-I_{ll}'\right|\left(I_{km}+I_{lm}\right)\right]\\
&+|I_{ij}I_{jk}-I_{ij}'I_{jk}'|+|I_{ij}I_{ik}-I_{ij}'I_{ik}'|+|I_{ik}I_{jk}-I_{ij}'I_{jk}'|+|I_{ij}I_{jl}-I_{ij}'I_{jl}'|\\
&+|I_{ij}I_{il}-I_{ij}'I_{il}'|+|I_{il}I_{jl}-I_{ij}'I_{jl}'|+|I_{ik}I_{kl}-I_{ik}'I_{kl}'|+|I_{ik}I_{il}-I_{ik}'I_{il}'|\\
&+|I_{il}I_{kl}-I_{ik}'I_{kl}'|+|I_{jk}I_{jl}-I_{jk}'I_{jl}'|+|I_{jl}I_{kl}-I_{jl}'I_{kl}'|+|I_{kl}I_{jk}-I_{kl}'I_{jk}'|\Bigg]^2 \Bigg\}^{1/2}\\
\leq&\frac{\sqrt{36(n-2)^2+\left(12(n-4)+12\right)^2}}{n^2}\\
=&\frac{\sqrt{180n^2-1008n+1440}}{n^2}.
\end{align*}
Therefore, by (\ref{4.5}) and \eqref{4.55},
\begin{align}
S_4\leq&\frac{\|g\|_{M}\cdot \sqrt{180n^2-1008n+1440}+\sqrt{73n^2-372n+477}}{3n^2}\leq \frac{\left(\sqrt{612}+\sqrt{178}\right)\|g\|_{M}}{3n}.\label{4.6}
\end{align}

The result now follows by  (\ref{4.2}), (\ref{4.3}), (\ref{s_3}), (\ref{4.6}).
\end{proof}
\subsection{Technical details of the proof of Theorem \ref{theorem_continuous}}

\begin{lemma}\label{lemma10_app}
Using the notation of \textbf{Step 2} of the proof of Theorem \ref{theorem_continuous},
\begin{align*}
&\mathbbm{E}\left\|\mathbf{Z}_n-\mathbf{Z}\right\|\leq\frac{8}{\sqrt{n}}+\frac{39\sqrt{\log n}}{\sqrt{n}}\\
&\mathbbm{E}\left\|\mathbf{Z}_n-\mathbf{Z}\right\|^3\leq\frac{49}{n^{3/2}}+\frac{8167(\log n)^{3/2}}{n^{3/2}}\\
&\mathbbm{E}\|\mathbf{Z}\|^2\leq \frac{4}{3}.
\end{align*}
\end{lemma}
\begin{proof}
Note the following:
\begin{enumerate}
\item By Doob's $L^2$ and $L^3$ inequalities,
\begin{align}
\hspace{-0.7cm}\text{A)}\,&\mathbbm{E}\left[\sup_{t\in[0,1]}\left|\mathbf{B}_3\left(\frac{\lfloor nt\rfloor(\lfloor nt\rfloor -1)(\lfloor nt\rfloor -2)}{n^3}\right)\right|\right]\leq 2\sqrt{\mathbbm{E}\left[\left|\mathbf{B}_3\left(\frac{n(n -1)(n-2)}{n^3}\right)\right|^2\right]}\leq 2;\notag\\
\hspace{-0.7cm}\text{B)}\,&\mathbbm{E}\left[\sup_{t\in[0,1]}\left|\mathbf{B}_3\left(\frac{\lfloor nt\rfloor(\lfloor nt\rfloor -1)(\lfloor nt\rfloor -2)}{n^3}\right)\right|^3\right]\leq \frac{27}{8}\mathbbm{E}\left[\left|\mathbf{B}_3\left(\frac{n(n -1)(n-2)}{n^3}\right)\right|^3\right]\leq \frac{27}{8}.
\label{second_in}
\end{align}
\item By Doob's $L^2$ and $L^3$ inequality, for all $t\in[0,1]$,
\begin{equation}\label{fourth_in_1}
\mathbbm{E}\left[\sup_{t\in[0,1]}|\mathbf{B}_1(t^2)|\right]\leq 2,\quad \mathbbm{E}\left[\sup_{t\in[0,1]}|\mathbf{B}_1(t^2)|^3\right]\leq \frac{27}{8} \quad\text{and}\quad\left|\frac{\lfloor nt\rfloor-2}{n}-t\right|\leq \frac{3}{n}.
\end{equation}
\item Using  \cite[Lemma 3]{ito_processes} and the fact that
\[\left|\frac{\lfloor nt\rfloor(\lfloor nt\rfloor-1)}{n^2}-t^2\right|\leq \left|\frac{(nt-\lfloor nt\rfloor)(nt+\lfloor nt\rfloor)}{n^2}\right|+\frac{1}{n^2}\leq \frac{3}{n},\]
we obtain
\begin{align}
&\mathbbm{E}\left[\sup_{t\in[0,1]}\left|\mathbf{B}_1\left(\frac{\lfloor nt\rfloor(\lfloor nt\rfloor -1)}{n^2}\right)-\mathbf{B}_1(t^2)\right|\right]\leq \frac{30\sqrt{3\log\left(\frac{2n}{3}\right)}}{n^{1/2}\sqrt{\pi\log(2)}};\notag\\
&\mathbbm{E}\left[\sup_{t\in[0,1]}\left|\mathbf{B}_1\left(\frac{\lfloor nt\rfloor(\lfloor nt\rfloor -1)}{n^2}\right)-\mathbf{B}_1(t^2)\right|^3\right]\leq\frac{1080\left(3\log\left(\frac{2n}{3}\right)\right)^{3/2}}{n^{3/2}\left(\pi\log(2)\right)^{3/2}}\label{fourth_in_2}.
\end{align}
\end{enumerate}

Now, we can bound $\mathbbm{E}\left\|\mathbf{Z}_n-\mathbf{Z}\right\|$ in the following way:
\begin{align*}
&\mathbbm{E}\left\|\mathbf{Z}_n-\mathbf{Z}\right\|\nonumber\\
\leq&\frac{\sqrt{p(1-p)}}{\sqrt{2+8p^2}}\mathbbm{E}\left[\sup_{t\in[0,1]}\left|\frac{\lfloor nt\rfloor -2}{n}\mathbf{B}_1\left(\frac{\lfloor nt\rfloor(\lfloor nt\rfloor -1)}{n^2}\right)-t\mathbf{B}_1(t^2)\right|\right]\nonumber\\
&+\frac{p\sqrt{2p(1-p)}}{\sqrt{1+4p^2}}\mathbbm{E}\left[\sup_{t\in[0,1]}\left|\frac{\lfloor nt\rfloor -2}{n}\mathbf{B}_2\left(\frac{\lfloor nt\rfloor(\lfloor nt\rfloor -1)}{n^2}\right)-t\mathbf{B}_2(t^2)\right|\right]\nonumber\\
&+\frac{p\sqrt{2p(1-p)}}{\sqrt{1+4p^2}}\mathbbm{E}\left[\sup_{t\in[0,1]}\left|\frac{\lfloor nt\rfloor -2}{n}\mathbf{B}_1\left(\frac{\lfloor nt\rfloor(\lfloor nt\rfloor -1)}{n^2}\right)-t\mathbf{B}_1(t^2)\right|\right]\nonumber\\
&+\frac{2p^2\sqrt{2p(1-p)}}{\sqrt{1+4p^2}}\mathbbm{E}\left[\sup_{t\in[0,1]}\left|\frac{\lfloor nt\rfloor -2}{n}\mathbf{B}_2\left(\frac{\lfloor nt\rfloor(\lfloor nt\rfloor -1)}{n^2}\right)-t\mathbf{B}_2(t^2)\right|\right]\nonumber\\
&+\frac{p(1-p)}{n^{1/2}}\mathbbm{E}\left[\sup_{t\in[0,1]}\left|\mathbf{B}_3\left(\frac{\lfloor nt\rfloor(\lfloor nt\rfloor -1)(\lfloor nt\rfloor -2)}{n^3}\right)\right|\right]\nonumber\\
\stackrel{(\ref{second_in})}\leq& \frac{(1+4p+4p^2)\sqrt{p(1-p)}}{\sqrt{2+8p^2}}\mathbbm{E}\left[\sup_{t\in[0,1]}\left|\frac{\lfloor nt\rfloor -2}{n}\mathbf{B}_1\left(\frac{\lfloor nt\rfloor(\lfloor nt\rfloor -1)}{n^2}\right)-t\mathbf{B}_1(t^2)\right|\right]\nonumber\\
&+\frac{2p(1-p)}{n^{1/2}}\nonumber\\
\leq&\frac{(1+4p+4p^2)\sqrt{p(1-p)}}{\sqrt{2+8p^2}}\left(\mathbbm{E}\left[\sup_{t\in[0,1]}\left|\left(\frac{\lfloor nt\rfloor -2}{n}-t\right)\mathbf{B}_1(t^2)\right|\right]\right.\nonumber\\
&\left.+\mathbbm{E}\left[\sup_{t\in[0,1]}\left|\mathbf{B}_1\left(\frac{\lfloor nt\rfloor(\lfloor nt\rfloor -1)}{n^2}\right)-\mathbf{B}_1(t^2)\right|\right]\right)+\frac{2p(1-p)}{n^{1/2}}\nonumber\\
\stackrel{(\ref{fourth_in_1}),(\ref{fourth_in_2})}\leq&\frac{(1+4p+4p^2)\sqrt{p(1-p)}}{\sqrt{2+8p^2}}\left(\frac{6}{n}+\frac{30\sqrt{3\log n}}{n^{1/2}\sqrt{\pi\log(2)}}\right)+\frac{2p(1-p)}{n^{1/2}}\nonumber\\
\leq&\frac{8}{\sqrt{n}}+\frac{39\sqrt{\log n}}{\sqrt{n}}.
\end{align*}

Similarly,
\begin{align*}
&\mathbbm{E}\|\mathbf{Z}_n-\mathbf{Z}\|^3\nonumber\\
\leq&4\sqrt{3}\left(\frac{p(1-p)}{2+8p^2}\right)^{3/2}\mathbbm{E}\left[\sup_{t\in[0,1]}\left|\frac{\lfloor nt\rfloor -2}{n}\mathbf{B}_1\left(\frac{\lfloor nt\rfloor(\lfloor nt\rfloor -1)}{n^2}\right)-t\mathbf{B}_1(t^2)\right|^3\right]\nonumber\\
&+4\sqrt{3}\left(\frac{2p^3(1-p)}{1+4p^2}\right)^{3/2}\mathbbm{E}\left[\sup_{t\in[0,1]}\left|\frac{\lfloor nt\rfloor -2}{n}\mathbf{B}_2\left(\frac{\lfloor nt\rfloor(\lfloor nt\rfloor -1)}{n^2}\right)-t\mathbf{B}_2(t^2)\right|^3\right]\nonumber\\
&+9\sqrt{3}\left(\frac{2p^3(1-p)}{1+4p^2}\right)^{3/2}\mathbbm{E}\left[\sup_{t\in[0,1]}\left|\frac{\lfloor nt\rfloor -2}{n}\mathbf{B}_1\left(\frac{\lfloor nt\rfloor(\lfloor nt\rfloor -1)}{n^2}\right)-t\mathbf{B}_1(t^2)\right|^3\right]\nonumber\\
&+9\sqrt{3}\left(\frac{8p^5(1-p)}{1+4p^2}\right)^{3/2}\mathbbm{E}\left[\sup_{t\in[0,1]}\left|\frac{\lfloor nt\rfloor -2}{n}\mathbf{B}_2\left(\frac{\lfloor nt\rfloor(\lfloor nt\rfloor -1)}{n^2}\right)-t\mathbf{B}_2(t^2)\right|^3\right]\nonumber\\
&+9\sqrt{3}\frac{p^3(1-p)^3}{n^{3/2}}\mathbbm{E}\left[\sup_{t\in[0,1]}\left|\mathbf{B}_3\left(\frac{\lfloor nt\rfloor(\lfloor nt\rfloor -1)(\lfloor nt\rfloor -2)}{n^2}\right)\right|^3\right]\nonumber\\
\stackrel{(\ref{second_in})}\leq& \frac{\sqrt{6}p^{3/2}(1-p)^{3/2}(1+26p^3+126p^6)}{(1+4p^2)^{3/2}}\notag\\
&\hspace{1.5cm}\cdot\mathbbm{E}\left[\sup_{t\in[0,1]}\left|\frac{\lfloor nt\rfloor -2}{n}\mathbf{B}_1\left(\frac{\lfloor nt\rfloor(\lfloor nt\rfloor -1)}{n^2}\right)-t\mathbf{B}_1(t^2)\right|^3\right]+\frac{243\sqrt{3}}{512n^{3/2}}\nonumber\\
\leq& \frac{4\sqrt{6}p^{3/2}(1-p)^{3/2}(1+26p^3+126p^6)}{(1+4p^2)^{3/2}}\left(\mathbbm{E}\left[\sup_{t\in[0,1]}\left|\left(\frac{\lfloor nt\rfloor -2}{n}-t\right)\mathbf{B}_1(t^2)\right|^3\right]\right.\nonumber\\
&\left.+\mathbbm{E}\left[\sup_{t\in[0,1]}\left|\mathbf{B}_1\left(\frac{\lfloor nt\rfloor(\lfloor nt\rfloor -1)}{n^2}\right)-\mathbf{B}_1(t^2)\right|^3\right]\right)+\frac{243\sqrt{3}}{512n^{3/2}}\nonumber\\
\stackrel{(\ref{fourth_in_1}),(\ref{fourth_in_2})}\leq&\frac{4\sqrt{6}p^{3/2}(1-p)^{3/2}(1+26p^3+126p^6)}{(1+4p^2)^{3/2}}\left(\frac{81}{8n^3}+\frac{1080\left(3\log n\right)^{3/2}}{n^{3/2}\left(\pi\log(2)\right)^{3/2}}\right)+\frac{243\sqrt{3}}{512n^{3/2}}\nonumber\\
\leq&\frac{49}{n^{3/2}}+\frac{8167(\log n)^{3/2}}{n^{3/2}}.
\end{align*}
Furthermore,
\begin{align*}
\mathbbm{E}\|\mathbf{Z}\|^3\leq&\sqrt{2}\mathbbm{E}\left[\sup_{t\in[0,1]}\left(\frac{\sqrt{p(1-p)}}{\sqrt{2+8p^2}}t\mathbf{B}_1(t^2)+\frac{p\sqrt{2p(1-p)}}{\sqrt{1+4p^2}}t\mathbf{B}_2(t^2)\right)^3\right]\nonumber\\
&+\sqrt{2}\mathbbm{E}\left[\sup_{t\in[0,1]}\left(\frac{p\sqrt{2p(1-p)}}{\sqrt{1+4p^2}}t\mathbf{B}_1(t^2)+\frac{2p^2\sqrt{2p(1-p)}}{\sqrt{1+4p^2}}t\mathbf{B}_2(t^2)\right)^3\right]\nonumber\\
\leq&\frac{2p^{3/2}(1-p)^{3/2}(1+2^{7/2}p^3+2^{11/2}p^6)}{(1+4p^2)^{3/2}}\mathbbm{E}\left[\sup_{t\in[0,1]}|\mathbf{B}_1(t^2)|^3\right]\nonumber\\
\leq&\frac{27p^{3/2}(1-p)^{3/2}(1+2^{7/2}p^3+2^{11/2}p^6)}{4(1+4p^2)^{3/2}}\leq \frac{4}{3}.
\end{align*}
This finishes the proof.
\end{proof}




\begin{thebibliography}{99}

\bibitem{barbour_mollison}
A.~Barbour and D.~Mollison, \emph{Epidemics and random graphs}, Stochastic
  Processes in Epidemic Theory (Berlin, Heidelberg) (J.-P. Gabriel,
  C.~Lef{\`e}vre, and Ph. Picard, eds.), Springer Berlin Heidelberg, 1990,
  pp.~86--89.

\bibitem{barbour}
A.D. Barbour, \emph{{Stein's method and Poisson process convergence}}, Journal
  of Applied Probability \textbf{25} (1988), 175--184.

\bibitem{diffusion}
\bysame, \emph{{Stein's Method for Diffusion Approximation}}, Probability
  Theory and Related Fields \textbf{84} (1990), 297--322.

\bibitem{janson}
A.D. Barbour, Lars Holst, and Svante Janson, \emph{Poisson approximation},
  Oxford Studies in Probability, Clarendon Press, 1992.

\bibitem{functional_combinatorial}
A.D. Barbour and S.~Janson, \emph{{A functional combinatorial central limit
  theorem}}, Electronic Journal of Probability \textbf{14} (2009), no.~81,
  2352--2370.

\bibitem{Basa}
A.~Basalykas, \emph{Functional limit theorems for random multilinear forms},
  Stochastic Process. Appl. \textbf{53} (1994), no.~1, 175--191. \MR{1290712}

\bibitem{decreusefond2}
E.~Besan{\c c}on, L.~Decreusefond, and P.~Moyal, \emph{{Stein's method for
  diffusive limit of Markov processes}}, arXiv:1805.01691, 2018.

\bibitem{Bill}
P.~Billingsley, \emph{Convergence of probability measures}, John Wiley \& Sons,
  Inc., New York-London-Sydney, 1968. \MR{0233396}

\bibitem{Blom}
G.~Blom, \emph{Some properties of incomplete {$U$}-statistics}, Biometrika
  \textbf{63} (1976), no.~3, 573--580. \MR{474582}

\bibitem{campese}
S.~Bourguin and S.~Campese, \emph{{Approximation of Hilbert-valued Gaussian
  measures on Dirichlet structures}}, Electron. J. Probab. \textbf{25} (2020), no.~150, 1--30.

\bibitem{BrKil}
B.~M. Brown and D.~G. Kildea, \emph{Reduced {$U$}-statistics and the
  {H}odges-{L}ehmann estimator}, Ann. Statist. \textbf{6} (1978), no.~4,
  828--835. \MR{491556}

\bibitem{cannings}
C.~Cannings, \emph{{Evolutionary Stable Strategies}}, {Encycl. Maths.
  Supplemen.}, vol.~1, Kluwer Academic Publishers, 1997.

\bibitem{models}
C.~Cannings and D.B. Penman, \emph{Ch. 2. models of random graphs and their
  applications}, Stochastic Processes: Modelling and Simulation, Handbook of
  Statistics, vol.~21, Elsevier, 2003, pp.~51 -- 91.

\bibitem{chatterjee}
S.~Chatterjee, P.~Diaconis, and E.~Meckes, \emph{Exchangeable pairs and poisson
  approximation}, Probab. Surveys \textbf{2} (2005), 64--106.

\bibitem{chatterjee1}
S.~Chatterjee, J.~Fulman, and A.~Röllin, \emph{Exponential approximation by
  stein’s method and spectral graph theory}, ALEA Lat. Am. J. Probab. Math.
  Stat (2011).

\bibitem{meckes}
S.~Chatterjee and E.~Meckes, \emph{Multivariate normal approximation using
  exchangeable pairs}, ALEA Lat. Am. J. Probab. Math. Stat. \textbf{4} (2008),
  257--283. \MR{2453473}

\bibitem{normal_approx}
L.H.Y Chen, L.~Goldstein, and Q.-M. Shao, \emph{Normal approximation by
  stein’s method}, {Probability and Its Applications}, Springer Verlag, 2011.

\bibitem{christofides}
T.C. Christofides, \emph{{Maximal probability inequalities for
  multidimensionally indexed semimartingales and convergence theory of
  u-statistics}}, Ph.D. thesis, Johns Hopkins University, 1987.

\bibitem{Coutin}
L.~Coutin and L.~Decreusefond, \emph{{Stein's method for Brownian
  Approximations}}, Communications on Stochastic Analysis \textbf{7} (2013),
  no.~3, 349--372.

\bibitem{decreusefond_higher}
\bysame, \emph{{Higher order expansions via Stein's method}}, Communications on
  Stochastic Analysis \textbf{8} (2014), no.~2, 155--168.

\bibitem{decreusefond_rough}
\bysame, \emph{{Stein's method for rough paths}}, Potential Analysis \textbf{53} (2020), 387-–406

\bibitem{CH88}
M.~Cs\"org\"o and L.~Horvath, \emph{Invariance principles for changepoint
  problems}, J. Multivariate Anal. \textbf{27} (1988), 151--168.

\bibitem{CH_book}
\bysame, \emph{Limit theorems in changepoint analysis}, Wiley, 1997.

\bibitem{deJo87}
P.~de~Jong, \emph{A central limit theorem for generalized quadratic forms},
  Probab. Theory Related Fields \textbf{75} (1987), no.~2, 261--277.
  \MR{885466}

\bibitem{Doe20}
C.~D\"obler, \emph{{Normal approximation via non-linear exchangeable pairs}}, arXiv:2008.02272, 2020.


\bibitem{dobler}
C.~D{\"o}bler, \emph{Stein's method of exchangeable pairs for the beta
  distribution and generalizations}, Electron. J. Probab. \textbf{20} (2015),
  34 pp.

\bibitem{DKP}
C.~D\"obler, M.~Kasprzak, and G.~Peccati, \emph{{Functional convergence of
  $U$-processes with size-dependent kernels}}, {\tt arXiv:1912.02705} (2019).

\bibitem{DP16}
C.~D\"obler and G.~Peccati, \emph{{Quantiative de Jong theorems in any
  dimension}}, Electron. J. Probab. \textbf{22} (2017), no. 2, 1--35.

\bibitem{DP18b}
C.~D\"{o}bler and G.~Peccati, \emph{The gamma {S}tein equation and noncentral
  de {J}ong theorems}, Bernoulli \textbf{24} (2018), no.~4B, 3384--3421.
  \MR{3788176}

\bibitem{DP19}
\bysame, \emph{Quantitative {CLT}s for symmetric {$U$}-statistics using
  contractions}, Electron. J. Probab. \textbf{24} (2019), Paper No. 5, 43.
  \MR{3916325}

\bibitem{Ferger94}
D.~Ferger, \emph{{An extension of the Cs{\"o}rg{\H o}-Horv{\'a}th functional
  limit theorem and Its applications to changepoint problems}}, J. Multivariate
  Anal. \textbf{51} (1994), no.~2, 338--351.

\bibitem{Ferger01}
\bysame, \emph{Analysis of change-point estimators under the null hypothesis},
  Bernoulli \textbf{7} (2001), no.~3, 487--506.

\bibitem{ito_processes}
M.~Fischer and G.~Nappo, \emph{{On the Moments of the Modulus of Continuity of
  Ito Processes}}, Stochastic Analysis and Applications \textbf{28} (2010),
  no.~1, 103--122.

\bibitem{Gombay2004}
E.~Gombay, \emph{{U-Statistics in Sequential Tests and Change Detection}},
  Sequential Analysis: Design Methods and Applications \textbf{23} (2004),
  no.~2, 257--274.

\bibitem{GH95}
E.~Gombay and L.~Horvath, \emph{{An application of U-statistics to change-point
  analysis}}, Acta Sci. Math. (Szeged) \textbf{60} (1995), 345--357.

\bibitem{gotze}
F.~G{\"o}tze, \emph{{On the rate of convergence in the multivariate CLT}}, The
  Annals of Probability \textbf{19} (1991), no.~2, 724--739.

\bibitem{GT}
F.~G\"{o}tze and A.~N. Tikhomirov, \emph{Asymptotic distribution of quadratic
  forms}, Ann. Probab. \textbf{27} (1999), no.~2, 1072--1098. \MR{1699003}

\bibitem{HR_survey}
L.~Horvath and G.~Rice, \emph{Extensions of some classical methods in change
  point analysis}, TEST \textbf{23} (2014), no.~2, 219--255.

\bibitem{Ja}
S.~Janson, \emph{The asymptotic distributions of incomplete {$U$}-statistics},
  Z. Wahrsch. Verw. Gebiete \textbf{66} (1984), no.~4, 495--505. \MR{753810}

\bibitem{Janson1991}
S.~Janson and K.~Nowicki, \emph{{The asymptotic distributions of generalized
  U-statistics with applications to random graphs}}, Probability Theory and
  Related Fields \textbf{90} (1991), no.~3, 341--375.

\bibitem{kasprzak1}
M.J. Kasprzak, \emph{{Diffusion approximations via Stein's method and time
  changes}}, arXiv:1701.07633, 2017.

\bibitem{kasprzak3}
M.J. Kasprzak, \emph{{Functional approximations via Stein's method of
  exchangeable pairs}}, Ann. Inst. H. Poincaré Probab. Statist. \textbf{56} (2020), no.~4, 2540--2564.

\bibitem{kasprzak2}
M.J. Kasprzak, \emph{{Stein's method for multivariate Brownian approximations
  of sums under dependence}}, {Stochastic Processes and their Applications} \textbf{130}
  (2020), no.~8, 4927--4967.
  
\bibitem{kasprzak}
M.J. Kasprzak, A.~B. Duncan, and S.J. Vollmer, \emph{{Note on A. Barbour's
  paper on Stein's method for diffusion approximations}}, Electron. Commun.
  Probab. \textbf{22} (2017), no.~23, 1--8.

\bibitem{reinert}
C.~Ley, G.~Reinert, and Y.~Swan, \emph{{Stein's method for comparison of
  univariate distributions}}, Probability Surveys \textbf{14} (2017), 1--52.

\bibitem{Major}
P.~Major, \emph{Asymptotic distributions for weighted {$U$}-statistics}, Ann.
  Probab. \textbf{22} (1994), no.~3, 1514--1535. \MR{1303652}

\bibitem{meckes09}
E.~Meckes, \emph{On stein's method for multivariate normal approximation},
  Collections, vol. Volume 5, pp.~153--178, Institute of Mathematical
  Statistics, Beachwood, Ohio, USA, 2009.

\bibitem{Mik}
T.~Mikosch, \emph{The rate of convergence in the functional central limit
  theorem for random quadratic forms with some applications to the law of the
  iterated logarithm}, Monatsh. Math. \textbf{107} (1989), no.~2, 137--153.
  \MR{994980}

\bibitem{nourdin}
I.~Nourdin and G.~Peccati, \emph{Normal approximations with malliavin
  calculus}, Cambridge tracts in Mathematics, Cambridge University Press, 2012.

\bibitem{NPR}
I.~Nourdin, G.~Peccati, and G.~Reinert, \emph{Invariance principles for
  homogeneous sums: universality of {G}aussian {W}iener chaos}, Ann. Probab.
  \textbf{38} (2010), no.~5, 1947--1985. \MR{2722791}

\bibitem{OnRe}
K.~A. O'Neil and R.~A. Redner, \emph{Asymptotic distributions of weighted
  {$U$}-statistics of degree {$2$}}, Ann. Probab. \textbf{21} (1993), no.~2,
  1159--1169. \MR{1217584}

\bibitem{RW19}
A.~Ra{\v c}kauskas and M.~Wendler, \emph{{Convergence of U-Processes in
  H\"older Spaces with Application to Robust Detection of a Changed Segment}}, Statistical Papers \text{61} (2020), 1409-–1435. 

\bibitem{reinert_roellin}
G.~Reinert and A.~R{\"o}llin, \emph{{Multivariate normal approximation with
  Stein's method of exchangeable pairs under a general linearity condition}},
  The Annals of Probability \textbf{37} (2009), no.~6, 2150--2173.

\bibitem{reinert_roellin1}
\bysame, \emph{Random subgraph counts and $u$-statistics: Multivariate normal
  approximation via exchangeable pairs and embedding}, Journal of Applied
  Probability \textbf{47} (2010), no.~2, 378--393.

\bibitem{RiUt}
M.~Rifi and F.~Utzet, \emph{On the asymptotic behavior of weighted
  {$U$}-statistics}, J. Theoret. Probab. \textbf{13} (2000), no.~1, 141--167.
  \MR{1744988}

\bibitem{RiRo97}
Y.~Rinott and V.~Rotar, \emph{On coupling constructions and rates in the {CLT}
  for dependent summands with applications to the antivoter model and weighted
  {$U$}-statistics}, Ann. Appl. Probab. \textbf{7} (1997), no.~4, 1080--1105.
  \MR{1484798 (99g:60050)}

\bibitem{roellin1}
A.~R{\"o}llin, \emph{Translated poisson approximation using exchangeable pair
  couplings}, The Annals of Applied Probability \textbf{17} (2007), no.~5/6,
  1596--1614.

\bibitem{ross}
N.~Ross, \emph{{Fundamentals of Stein's Method}}, Probability Surveys
  \textbf{8} (2011), 210--293.

\bibitem{Rot1}
V.~I. Rotar, \emph{Certain limit theorems for polynomials of degree two}, Teor.
  Verojatnost. i Primenen. \textbf{18} (1973), 527--534. \MR{0326803}

\bibitem{rubin_vitale}
H.~Rubin and R.A. Vitale, \emph{{Asymptotic Distribution of Symmetric
  Statistics}}, Ann. Statist. \textbf{8} (1980), no.~1, 165--170.

\bibitem{Ruc}
A.~Ruci\'{n}ski, \emph{When are small subgraphs of a random graph normally
  distributed?}, Probab. Theory Related Fields \textbf{78} (1988), no.~1,
  1--10. \MR{940863}

\bibitem{serfling}
R.J. Serfling, \emph{{Approximation Theorems of Mathematical Statistics}},
  Wiley Series in Probability and Statistics, John Wiley and Sons, Inc., 1980.

\bibitem{shih}
H.-H. Shih, \emph{{On Stein's method for infinite-dimensional Gaussian
  approximation in abstract Wiener spaces}}, Journal of Functional Analysis
  \textbf{261} (2011), no.~5, 1236 -- 1283.

\bibitem{stein}
Ch. Stein, \emph{{A bound for the error in the normal approximation to the
  distribution of a sum of dependent random variables}}, Proc. Sixth Berkeley
  Symp. on Math. Statist. and Prob. \textbf{2} (1972), 583--602.

\bibitem{stein1}
\bysame, \emph{{Approximate Computation of Expectations}}, Institute of
  Mathematical Statistics Lecture Notes, Monograph Series, 7. Hayward, Calif.,
  Institute of Mathematical Statistics, 1986.

\bibitem{swan}
Y.~Swan, \emph{{A gateway to Stein's Method}},
  \url{https://sites.google.com/site/steinsmethod/home}, 2016, Accessed on
  19/05/2016.

\bibitem{vitale}
A.R. Vitale, \emph{An expansion for symmetric statistics and the efron-stein
  inequality}, Lecture Notes--Monograph Series, vol.~5, pp.~112--114, Institute
  of Mathematical Statistics, 1984.

\end{thebibliography}



\ACKNO{The authors would like to thank Gesine Reinert, Giovanni Peccati and Alison Etheridge for helpful discussions and comments on the early versions of this work. Mikołaj Kasprzak was supported by the \textbf{FNR grant FoRGES (R-AGR- 3376-10)} at Luxembourg University.}


\end{document}